\documentclass[twoside,a4paper,reqno,11pt]{amsart}
\usepackage{amsfonts, amsbsy, amsmath, amsthm, amssymb, latexsym, verbatim, enumerate}
\usepackage{mathrsfs}
\usepackage[top=30mm,right=30mm,bottom=30mm,left=30mm]{geometry}

\usepackage{hyperref}
\usepackage{amssymb,amsfonts,amsmath,amsopn,amstext,amscd,latexsym,amsthm}
\usepackage{bbm}
\usepackage[all,cmtip]{xy}

\usepackage{enumerate}
\usepackage[shortlabels]{enumitem}

\theoremstyle{plain}

\begin{document}

\def\a{\alpha}
\def\om{\omega}
 \def\b{\beta}
 \def\e{\epsilon}
 \def\d{\delta}
  \def\D{\Delta}
 \def\c{\chi}
 \def\k{\kappa}
 \def\g{\gamma}
 \def\t{\tau}
\def\ti{\tilde}
 \def\N{\mathbb N}
 \def\Q{\mathbb Q}
 \def\Z{\mathbb Z}
 \def\C{\mathbb C}
 \def\F{\mathbb F}
 \def\ovF{\overline\F}
 \def\bfN{\mathbf N}
 \def\cG{\mathcal G}
 \def\cT{\mathcal T}
 \def\cX{\mathcal X}
 \def\cY{\mathcal Y}
 \def\cC{\mathcal C}
 \def\cD{\mathcal D}
 \def\cZ{\mathcal Z}
 \def\cO{\mathcal O}
 \def\cW{\mathcal W}
 \def\cL{\mathcal L}
 \def\bfC{\mathbf C}
 \def\bfZ{\mathbf Z}
 \def\bfO{\mathbf O}
 \def\G{\Gamma}
 \def\bO{\boldsymbol{\Omega}}
 \def\bgo{\boldsymbol{\omega}}
 \def\go{\rightarrow}
 \def\do{\downarrow}
 \def\ra{\rangle}
 \def\la{\langle}
 \def\fix{{\rm fix}}
 \def\ind{{\rm ind}}
 \def\rfix{{\rm rfix}}
 \def\diam{{\rm diam}}
 \def\uni{{\rm uni}}
 \def\diag{{\rm diag}}
 \def\Irr{{\rm Irr}}
 \def\Syl{{\rm Syl}}
 \def\Gal{{\rm Gal}}
 \def\Tr{{\rm Tr}}
 \def\M{{\cal M}}
 \def\cE{{\mathcal E}}
\def\td{\tilde\delta}
\def\tx{\tilde\xi}
\def\DC{D^\circ}

\def\Ker{{\rm Ker}}
 \def\rank{{\rm rank}}
 \def\soc{{\rm soc}}
 \def\Cl{{\rm Cl}}
 \def\A{{\sf A}}
 \def\sP{{\sf P}}
 \def\sQ{{\sf Q}}
 \def\SSS{{\sf S}}
  \def\SQ{{\SSS^2}}
 \def\St{{\sf {St}}}
 \def\p{\ell}
 \def\ps{\ell^*}
 \def\SC{{\rm sc}}
 \def\supp{{\sf{supp}}}
  \def\cR{{\mathcal R}}
 \newcommand{\tw}[1]{{}^#1}

 \def\Sym{{\rm Sym}}
 \def\PSL{{\rm PSL}}
 \def\SL{{\rm SL}}
 \def\Sp{{\rm Sp}}
 \def\GL{{\rm GL}}
 \def\SU{{\rm SU}}
 \def\GU{{\rm GU}}
 \def\SO{{\rm SO}}
 \def\PO{{\rm P}\Omega}
 \def\Spin{{\rm Spin}}
 \def\PSp{{\rm PSp}}
 \def\PSU{{\rm PSU}}
 \def\PGL{{\rm PGL}}
 \def\PGU{{\rm PGU}}
 \def\Iso{{\rm Iso}}
 \def\Stab{{\rm Stab}}
 \def\GO{{\rm GO}}
 \def\Ext{{\rm Ext}}
 \def\E{{\cal E}}
 \def\l{\lambda}
 \def\ve{\varepsilon}
 \def\Lie{\rm Lie}
 \def\s{\sigma}
 \def\O{\Omega}
 \def\o{\omega}
 \def\ot{\otimes}
 \def\op{\oplus}
 \def\oc{\overline{\chi}}
 \def\pf{\noindent {\bf Proof.$\;$ }}
 \def\Proof{{\it Proof. }$\;\;$}
 \def\no{\noindent}
\def\hal{\unskip\nobreak\hfil\penalty50\hskip10pt\hbox{}\nobreak
 \hfill\vrule height 5pt width 6pt depth 1pt\par\vskip 2mm}

 \renewcommand{\thefootnote}{}

\newtheorem{theorem}{Theorem}
 \newtheorem{thm}{Theorem}[section]
 \newtheorem{prop}[thm]{Proposition}
 \newtheorem{lem}[thm]{Lemma}
 \newtheorem{lemma}[thm]{Lemma}
 \newtheorem{defn}[thm]{Definition}
 \newtheorem{cor}[thm]{Corollary}
 \newtheorem{coroll}[theorem]{Corollary}
\newtheorem*{corB}{Corollary}
 \newtheorem{rem}[thm]{Remark}
 \newtheorem{exa}[thm]{Example}
 \newtheorem{cla}[thm]{Claim}

\numberwithin{equation}{section}
\parskip 1mm

\title{Multiplicity-free representations of algebraic groups II}

\author[M.W. Liebeck]{Martin W. Liebeck}
\address{M.W. Liebeck, Department of Mathematics,
    Imperial College, London SW7 2BZ, UK}
\email{m.liebeck@imperial.ac.uk}

\author[G.M. Seitz]{ Gary M. Seitz}
\address{G.M. Seitz, University of Oregon, Eugene, Oregon  97403, USA} 
\email{seitz@uoregon.edu}

\author[D.M. Testerman]{Donna M. Testerman}
\address{D.M. Testerman, EPFL, Lausanne,  CH-105  Switzerland}
\email{donna.testerman@epfl.ch}
\date{}

\maketitle

\medskip
\begin{center} {\it Dedicated to the memory of our friend and colleague Jan Saxl}
\end{center}
\medskip

\begin{abstract}
We continue our work, started in \cite{MF}, on the program of classifying triples $(X,Y,V)$, where $X,Y$ are simple algebraic groups over an algebraically closed field of characteristic zero with $X<Y$, and $V$ is an irreducible module for $Y$ such that the restriction $V\downarrow X$ is multiplicity-free. In this paper we handle the case where $X$ is of type $A$, and is irreducibly embedded in $Y$ of type $B,C$ or $D$. It turns out that there are relatively few triples for $X$ of arbitrary rank, but a number of interesting exceptional examples arise for small ranks.
\end{abstract}


 \footnotetext{


The authors would like to thank the referee for reading the paper very carefully and suggesting many improvements. 

The authors acknowledge the support of the National Science Foundation under Grant No. DMS-1440140 while they were in residence at the Mathematical Sciences Research Institute in Berkeley, California, during the Spring 2018 semester, and also support from the Cecilia Tanner Fund of the Mathematics Department at Imperial College London. The third author also 
acknowledges support from the  Swiss Science Foundation grants 200021-156583 and 200021-175571. }



\section{Introduction}

Let $K$ be an algebraically closed field of characteristic zero. A finite-dimensional module for a connected reductive algebraic group over $K$ is said to be {\it multiplicity-free} (abbreviated as MF) if each of its composition factors appears with multiplicity 1. Many classical results in the literature can be framed as part of the study of triples $(X,Y,V)$ satisfying the following properties:
\begin{itemize}
\item[(1)] $Y$ is a connected reductive group over $K$, and $V$ is a finite-dimensional irreducible $Y$-module such that 
the action of $Y$ on $V$ does not contain a full classical group $SL(V)$, $Sp(V)$ or $SO(V)$;
\item[(2)] $X$ is a proper connected reductive subgroup of $Y$;
\item[(3)] the restriction $V \downarrow X$ is multiplicity-free. 
\end{itemize}
For example, as explained in the Introduction to \cite{MF}, there are celebrated results of Dynkin \cite{dynk}, Howe \cite{Howe}, Kac \cite{Kac},  Stembridge \cite{stem} and Weyl \cite{weyl} that fall under this paradigm. In \cite{MF}, we initiated the program of classifying the triples $(X,Y,V)$ as above that satisfy the further conditions that $X$ and $Y$ are simple, and $X$ is irreducibly embedded in $Y$ (that is, $X$ is contained in no parabolic subgroup of $Y$). The main result of \cite{MF} determines such triples in the case where $X$ and $Y$ are both of type $A$, and \cite{LST} solves the problem in general in the case where $X=A_1$.

In this paper we continue this program, and deal with the case where $X$ has type $A$ and $Y$ is a classical group of type $B$, $C$ or $D$. Our main result classifies all triples $(X,Y,V)$ in this case when $X$ acts irreducibly on the natural module for $Y$. In the statement, we use the usual notation for irreducible modules: $\l_i$ and $\o_i$ are fundamental dominant weights for $Y$ and $X$, respectively. Also for a dominant weight $\l$, $V_Y(\l)$ denotes the irreducible $KY$-module of highest weight $\l$.


\begin{theorem}\label{mainthm} Let $X = A_{l+1}$ with $l\ge 0$, and let $W = V_X(\d)$ be an irreducible self-dual module, such that $X$ embeds via $\d$ as a proper subgroup of $Y = SO(W)$ or $Sp(W)$. Suppose $V = V_Y(\l)$ is an irreducible $Y$-module such that $V\downarrow X$ is multiplicity-free, and assume $\l \ne 0,\,\l_1$. Then $\l,\d$ are as in Table $\ref{tab1}$ or $\ref{tbl2}$.
Conversely, for each possibility in the tables, $V\downarrow X$ is multiplicity-free.
\end{theorem}


\begin{table}[h!]
\caption{MF examples for $X$ of arbitrary rank $\ge 2$}\label{tab1}
\[
\begin{array}{|l|l|}
\hline
\l & \d \\
\hline
\l_2, 2\l_1 & b\om_{s+1}\,(l \hbox{ even}, s=\frac{l}{2}) \\
                 & \om_1+\om_{l+1} \\
                 & \om_s+\om_{s+1}\,(l \hbox{ odd}, s = \frac{l+1}{2}) \\
\hline
\l_2 & 2\om_1+2\om_{l+1} \\
       & \om_2+\om_l\,(l\ge 3) \\
\hline
\l_3, 3\l_1 &  \om_1+\om_{l+1} \\
\hline
\end{array}
\]
\end{table}

\begin{table}[h!]
\caption{Exceptional MF examples for $X$ of small rank}\label{tbl2}
\[
\begin{array}{|l|l|l|l|}
\hline
X & \d & Y &  \l   \\
\hline
A_1      & 3 & C_2 & 01,\,b0,\,0b\,(2\le b\le 5),\,11,\,12,\,21 \\
     & 4 & B_2 & 01,\,b0,\,0b\,(2\le b\le 5),\,11,\,12,\,21 \\
     & 5 & C_3 & 010,\,200,\,001,\,300 \\
     & 6 & B_3 & 010,\,200,\,001,\,101,\,002,\,300 \\

 & 2n-1\,(n\ge 4) & C_n & \l_2,\,2\l_1,\,\l_3\,(n=4,5),\,\l_n\,(n=4,5) \\
         & 2n & B_n & \l_2,\,2\l_1,\,\l_n\,(n\le 8) \\
 \hline
A_2 & 11 & D_4 & c\l_i\,(c\le 5,\,i=1,3,4) \\
        &&& \l_2 \\
        &&& \l_i+\l_j\,(i,j \in \{1,3,4\}) \\
\hline
A_3 & 020 & D_{10} & c\l_1\,(c\le 3) \\
            &&& \l_2,\l_3 \\
            &&& \l_9,\l_{10} \\
\hline
A_5 & \om_3 &  C_{10} & \l_i \\
             &&& c\l_1 \,(c \le 5) \\
            &&& 2\l_2 \\
         &&& \l_1+\l_i\, (i = 2,3,4)   \\
        &&& \l_1+\l_{10} \\
        &&& 2\l_1 + \l_2   \\
\hline
A_7 & \om_4 &  D_{35}& \l_3,\,\l_4 \\
              &&&  3\l_1,\,4\l_1 \\
              &&&  \l_{34}, \l_{35}\\
\hline 
A_9 & \o_5 &  C_n\,,n = \frac{1}{2}{10\choose 5} & \l_3 \\
       &&& 3\l_1 \\
\hline 
A_{11} & \o_6 & D_n\,,n = \frac{1}{2}{12\choose6} & \l_3 \\
          &&& 3\l_1 \\
\hline
 A_{13} & \om_7 & C_n\,,n = \frac{1}{2}{14\choose7}&\l_3\\      
\hline
\end{array}
\]
\end{table}


Note that the hypothesis that $X$ is proper in $Y$ in the theorem is necessary, in order to avoid the case where $X=A_3$ and $\d = \o_2$, in which case $Y=D_3=X$. Note also that the entries in Table \ref{tbl2} with $X=A_1$ are obtained from \cite{LST}, so this case will not be considered here.

It is interesting to note that when the rank of $X$ is large (greater than 13), there are relatively few triples (those in Table \ref{tab1}), and these are all consequences of examples in \cite{MF}. 
For small ranks, many of the exceptional triples in Table \ref{tbl2} also follow from, or are closely related to, examples in \cite{MF}, with two notable exceptions where $(X,Y) = (A_3, D_{10})$ or $(A_7,D_{35})$  and $V_Y(\l)$ is a spin module. Particularly striking is the latter case, where $V_Y(\l)\downarrow X$ is MF with 36 composition factors (see Lemma \ref{MFAD}).

Our proof of Theorem \ref{mainthm} follows a similar strategy to that in \cite{MF}. We choose a parabolic subgroup 
$P_X = Q_XL_X$ of $X$ (with Levi factor $L_X' = A_l$ and unipotent radical $Q_X$), and a parabolic subgroup $P_Y = Q_YL_Y$ of $Y$, such that $Q_X\le Q_Y$ and $L_X\le L_Y$. The proof proceeds by analysing the action of $L_X'$ on the {\it levels} $V^i(Q_Y)$, where $V^1(Q_Y) = V/[V,Q_Y]$, $V^2(Q_Y) = [V,Q_Y]/[V,Q_Y,Q_Y]$, and so on (see Section \ref{levs}). 

After some preliminaries in Sections 2-4,  the theorem is proved for a series of specific weights $\d$ in Sections 5-11. Based on this, the general proof for arbitrary weights $\d$ is then carried out in Sections 12 and 13.

Throughout the paper, for various simple complex Lie algebras $L$ of small rank, we make substantial use of the Lie theoretic representation theory packages in Magma \cite{Magma} that enable one to decompose tensor products and  symmetric and alternating powers of irreducible representations of $L$ with given highest weights.

\section{General Setup}\label{setupp}

Before describing our specific notation for $X$ and $Y$, we give a general definition of the $S$-{\it value} and the $L$-{\it value} of a module for a semisimple algebraic group $G$. Let $\mu_1,\ldots ,\mu_n$ be fundamental dominant weights for $G$, and for a dominant weight $\mu = \sum_{i=1}^nc_i\mu_i$, define $S(\mu) = \sum c_i$ and also let $L(\mu)$ be the number of values of $i$ such that $c_i\ne 0$. For a $G$-module $Z$, set $S(Z)$ to be the maximum of $S(\mu)$ over all irreducible summands $V_G(\mu)$ of $Z$. Also, we shall sometimes denote $\mu$ by the sequence of integers $c_1c_2\cdots$.

Now let $X < Y$ be simple algebraic groups over an algebraically closed field $K$ of characteristic $0$, with $X = A_{l+1}$ and $Y$ classical.  Let $W$ be the natural module for $Y$ and assume that $X$ acts irreducibly on $W$ via a representation of highest weight $\d$. 
Let $\Sigma(X)$ denote the root system of $X$ with respect to a 
maximal torus $T_X$, and let $\Pi(X)$ be a fundamental system in $\Sigma(X)$. 
Write $\Pi(X) = \{\a_1,\ldots,\a_{l+1}\}$, and let $\omega_1, \ldots, \omega_{l+1}$ be the corresponding fundamental dominant weights of $X$. For a dominant weight $\d = \sum_1^{l+1}c_i\o_i$, denote by $V_X(\d)$ the irreducible $KX$-module of highest weight $\d$.

Next, let $P_X = Q_XL_X$ be a maximal parabolic subgroup of $X$ corresponding to the end-node $\a_{l+1}$, with unipotent radical $Q_X$ and Levi factor $L_X$. Write $L_X = L_X'T$, where $T$
is the $1$-dimensional central torus. We take $Q_X$ to be a product of root subgroups for negative roots.
We have $\Pi(L_X') = \{\alpha_1, \ldots, \alpha_l\}$ and  $\Pi(X) = \Pi(L_X') \cup \{\alpha\}$, 
where  $\alpha= \alpha_{l+1}$. 
Write $T_X = S_XT$, where $S_X = T_X \cap L_X'$. 


Now let $T_Y$ be a fixed maximal torus of $Y,$ let $\Sigma(Y)$ be the corresponding root system, $\Pi(Y) = \{\b_1, \ldots, \b_n\}$ a fundamental system of positive roots, and  $ \l_1, \ldots, \l_n $  the corresponding set of fundamental dominant weights for $Y$.  Let $P_Y = Q_YL_Y$ be a parabolic subgroup of $Y$ with  Levi factor $L_Y$ containing $T_Y$, and 
unipotent radical $Q_Y$  a product of root groups for negative roots. Write $L_Y' = C^0 \times \cdots \times C^k$ with each $C^i$ simple. For each $i$,  $\Pi(C^i) =  \la \b_1^i, \ldots, \b_{r_i}^i \ra$ is a string of  fundamental roots of 
$\Pi(Y)$ and we  order the factors so that the string for $\Pi(C^i)$ comes before the string for $\Pi(C^{i+1})$. 
For each $j$ with $1 \le j \le r_i$,  let $\l_j^i$ denote the fundamental dominant weight corresponding to $\b_j^i$.
 
In situations to follow, we will choose a specific parabolic $P_Y$ associated to our embedding $X<Y$ and the choice of $P_X$. It then turns out that there will be just one 
fundamental root separating  $C^{i-1}$ and  $C^i$ which is designated by $\gamma_i$. Also when $C^k$ has type $A$, the  last root $\b_n$ lies outside $\Pi(C^k)$, and will be written as $\g_{k+1}$. These roots $\g_i$ comprise all of the fundamental roots lying outside $\Pi(L_Y')$.

Let $\l$ be a dominant weight for $Y$ and set $V = V_Y(\l)$.  Write $\mu^i = \l \downarrow (T_Y \cap C^i)$, 
so that $V_{C^0}(\mu^0) \otimes \cdots \otimes V_{C^k}(\mu^k)$ is a composition factor of $L_Y$ on $V_Y(\l)$.

We shall denote the Lie algebra of $X$ by $L(X)$.
For a positive root $\a$ in $\Sigma(X)$, let $e_\a$ and $f_\a = e_{-\a}$ be corresponding root elements in $L(X)$, and for $1\le i\le l+1$,  let $f_i = f_{\alpha_i + \cdots + \alpha_{l+1}}$. Then $\{f_1, \ldots, f_{l+1} \}$ is a basis of  commuting root elements  of the Lie algebra $L(Q_X)$.

Finally, a few general pieces of terminology: for a semisimple algebraic group $G$ and a dominant weight $\l$, we shall often denote the module $V_G(\l)$ just by the weight $\l$ -- and if there is any danger of confusion, we shall write $\l \oplus \mu$ for the sum $V_G(\l)+V_G(\mu)$ (rather than just $\l+\mu$, which could refer to the module $V_G(\l+\mu)$); 
if we have modules $A = B\oplus C$ we shall write $C = A-B$; and finally, throughout the rest of the paper we shall abbreviate the term ``multiplicity-free" by the initials MF.

\section{ Levels}\label{levs}

In this section we introduce the notion of levels of a module with respect to a parabolic subgroup, and 
establish some basic results on these.  For these
results we allow $X$ and $Y$ to be classical groups of any type.  Later we will restrict to the case where $X$ has
type $A$ and $Y$ has type $B,C,$ or $D$.

We shall make use of a few basic results in \cite{mem1}: specifically (1.2), (2.3) and (3.6). Although these results are proved there under the assumption of positive characteristic, their proofs make no use of this assumption, and thus the results also hold in our situation of characteristic zero.

Set $[V,Q_Y^0]=V$, define $[V, Q_Y^1]$ to be the commutator space $[V,Q_Y]$, and for  $d>1$  inductively define  $[V, Q_Y^d] = [[V,Q_Y^{d-1}],Q_Y]$.  Now set 
\[
V^{d+1}(Q_Y) = [V,Q_Y^d]/[V,Q_Y^{d+1}].
\]
 This quotient is $L_Y$-invariant and will be called the $d$th level.  Thus $V^1(Q_Y) = V/[V, Q_Y]$ is level $0$; 
$V^2(Q_Y) = [V,Q_Y]/[V,Q_Y^2]$ is level $1$, and so on. It is well-known (see \cite[(1.2)]{mem1} for example) that $V^1(Q_Y)$ is an irreducible $L_Y'$-module.

  We have $V^1(Q_Y) = V/[V,Q_Y]$ and $V^2(Q_Y) = [V,Q_Y]/[V,Q_Y^2].$  It follows from \cite[(2.3)(ii)]{mem1} that  $V^2(Q_Y) = [V,Q_Y]/[V,Q_Y^2]$ can be regarded as the direct sum of weight spaces 
  of $V$ corresponding to weights of the form $\l -\psi -\gamma_j$, where $\psi$ is a sum of positive roots in $\Sigma (L_Y')$ and
 $\g_j$ is as defined above.
  Therefore we can write $V^2(Q_Y) = \sum_{j=1}^{k+1} V_{\g_j}^2,$ where $V_{\gamma_j}^2 = V_{\gamma_j}^2(Q_Y)$ is the sum of such weight spaces for a fixed value of $j$. This is an $L_Y'$-module, and we set 
$S_j^2 = S(V_{\g_j}^2(Q_Y))$ (the $S$-value, as defined in Section 2).

\begin{lem}\label{parabembed}  There is a parabolic subgroup $P_Y = Q_YL_Y$  containing $P_X$ such that the following conditions hold:
\begin{itemize}
\item[{\rm (i)}] $Q_X < Q_Y$;
\item[ {\rm (ii)}] $L_X \le L_Y = C_Y(T)$ and $T_X \le T_Y$, a maximal torus of $L_Y$;
 \item[{\rm (iii)}] $\eta \downarrow T = \alpha \downarrow T$ for each fundamental root $\eta \in  \Pi(Y) \setminus \Pi(L_Y')$;
 \item[{\rm (iv)}] $[W, Q_X^d] = [W,Q_Y^d]$ for $d  \ge 0$.
\end{itemize}
\end{lem}

\pf  Parts (i)-(iii) follow from (3.6) of \cite{mem1}, and (iv) follows from the proof of that result. \hal

Recall from the previous section that $W$ is the natural module for $Y$, and $W = V_X(\d)$ with $\d = \sum_1^{l+1}c_i\o_i$. We have $L_Y' = C^0\times \cdots \times C^k$, and for each $i$, the embedding of $L_X'$ in $C^i$ corresponds to the action of $L_X'$ on the $i$th level $W^i(Q_X)$ (defined as above). The composition factors of these levels are given by \cite[Thm. 5.1.1]{MF}. In particular,  the first level $W^1(Q_X)$ is irreducible for $L_X'$, with highest weight 
\[
\d' = \sum_1^lc_i\o_i.
\]

For the rest of this section  assume that $V = V_Y(\l)$ is such that $V\downarrow X = V_1+\cdots +V_s$ is MF, and $V_i = V_X(\theta_i)$ for dominant weights $\theta_i$. As above, for each $i$ and $d \ge 0$ we write  $V_i^{d+1}(Q_X) = [V_i,Q_X^d]/[V_i,Q_X^{d+1}]$.
Also, $V_i^1(Q_X) = V_i/[V_i,Q_X]$ is an irreducible $L_X'$-module for each $i$.

Throughout we will use the following notation
\[
V^i = V^i(Q_Y)\downarrow L_X'.
\]
In particular, $V^1 = \left(V/[V,Q_Y]\right) \downarrow L_X'$.

\begin{lem} \label{Twts} {\rm (i)} For $d \ge 0$, $T$ induces scalars on $V^{d+1}(Q_Y)$  via the weight $(\lambda - d\alpha) \downarrow T$.

{\rm (ii)}  For $1\le i \le s$ and $t \ge 0$, $T$ induces scalars on $V_i^{t+1}(Q_X)$  via the weight $(\theta _i - t\alpha) \downarrow T$.

\end{lem} 

\pf  Here we refer to (2.3)(ii) of \cite{mem1}.  Applying that result to the action of $Y$ on $V$  shows that
$V^{d+1}(Q_Y)$ is isomorphic as a vector space to the direct sum of weight spaces of $V$ having
$Q_Y$-level $d$.   As $T$ centralizes $L_Y$, part (i) follows from Lemma \ref{parabembed}(iii). Part (ii) is similar. \hal

\begin{lem} \label{locateVi}  Let $1 \le i \le s.$
\begin{itemize}
\item[{\rm (i)}] There is a unique $n_i \ge 0$ such that $\theta_i \downarrow T =  (\lambda - n_i \alpha) \downarrow T.$
\item[{\rm (ii)}] $n_i$ is maximal subject to $V_i \le [V,Q_Y^{n_i}].$
\item[{\rm (iii)}] $\left(V_i + [V,Q_Y^{n_i + 1}]\right)/ [V,Q_Y^{n_i + 1}]$ is irreducible under the action of $L_X'$.
\end{itemize}
\end{lem} 

\pf    Lemma \ref{Twts}(ii) shows that  the weights of $T$
on $V_X(\theta _i)$ have the form $\theta _i - t\alpha$ for $0\le t \le l_i$ where $l_i$ is
maximal among values $j$ with  $[V_i, Q_X^j] \ne 0.$  On the other hand, from Lemma \ref{Twts}(i) we see that all weights
of $T$ on $V$ have form $(\lambda - s\alpha)\downarrow T$ for some $s\ge 0$.  Therefore there is a unique value, say $n_i$,
such that $\theta_i \downarrow T =  (\lambda - n_i \alpha) \downarrow T.$  This gives (i) and (ii).

For (iii), first note that $(V_i + [V,Q_Y^{n_i + 1}])/ [V,Q_Y^{n_i + 1}] \cong V_i/(V_i \cap  [V,Q_Y^{n_i + 1}] )$
as $L_X'$-modules. Now  $Q_X \le Q_Y$ implies that  $[V_i, Q_X] \le V_i \cap [V,Q_Y^{n_i + 1}]$, 
which is properly contained in $V_i$ by the maximality of $n_i$. 
 The result follows since we know that $V_i/[V_i,Q_X]$ is an irreducible $L_X'$-module. \hal

\begin{lem}\label{resttoLX'}  Assume $i,j \in \{1, \ldots, s\}$ with $i\ne j$ satisfy $n_i = n_j = m$. Then  
\begin{itemize}
\item[{\rm (i)}] $\theta_i  \downarrow S_X \ne \theta _j \downarrow S_X$;
\item[{\rm (ii)}] $(V_i + [V,Q_Y^{m + 1}])/ [V,Q_Y^{m + 1}] \not  \cong (V_j + [V,Q_Y^{m + 1}])/ [V,Q_Y^{m + 1}]$
as $L_X'$-modules.
\end{itemize}
 \end{lem}
 
 \pf  Lemma \ref{locateVi} implies that both $(V_i + [V,Q_Y^{m + 1}])/ [V,Q_Y^{m + 1}]$ and  
$(V_j + [V,Q_Y^{m + 1}])/ [V,Q_Y^{m + 1}]$ are irreducible
$L_X'$-summands of $V^{m+1}(Q_Y)$. The highest weights of these summands 
are $\theta_i \downarrow S_X$ and $\theta_j \downarrow S_X$, respectively.  Now $T_X = S_XT$
and by hypothesis and Lemma \ref{locateVi}, $\theta_i \downarrow T = \theta_j \downarrow T$.  
As $\theta_i \ne \theta_j$ we conclude that they have different restrictions to $S_X$ and the
assertions follow. \hal

\begin{prop}\label{induct}  The following assertions hold.
\begin{itemize}
\item[{\rm (i)}] $V^1 =  \sum _{i: n_i = 0} V_i/[V_i,Q_X]$ is MF.
\item[{\rm (ii)}] $V^2 = \sum _{i: n_i = 0} V_i^2(Q_X) + \sum _{i: n_i = 1} V_i/[V_i,Q_X]$.  Moreover, the second sum is MF.
\item[{\rm (iii)}] For any $d\ge 0$,
\[
V^{d+1}  = \sum _{i\,:\,0\le n_i\le d-1} V_i^{d+1-n_i}(Q_X) + \sum _{i\,:\, n_i = d} V_i/[V_i,Q_X],
\]
and the last summand is MF.
\end{itemize}
\end{prop}

\pf  Note that (i) and (ii) are the special cases $d=0,1$ of (iii), so it suffices to prove (iii). 
Fix $d \ge 0$. As mentioned in the proof of Lemma \ref{Twts}, $V^{d+1}(Q_Y) = [V,Q_Y^d]/[V,Q_Y^{d+1}]$ is isomorphic as a vector space to the direct sum of weight spaces of $V$ having
$Q_Y$-level $d.$ By Lemma \ref{Twts}(iii) these weights restrict to $T$ as does $\lambda - d\alpha.$  On the other hand,
$V \downarrow X = \sum V_i$ and for each $i$ the weights of $T$ on $V_i$ have the form $\theta _i - t\alpha$ for  
non-negative integers $t$. Now we use a weight space comparison and Lemma \ref{locateVi} to complete the proof:  indeed, in order to get  weight $\lambda - d\alpha$ we can start with $\theta_j$ with $n_j \le d$
and take the $T$-weight space for weight $\theta_j - (d-n_j)\alpha$ of $V_j.$  Lemma \ref{Twts} shows that this weight space is $L_X$-isomorphic to $V_j^{d-n_j}(Q_X)$.  Furthermore, the sum of such weight spaces
yields  $V^{d+1}(Q_Y)$.  This together with Lemma \ref{resttoLX'} gives (iii).   \hal

If $Z = V_i^j(Q_X)$ is a summand of $V^c$ as in Proposition \ref{induct}(iii), then $Z'=V_i^{j+1}(Q_X)$ is a summand of 
$V^{c+1}$, and $S(Z')\le S(Z)+1$ (see \cite[Lemma 3.7]{MF}). We say that $Z'$ {\it arises from} $Z$, and note that given 
$Z$, the possible highest weights of  composition factors of $Z'$ can be computed using Corollary 5.1.2 of \cite{MF}.

\begin{cor}\label{conseq}  Let $\rho$ be a dominant weight for the maximal torus $S_X$  of $L_X'$.
\begin{itemize}
\item[{\rm (i)}]  If $V_{L_X'}(\rho)$ appears with
multiplicity $t$ in $\sum _{i: n_i = 0} V_i^2(Q_X),$  then $V_{L_X'}(\rho)$ appears with multiplicity
at most $t+1$ in $V^2$.  
\item[{\rm (ii)}] If $V_{L_X'}(\rho)$ appears with
multiplicity $t$ in $\sum _{i: n_i = 0} V_i^3(Q_X) +\sum _{j: n_j= 1} V_j^2(Q_X),$  then $V_{L_X'}(\rho)$ appears with multiplicity at most $t+1$ in $V^3$.
\end{itemize}
\end{cor}

\pf  This follows from Proposition \ref{induct}(ii),(iii). \hal

\begin{lem}\label{prop38} Let $d\ge 1$, and suppose $V_{L_X'}(\nu)$ is a composition factor of multiplicity at least $2$ in $V^{d+1}(Q_Y)\downarrow L_X'$. Then $S(\nu) \le S(V^{d}) + 1$.
\end{lem}

\pf This follows from  \cite[Prop. 3.8]{MF}.  \hal

\section{Preliminary lemmas}\label{prel}

As in the hypothesis of Theorem \ref{mainthm}, let $X = A_{l+1}$ with $l\ge 0$, and let $W = V_X(\d)$ be a self-dual module, such that $X$ embeds via $\d$ in $Y = SO(W)$ or $Sp(W)$. Let $\d = \sum_1^{l+1}d_i\o_i$, and write this also as $\d = (d_1,\ldots,d_{l+1})$. Take $P_X = Q_XL_X$ and $P_Y=Q_YL_Y$ to be parabolic subgroups of $X$ and $Y$ as in Section \ref{setupp}, with $L_Y' = C^0\times \cdots \times C^k$. 

\begin{lem}\label{setup} Let $X < Y$ be as above.
\begin{itemize}
\item[{\rm (i)}] Suppose $l+1=2s$ and $\d = (a_1, \cdots ,a_s,a_s, \cdots ,a_1)$. Then
$Y$ is an orthogonal group, $k = \sum_{i=1}^s a_i $ and $V_{C^k}(\l_1^k) \downarrow L_X' \supseteq
(a_1, \cdots ,a_{s-1},2a_s,a_{s-1}, \cdots ,a_1).$ 

\item[{\rm (ii)}] Suppose $l+1=2s+1$ and $\d = (a_1, \cdots ,a_s,c,a_s, \cdots ,a_1)$. Then one of the following holds:
\begin{itemize}
\item[{\rm (a)}]  $c$ is even, $Y$ is orthogonal, $k = a_1 + \cdots +a_s + \frac{c}{2}$ and $V_{C^k}(\l_1^k) \downarrow L_X' \supseteq (a_1, \cdots ,a_{s-1},a_s +\frac{c}{2},\,a_s +\frac{c}{2},\,a_{s-1}, \cdots ,a_1).$ 
\item[{\rm (b)}] $c$ is odd, $Y$ is symplectic or orthogonal according to whether $s$ is even or odd, $k = a_1 + \cdots +a_s + \frac{c-1}{2}$ and $V_{C^k}(\l_1^k) \downarrow L_X' \supseteq (a_1, \cdots ,a_{s-1},a_s +\frac{c-1}{2},\,a_s +\frac{c+1}{2},\,a_{s-1}, \cdots ,a_1).$
\end{itemize}
\end{itemize}
In each case in {\rm (i), (ii)},  the indicated composition factor has maximal $S$-value among composition factors of 
$V_{C^k}(\l_1^k) \downarrow L_X'$.
\end{lem}

\pf   The actions of $L_X'$ on the levels of $W$ are given by \cite[Thm. 5.1.1]{MF}, and the expressions for $k$ follow from this.
 To determine whether $Y$ is symplectic or orthogonal we use a result of Steinberg \cite[Lemma 79]{St}.  Let $z = \prod_{\a}h_{\a}(-1),$ the product over positive roots, viewed as an element of $SL_{l+2}$. Then  $z$ is in the center of $SL_{l+2},$ $z^2 = 1$, and $z$ acts as $1$ or $-1$ on the natural module according to whether $l$ is odd or even. Steinberg's lemma  shows that $\d$ affords an  orthogonal or symplectic representation according to whether $z$ induces $1$ or $-1$ on $W$. In the situation of (i), $l$ is odd, so $z$ induces $1$ on the natural module and on $W$ so that $Y$ is an orthogonal group.  In the situation of (ii), $l$ is even so that $z $ induces $-1$ on the natural module and hence $(-1)^i$ on the $i$th wedge of the natural module.  It follows that $z$ induces $(-1)^{c(s+1)}$ on $W$  and the nature of $Y$ (symplectic or orthogonal)  asserted in (ii) follows.  

To simplify notation, recall the definition of $f_i$ in Section \ref{setupp} and  set $\bar{f_i} = f_{l+2-i}$. Then 
a maximal vector of largest $S$-value 
at level $k$ is achieved by applying $\bar {f_s}^{a_s}\cdots \bar {f_1}^{a_1}$, $\bar {f}_{s+1}^{\frac{c}{2}}\bar {f_s}^{a_s}\cdots \bar {f_1}^{a_1}$, or 
$\bar {f}_{s+1}^{\frac{c-1}{2}}\bar {f_s}^{a_s}\cdots \bar {f_1}^{a_1}$ to a maximal vector (of weight $\d$), according to whether we are in the situation of (i), (ii)(a), or (ii)(b), respectively. Corollary 5.1.4 of  \cite{MF} shows this is of maximal $S$-value.
\hal

We next state a useful result of Koike and Terada \cite{tensor1} (see p.509).

\begin{lem}\label{tensorl1} 
{\rm (i)} Let $Y = C_n\,(n \ge 2)$, and let $\l$ be a dominant weight for $Y$.  Then
$V_Y(\l_1) \otimes V_Y(\l) = \sum V_Y(\mu)$, where the sum runs over dominant weights $\mu$ with Young diagrams obtained by adding or deleting one square from that of $\l.$ In particular, 
\[
V_Y(\l_1) \otimes V_Y(\l_i) = V_Y(\l_1 + \l_i) +V_Y(\l_{i-1}) +V_Y(\l_{i+1}),
\]
where the last term on the right hand side is omitted if $i = n.$

{\rm (ii)} Let $Y = B_n$ or $D_n$ . Then for $i<n-1$,
\[
V_Y(\l_1) \otimes V_Y(\l_i) = V_Y(\l_1 + \l_i) +V_Y(\l_{i-1}) +V_Y(\l_{i+1}),
\]
except that for $Y=D_n$ and $i=n-2$, the last term should be replaced by $V_Y(\l_{n-1}+\l_n)$. 
\end{lem}


We shall also need the following result of Stembridge \cite{stem}.

\begin{lem}\label{stemb}  Assume that $D$ is a simple algebraic group over $K$, and $\mu,\nu$ are dominant weights such that $V_D(\mu)\otimes V_D(\nu)$ is MF.  Then either $\mu$ or $\nu$ is a multiple of a fundamental dominant weight.
\end{lem} 

The next result is a well-known fact about symmetric and wedge powers for symplectic and orthogonal groups.

\begin{lem}\label{symwed} Let $Y=C_n$, $B_n$ or $D_n$ be a symplectic or orthogonal group with natural module $W$. Let $a\ge 2$.
\begin{itemize}
\item[{\rm (i)}] If $Y=C_n$, then $V_Y(a\l_1) = S^a(W)$, and for $a\le n$ we have $V_Y(\l_a) = \wedge^a(W)-\wedge^{a-2}(W)$.
\item[{\rm (ii)}] If $Y=B_n$ or $D_n$, then $V_Y(a\l_1) = S^a(W)-S^{a-2}(W)$; and for $a<n-1$ we have  $V_Y(\l_a) = \wedge^a(W)$; also for $Y=B_n$ we have $V_Y(\l_{n-1}) = \wedge^{n-1}(W)$, and for $Y=D_n$ we have $V_Y(\lambda_{n-1}+\lambda_n) = \wedge^{n-1}(W)$.
\end{itemize}
\end{lem}

Finally, we shall need some special facts about certain symmetric and exterior powers of modules for groups of type $A$. Recall that $X = A_{l+1}$; for convenience, write $m=l+1$.

\begin{lem}\label{extab} Let $X=A_m$ with $m\ge 2$, and let $\nu = (\cdots,a,b,\cdots)$ be a dominant weight for $X$ with $a,b$ in the $j^{th}$ and $(j+1)^{st}$ positions, where $a\ge 1, b\ge 2$. Then both $\wedge^3\nu$ and $(\wedge^2\nu \otimes \nu) - \wedge^3\nu$ have a repeated summand of highest weight $\e = 3\nu - \a_j-2\a_{j+1} = (\cdots,3a,3b-3,\cdots)$, and $S(\e) \ge S(\nu)+2a+2b-3$.
\end{lem}

\pf Any weight affording an irreducible summand of $\wedge^3 \nu$ must have the form $3\nu-\rho$, where 
$\rho=\sum_{i=1}^{m}c_i\alpha_i$ with $\sum_{i=1}^{m}c_i \geq 2$. Any such summand containing the weight $\epsilon$ is either $3\nu-2\alpha_{j+1}$, $3\nu-\alpha_j-\alpha_{j+1}$ or $\epsilon$. Now it is straightforward to see that the first weight does not afford a summand of $\wedge^3 \nu$, while the second affords a summand of multiplicity 1. Moreover, the weight $\epsilon$ occurs with multiplicity 1 in this summand. On the other hand $\epsilon$ occurs with multiplicity 3 in $\wedge^3 \nu$ as $\nu\wedge(\nu-\alpha_j)\wedge(\nu-2\alpha_{j+1})$ and $\nu\wedge(\nu-\alpha_{j+1})\wedge
(\nu-\alpha_j-\alpha_{j+1})$. (Here we abuse notation and write the wedge of weights to mean the wedge of the associated weight spaces.) This establishes the claim for $\wedge^3 \nu$.

To treat the second case, we need to show that $\epsilon$ affords a summand of multiplicity at least 4 in the tensor product 
$\wedge^2 \nu\otimes \nu$. To this end, we note that $\wedge^2 \nu$ has summands with highest weights $2\nu-\alpha_j$, $2\nu-\alpha_{j+1}$, $2\nu-\alpha_j-\alpha_{j+1}$ and $2\nu-\alpha_j-2\alpha_{j+1}$. (We use a weight count similar to the above.)
Tensoring each of these summands with $\nu$ we find that the first gives rise to one summand $\epsilon$, the second (using \cite[7.1.5]{MF}) gives rise to two summands $\epsilon$ and the last tensor product gives rise to one summand $\epsilon$, hence there are at least four such summands, as required.
 \hal

\begin{lem}\label{sympowers} Let $X = A_{m}$ with $m\ge 2$.
\begin{itemize}
\item[{\rm (i)}] If $m=2$, then $S^2(11) = (22)+(11)+0$ and $\wedge^2(11) = (30) + (03) +(11).$

\item[{\rm (ii)}] If $m \ge 3$, then $S^2(10\cdots01) = (20\cdots02)+(10\cdots01)+(010\cdots010)+0.$

\item[ {\rm (iii)}] If $m\ge 3$, then $\wedge^2(10\cdots01) = (010\cdots 02)+(20\cdots010)+(10\cdots01).$
 
\item[{\rm (iv)}]  If $m\ge 4$, then 
\[
\begin{array}{ll}
\wedge^3(10\cdots01)= & (110\cdots011)+(20\cdots02)+(10\cdots 01)+ (010\cdots 010)+ \\
   & (30\cdots 0100)+(0010\cdots 03)+(20\cdots010)+(010\cdots02)+0.
\end{array}
\]
If $m=3$ then $\wedge^3(101) = (101)+(020)+(121)+(202)+(210)+(012)+(400)+(004)+(000)$. And if $m=2$, then 
$\wedge^3(11) = (11)+(22)+(30)+(03)+(00)$.
  
\item[{\rm (v)}] If $m\ge 3$ and $a \ge 4$, then $S^a(10\cdots01) \supseteq ((a-2)0\cdots 0(a-2))^3.$

\item[{\rm (vi)}] $S^3(11) = (33)+(22)+(11) +(30)+(03)+(00).$

\item[{\rm (vii)}] $S^3(101) = (303)+(202)+(101)^2+(121)+(012)+(210)+(000).$

\item[{\rm (viii)}] $S^3(1001) = (3003)+(2002)+(1001)^2+(1111)+(0102)+(2010)+(0110)+(0000).$

\item[{\rm (ix)}] If $m \ge 5$, then 
\[
\begin{array}{ll}
S^3(10\cdots01)= & (30\cdots03)+(20\cdots02)+(10\cdots 01)^2+ (010\cdots 010)+ \\
   & (0010\cdots0100)+(110\cdots011)+(20\cdots010)+(010\cdots02)+0,
\end{array}
\]
where if $m=5$, the fifth term of the right hand side is replaced by $(00200)$.
\end{itemize}
 \end{lem}

\pf Parts (i), (vi), (vii), and (viii) are immediate from Magma computations, while (ii) and (iii) can be easily established using the domino technique of \cite{domi} (see \cite[\S 4.2]{MF}). Also (iv) follows from \cite[6.3.5]{MF}, the proof of which gives the full decomposition of $\wedge^3(10\cdots01)$.

Now consider parts (v) and (ix). 
View $X < Z = SL(W)$ and set $V = S^b(W)$, so $\l = b\l_1$.
We consider levels as usual. Then $V^1 = (b0\cdots0)$ and all other irreducible summands of $V \downarrow X$ involve highest weights at lower levels. We will determine these summands using level analysis combined with the method where  the highest weights of irreducible modules for $L_X'$ at certain levels can be used to find highest weights for irreducible summands of  $V \downarrow X$. 

For this we make use of the torus $T$ in $X$ consisting of elements $T(c) = h_1(c)h_2(c^2)\cdots h_{m}(c^{m})$.    Each of the levels $V^i(Q_Z)$ is a weight space  for the maximal torus $T_Z$, and $T(c)$ commutes with $L_X'$ and   induces scalars via the restriction of the weight $\l-(i-1)\a_{m}.$ 

  For (v) and (ix), $\l = b\l_1 $ restricts to $(b0\cdots0b)$ for $T_X$, so $T(c)$ acts on $V^1$ as $c^{b+bm} = c^{(m+1)b}$.  Also $T(c)$  acts via $c^{2m-(m-1)} = c^{m+1}$ on the root vector for the root $\a_{m}.$ Therefore $T(c)$ induces 
$c^{(m+1)b -(m+1)(i-1)}= c^{(m+1)(b-i+1)}$ on $V^i.$
  
 We now focus on part (v). Assume $V = S^a(\d)$ for $a \ge 4.$  Of particular interest is the multiplicity of $\mu = ((a-2)0\cdots0)$ in $V^3.$
Here $V^1 = (a0\cdots0)$ and $V^2 = ((a-1)0\cdots0) \otimes ( (10\cdots 01 ) + 0).$ Therefore
$V^2 = (a0\cdots 01)+ ((a-2)10\cdots01) + ((a-1)0\cdots0)^2 + ((a-3)10\cdots0),$ and at most 3 summands of $V^3$ of highest weight $\mu$ can arise from this.
(Here, and many times subsequently, we compute the summands of $V^3$ ``arising" from $V^2$ as explained in the remark after Proposition \ref{induct}.) 

 Now $V^3$ has  summands $((a-2)0\cdots 0) \otimes (S^2(10\cdots 01) + 0)$ and $((a-1)0\cdots 0) \otimes (0\cdots01).$ From (ii) we see that the first summand is $((a-2)0\cdots 0) \otimes ((20\cdots02)+(10\cdots01)^2+(010\cdots010)+0^2)$. Tensoring 
 $((a-2)0\cdots 0)$ with the four summands of the second tensor factor we see that $\mu$ appears $1,2,0,2$ times, respectively. And $\mu$ appears once in the second summand of $V^3.$ Therefore $\mu^6$ appears in $V^3.$ It follows from the above that   $V\downarrow X$ has a summand of the form
$((a-2)0\cdots0x)^3.$  

 Now  $T(c)$ induces $c^{(m+1)(a-3+1)} = c^{(m+1)(a-2)}$ on $V^3.$ On the other hand this must equal 
$c^{(a-2)+mx}$ so that $(m+1)(a-2) = a-2+mx$ which forces $x = a-2,$ as required for (v).
 
%

It remains to prove (ix). To avoid special cases we use Magma to establish the result for $m \le 6$ so now assume $m > 6.$ Consideration of dominant weights shows that the possible weights of composition factors of $S^3(\om_1+\om_m)$ are among $(c0\cdots 0c)$ ($c \le 3$), $(010\cdots 010),$ $(0010\cdots 0100)$, $(110\cdots 011)$, $(20 \cdots 010)$, and $(010\cdots 02).$  Other dominant weights are easily ruled out as highest weights of composition factors.

 We have $V^1 = (30\cdots0)$ and  $V^2 = (20\cdots0) \otimes ( (10\cdots 01 ) + 0) = (30\cdots01) + (110\cdots01) + (20\cdots 0)^2+ (010\cdots0). $ Only $(20\cdots 0)^1$ and $(30\cdots 01)$ can arise from $V^1$, so it follows that there are composition factors of $V\downarrow X$ with highest weights $(20\cdots0x)$, $(010\cdots0x)$, and $(110\cdots01x)$.  These each occur with multiplicity 1 in $V \downarrow X$, and using $T(c)$ we show that $x = 2,2,1$, respectively. For example in the first case we use the above to get  the equality $c^{(m+1)(3-2+1)} = c^{2+mx}$ which forces $x = 2.$ Similarly for the other cases. 

Now $V^3$ is the sum of $(10\cdots0) \otimes S^2((10\cdots01) + 0)$ and $(20\cdots0) \otimes (0\cdots 01).$ Expanding these with the help of (ii) we find that there must exist composition factors of highest weights $(010\cdots01x)$, $(10\cdots0x)$, $(20\cdots 01x)$ and $(0010\cdots010x)$ appearing with multiplicities $1,2,1,1$, respectively.  We then show that $x =0,1,0, 0$ respectively.
At this point we have dealt with the multiplicity of all the weights above
with the exception of the trivial module.  

Here we consider $V^4 $ which contains summands 
$0\otimes S^3(10\cdots01 + 0)$  and $(10\cdots0) \otimes (10\cdots01)+0)\otimes (0\cdots01),$ where the second summand is afforded by $3\l_1 -2 (\b_1^0+\cdots+\b_{r_0}^0)-2\g_1-(\b_1^1+
\cdots + \b_{r_1}^1)-\g_2.$ Using induction along with
(i) and (ii) we expand the terms and find that $0^5$ occurs.  On the other
hand  the composition factors obtained so far contribute precisely 4 copies of the trivial module for $L_X'$.  Therefore there is a composition factor for $X$ having highest weight $(0\cdots0x)$ and appearing with multiplicity 1. The usual argument shows that $x = 0$, completing the proof of (ix).  \hal

\begin{lem}\label{extbd} Let $X=A_m$ with $m\ge 4$.
\begin{itemize}
\item[{\rm (i)}]  For $\mu = (b0\cdots 1d)$, a dominant weight for $X$ with $b\ge 1, d\ge 0$, $\wedge^3\mu$ has a repeated summand of highest weight $\nu = 3\mu - \a_1-\a_2-\cdots -\a_{m-2}-2\a_{m-1}-\a_m$,
 which has $S$-value $S(\nu) \ge 3b+3d+1$.
\item[{\rm (ii)}] For $\mu = (b0\cdots010)$, $b\geq 1$, $\wedge^2(\mu)\otimes \mu - \wedge^3(\mu)$ has a repeated summand of   highest weight
  $\nu = 3\mu-\alpha_1-\cdots-\alpha_{m-1}$ which has $S$-value $S(\nu) = 3b+2$.
\item[{\rm (iii)}] For $\mu = (b0\cdots01d)$, $b,d\geq 1$, $\wedge^2(\mu)\otimes \mu - \wedge^3(\mu)$ has a repeated summand of  highest weight
  $\nu = 3\mu-\alpha_{m-1}-\alpha_m$ which has $S$-value $S(\nu) = 3b+3d+2$.
\end{itemize}
\end{lem}

\pf (i) We view $X<Z = SL(W)$, where $W = V_X(\mu)$. 
Let $V = \wedge^3\mu$. By \cite[Prop. A]{cavallin} we can assume $b = 1$ and $d \le 1$.  Then $V^1 = \wedge^3(10 \cdots 01)$, which contains $(010 \cdots 02)$ by
Lemma \ref{sympowers}.  Moreover, this weight  arises from
$(30 \cdots 03) - (2\a_1 +\a_2 + \cdots + \a_{m-1})$.

 Lemma \ref{sympowers}(iii) implies that $V^2  \supseteq
((010 \cdots 02) + (20 \cdots 010) + (10 \cdots 01)) \otimes ((10 \cdots 010)+(0 \cdots 01)+(10 \cdots 02))$, where the last summand of the second tensor factor does not occur if $d = 0$. 
We have
\[
\begin{array}{l}
(010 \cdots 02) \otimes (10 \cdots 010) \supseteq (20 \cdots 011)^1 \\
(20 \cdots 010) \otimes (0 \cdots 01)  \supseteq (20 \cdots 011)^1 \\
(20 \cdots 010) \otimes (10 \cdots 010) \supseteq (20 \cdots 011)^2 \\
(10 \cdots 01) \otimes (10 \cdots 010) \supseteq (20 \cdots 11)^1 \\
(10 \cdots 01) \otimes (10 \cdots 02) \supseteq (20 \cdots 011)^1.
\end{array}
\]
Therefore  $V^2  \supseteq (20\cdots 011)^5$ if $d = 0$ and $(20\cdots 011)^6$
if $d = 1 $.

Now $(20\cdots 011)$ can only arise from the summands $(110 \cdots 011)$,
$(20 \cdots 010)$, $(20\cdots 02)$, and $(30 \cdots 011)$  in $V^1.$
The last one cannot occur, since otherwise there would be a composition factor of highest weight $(30 \cdots 011)$ in $\wedge^3(10 \cdots 01)$, which is not the case by Lemma \ref{sympowers}(iv). 

As $V^1$ is MF it follows from the above that  there is a composition factor of $V$ of highest weight  $(20\cdots 011x)$ and multiplicity at least 2.

Moreover, from the first paragraph we see that there is also a composition factor of highest weight $(010\cdots 02y)$.  This must have the form $(30 \cdots 03\,(3d)) - (2\a_1 +\a_2 + \cdots + \a_{m-1}) =
(010 \cdots 0 2(1+3d))$.

At this point we use the argument with the torus $T(c)$ given in the proof of Lemma \ref{sympowers}  to find  that $x = 0$ or $3$ according as $d = 0$ or $d = 1$.
Therefore $V$ contains $(20\cdots 0110)^2$ or $(20\cdots 0113)^2$, completing the proof.

(ii) Let $\mu = (b0\cdots 010)$. We first claim that the weight $\nu = 3\mu-(11\cdots 10)$ affords a multiplicity 1
summand of $\wedge^3(\mu)$. The only summand whose highest weight lies above $\nu$ has highest weight
$3\mu-\alpha_1-\alpha_{m-1}$ and therein the weight $\nu$ occurs with multiplicity $m-3$. While in
$\wedge^3(\mu)$ the weight $\nu$ occurs as follows:

$\mu\wedge(\mu-\alpha_1-\cdots-\alpha_i)\wedge(\mu-\alpha_{i+1}-\cdots-\alpha_{m-1})$, for $i=2,\dots,m-2$,

\noindent  and also as

$\mu\wedge(\mu-\alpha_1)\wedge(\mu-\alpha_2-\cdots\alpha_{m-1})$, 

\noindent so a total of $m-2$ times, hence the claim.

Now we show that there are three summands of highest weight $3\mu-\alpha_1-\cdots-\alpha_{m-1}$ in
$\wedge^2(\mu)\otimes \mu$. To this end, note that $\wedge^2(\mu)$ has summands of highest weights
$2\mu-\alpha_{m-1}$ and $2\mu-\alpha_1$. Each of these when tensored with $\mu$ gives rise to a summand
$3\mu-\alpha_1-\cdots-\alpha_{m-1}$ (using \cite[4.1.4]{MF} for example). As well, a counting argument similar to
that of the previous paragraph shows that $\wedge^2(\mu)$ has a summand with highest weight
$2\mu-\alpha_1-\cdots-\alpha_{m-1}$. Then tensoring with $\mu$ gives a third summand of highest weight $\nu$.
Finally we note that $S(\nu) = 3b+2$.

(iii) We note that $\wedge^2(\mu)$ has summands with highest weights $2\mu-\alpha_{m-1}$, $2\mu-\alpha_m$ and $2\mu-\alpha_{m-1}-\alpha_m$. Now tensoring each of these with $\mu$ gives rise to a summand of highest weight $\nu =3\mu-\alpha_{m-1}-\alpha_m$. Moreover it is clear that this weight occurs with multiplicity 1 in $\wedge^3(\mu)$ and
$S(\nu) = 3b+3d+2$.
\hal

\section{The case $\d = b\om_1+b\om_{l+1}$}

In this section we assume that $X = A_{l+1}$ with $l \ge 1$, and let $X < Y = SO(W)$ or $Sp(W)$, where $W=V_X(\d)$ and    $\d = b0\cdots 0b = b\om_1 + b\om_{l+1}.$ It follows from Lemma \ref{setup} that in fact $Y = SO(W)$.   We will prove

\begin{thm} \label{b0...0b} Let $\l \ne 0, \l_1$ be a dominant weight for $Y$.  Then  $V_Y(\l)\downarrow X$ is MF if and only if one of the following holds:
\begin{itemize}
\item[{\rm (i)}] $b = 1$, $l = 1$, and  $\l = \l_2$, $\l_1+\l_4$ or $c\l_1\,(c \le 5)$, up to a twist by an automorphism of $Y=D_4.$
\item[{\rm (ii)}] $b = 1,$ $ l\ge 2 $ and  $\l = \l_2$, $2\l_1$, $3\l_1$, or $\l_3.$
\item[{\rm (iii)}] $b = 2$ and $\l = \l_2$.
\end{itemize} 
\end{thm}

\subsection{The case $b=1$}


Here we assume that $\d = \o_1+\o_{l+1}$, so that $X < Y = SO(W) = SO_n$ for $n = (l+2)^2-1$. Also 
 $L_X'$ is contained in $L_Y' = C^0 \times C^1,$ where $L_X'$ maps to  $ C^0$ via the representation $\om_1$, 
and to  $C^1$ via the sum of the representations  $\om_1 + \om_l$ and $0.$ We will frequently use the fact that this sum can also be realized as $\om_1\otimes \om_l.$  Note that $C^1$ is an orthogonal group with natural module of dimension $(l+1)^2.$ 

If $l$ is odd  then $C^1 = D_t$, and if $l$ is even  then $C^1 = B_t$; in the first case $L_X'$ projects into $B_{t-1}$ and in the second into $D_t$.  In either case let $D$ be the smaller orthogonal group $B_{t-1}$ or $D_t$. The embedding of $L_X'$ in $D$ is given by the representation
$\om_1 + \om_l$ so that we are in a position to use induction.   The fundamental system for $C^1$ is $\la \b_1^1, \cdots, \b_t^1\ra$. For the purposes of restricting weights of $C^1$ to $D$, we  may take the fundamental system for $D$ to be $\la \b_1^1, \cdots,\b_{t-1}^1, \b_{t-1}^1 + 2\b_t^1 \ra$ if $C^1 = B_t$ and $\la \b_1^1, \cdots,\b_{t-2}^1, \b_{t-1}^1 + \b_t^1 \ra$ if $C^1 = D_t.$

Our proof of Theorem \ref{b0...0b} for this case $(b=1$) will go by induction on the rank $l+1$. The next two results handle the cases $l=1$ or 2. 

\begin{prop}\label{l=1}  Assume $X = A_2$ and $X < Y = D_4$ via the representation $\d = (11).$ Let $\l$
be a dominant weight for $Y$, and assume $\l \ne 0$ or $\l_i$ for $i=1,3,4$.  Then, up to a twist by an automorphism of $D_4,$ $V_Y(\l)\downarrow X$ is MF if and only if $\l = \l_2$, $\l_1+\l_4$, or $c\l_1\, (c \le 5)$.
\end{prop}

 \pf The fact that $V_Y(\l)\downarrow X$ is MF for the values of $\l$ in the conclusion can be verified using Magma together with Lemma \ref{symwed}.

Now assume that $V_Y(\l)\downarrow X$ is MF. We have 
 $L_Y' = A_1^3$. Let $\l = adec$, with labels $a,e,c$ on the factors of $L_Y'$. Then $V^1 = a\otimes e \otimes c$ is multiplicity-free and this has highest weight $m = a+e+c.$  At least one of $a,e,c$ is 0. Say $c=0$.   
Assume first that $a=e=0$. If $d=1$ then $\l=\l_2$ is as in the conclusion, so take $d\ge 2$. We have $m=0$ and $V^2(Q_Y)$ is afforded by $\l-\b_2$, and so $V^2 = 
1\otimes 1\otimes 1 = 3+1^2$. Now $V^3$ contains $2\otimes 2\otimes 2$ (afforded by $\l-2\b_2$), and also three summands $2$ (afforded by $\l-2\b_2-\b_i-\b_j$ for $ij = 13,14$ or $34$). Hence $V^3$ contains $2^6$, which is a contradiction by Proposition \ref{induct}.

Now suppose $a\ge e>0$. Here $m=a+e$. If $d\ne 0$, then $\l-\b_2$ affords the  $L_X'$-summand $(a+1)\otimes (e+1)\otimes 1$ of $V^2$, while $\l-\b_2-\b_1$ affords $(a-1)\otimes (e+1)\otimes 1$. Between them these contain $(a+e+1)^3$, a contradiction by Proposition \ref{induct}(ii).  Therefore $d=0$. Then $\l-\b_1-\b_2,\,\l-\b_2-\b_3,\,\l-\b_1-\b_2-\b_3$ afford  $V^2$-summands $(a-1)\otimes (e+1)\otimes 1$, $(a+1)\otimes (e-1)\otimes 1$, $(a-1)\otimes (e-1)\otimes 1$ respectively. If $(a,e) \ne (1,1)$, then between them these contain $(a+e-1)^4$, contradicting Proposition \ref{induct}(ii). And if $(a,e)=(1,1)$ then  $\l$ is as in the conclusion.

Next assume $ad \ne 0$ and $e = 0$.  Then $V^2 $ contains $(a+1) \otimes 1 \otimes 1$ 
and $(a-1) \otimes 1 \otimes 1.$  The first summand contains $(a+1)^2$ and
the second summand contains $a+1$.  Therefore $V^2 \supseteq (a+1)^3$ and hence 
$V_Y(\l) \downarrow X$ is not MF by Proposition \ref{induct}(ii) again.

The final case now is where $e=d=0$, so $\l = a\l_1$. If $a\le 5$ then $\l$ is in the conclusion, so suppose $a \ge 6.$
Here we claim that $S^a(11)$ contains $(a-3,a-3)^3.$  Lemma \ref{symwed} together with a Magma 
computation gives this for $a = 6$, 
and a weight consideration shows that it also holds for $a > 6.$ Hence $V_Y(\l)\downarrow X$ is not MF, as 
$S^{a-2}(11)$ only contains $(a-3,a-3)^1.$ This completes the proof. \hal

The next result settles the case $l = 2$, where $\d=101$ and $Y = B_7 = SO_{15}.$ 

\begin{prop}\label{l=2} Assume $X = A_3$ and $\d = \om_1 + \om_3.$ Let $\l$ be a dominant weight for $Y=B_7$, with $\l \ne 0,\,\l_1$. 
Then $V_Y(\l)\downarrow X$ is MF if and only if $\l = \l_2$, $2\l_1$, $3\l_1$, or $\l_3.$
\end{prop}

\pf Assume the hypotheses of the proposition. A strightforward check using Lemma \ref{symwed} and a Magma 
computation shows that $V_Y(\l)\downarrow X$ is MF for all the weights $\l$ in the conclusion.

Now assume that $V_Y(\l)\downarrow X$ is MF.  We have $L_X' = A_2$ and  $L_X'$ is embedded in $C^0 =A_2$ and $C^1= B_4$ via the representations $(10)$ and $(11)+00$, respectively.  

\vspace{2mm} {\it Step 1. } $\mu^1$ is either the trivial or natural representation of $C^1 = B_4$.

\vspace{2mm} To prove this, assume false and write $\mu^1 = a_1\l_1^1 + a_2\l_2^1 +a_3\l_3^1 +a_4\l_4^1.$ Now $L_X'$ maps into the $D_4$ subsystem subgroup of $B_4$ with
fundamental system $\{ \b_1^1, \b_2^1, \b_3^1, \b_3^1+2\b_4^1\}$.  If $\nu_1, \nu_2, \nu_3, \nu_4$ are the corresponding fundamental weights for $D_4,$ then $\mu^1$ restricts to the
dominant weight $a_1\nu_1 + a_2\nu_2 + a_3\nu_3 +(a_3+a_4)\nu_4$ for $D_4.$ Applying Proposition \ref{l=1}
and recalling that $V^1$ is MF, 
 we see that $\mu^1 = \l_2^1$, $\l_1^1+\l_4^1$, $\l_4^1$, $\l_3^1$, $c\l_4^1$, or $c\l_1^1,$ where $2 \le c \le 5$. We will restrict these representations to $D_4$.  Note that there is an involution $\tau$ in $B_4$ inducing a graph automorphism of $D_4$ and fixing $L_X'$.  Indeed, we can take $\tau = s_4$, the fundamental reflection corresponding to the root $\b_4^1$. This can be used  to our advantage.  In addition we can make use of a triality morphism of $D_4$ fixing $L_X',$ although this is not in $ B_4.$

If $\mu^1 = \l_2^1$, then $V_{C^1}(\mu^1) \downarrow L_X' = \wedge^2(11+00) \supseteq (11)^2$ by Lemma \ref{sympowers},  
whereas this restriction  must be MF by Proposition \ref{induct}(i). Similarly if $\mu^1 =  \l_3^1$, then $V_{C^1}(\mu^1) \downarrow L_X' = \wedge^3(11+00) \supseteq 11^2$, again a contradiction. If $\mu^1 = \l_4^1,$ then using $\tau$ we see that  $V_{C^1}(\mu^1) \downarrow D_4 \supseteq V_{D_4}(\nu_3) + V_{D_4}(\nu_4).$ These are both conjugate to $V_{D_4}(\nu_1)$ under triality and so the restriction to $L_X'$ contains $(11)^2.$  Next suppose $\mu^1 = \l_1^1 + \l_4^1.$ Then  using $\tau$  we see that $V_{C^1}(\mu^1) \downarrow D_4 \supseteq (\nu_1+\nu_4) + (\nu_1 + \nu_3)$. Triality shows that each of these restricts to $L_X'$ as does  $(\nu_3+\nu_4),$  a contradiction.  If $\mu^1 = c\l_4^1$,
then the restriction to $D_4$ contains $c\nu_4 $ and $c\nu_3,$ both of which are conjugate under triality to $c\nu_1$, and again this is a contradiction.  Finally suppose $\mu^1 = c\l_1^1$ for $2 \le c \le 5.$ Then $V_{C^1}(\mu^1) = S^c(\mu^1) - S^{c-2}(\mu^1)$ and restricting this to $L_X'$ we obtain $S^c(11+00) - S^{c-2}(11+00)$, which gives
 $((c-1)(c-1))^2.$ So this establishes Step 1.

\vspace{2mm} {\it Step 2. } We have $\mu^1 = 0.$

\vspace{2mm} Suppose false, so that  Step 1 shows that $\mu^1 = \l_1^1.$  Then $V_{C^1}(\mu^1) \downarrow L_X' = (11) + (00)$, so Lemma \ref{stemb} implies that  $\mu^0 = (a0)$ or $(0a)$.  If $a\ne 0$, then tensoring with $(11) + (00)$ we get $(a0)^2$ or $(0a)^2$ in $V^1$, respectively, so this is a contradiction.  Therefore $\mu^0 = 0.$ Then $V^1 = (11) + (00)$ and using Lemma \ref{sympowers} we see that $V^2$ has a summand $(01) \otimes \wedge^2(11 + 00) = (01) \otimes (03 + 30 + 11^2) \supseteq (12)^3.$ But only one copy of $(12)$ can arise from $V^1$, so this is a contradiction. This completes Step 2.

\vspace{2mm} 
We now proceed with the proof of Proposition \ref{l=2}.
Write $\mu^0 = (ac)$ so that $V^1 = (ac).$
If $c \ne 0$, then $V^2 \supseteq ((a+1)(c-1)) \otimes ((11) + (00))$.
Similarly if $a \ne 0$, then $V^2 \supseteq ((a-1)c) \otimes ((11) + (00))$.

Suppose $c \ge 2$.  Then \cite[Lemma 7.1.5]{MF} implies that $V^2 \supseteq ((a+1)(c-1))^3$,
whereas only one such term can arise from $V^1.$ So $ c \le 1.$ Similarly we get a contradiction if $a \ge 2$ and $c \ne 0.$

Suppose $a = c = 1.$ Then  $V^2 \supseteq (20 + 01) \otimes (11 + 00) \supseteq (20)^3.$ This is a contradiction since only one such term can arise from $V^1.$

Next we rule out the cases $(ac) = (10),$ $(01)$, and $(00).$ 
If $\la \l,\g_1 \ra =0$ then $\l$ is as in the conclusion, so assume $\la \l, \g_1 \ra \ne 0$. If 
$(ac) = (10)$, then $V^2 \supseteq (11 +00) \otimes (11 +00) \supseteq
(11)^4$, and only one such term arises from $V^1.$  And if $(ac) = (01),$
then $V^2 \supseteq (10 + 02) \otimes (11+00) \supseteq (02)^3$, again a contradiction.  Finally assume that  $(ac) = (00).$  If $\la \l,\g_1 \ra =1$ then $\l=\l_3$ as in the conclusion, so assume that $\la \l, \g_1 \ra \ge 2.$ 
Then $V^2 = (01) \otimes (11 +00) = 12 + 20 + (01)^2$ and using Lemma \ref{sympowers}, $V^3 \supseteq (02) \otimes (S^2(11+00)-00) = (02) \otimes (22+(11)^2 +00) \supseteq (13)^3.$ Only one such term can arise from $V^2$, so this is again a contradiction.

So at this point we can assume that $c = 0$ and $a \ge 2.$ We next argue that $\la \l, \g_1 \ra = 0.$ Otherwise $V^2 \supseteq (a1+(a-1)0) \otimes (11 +00) \supseteq (a1)^4$, a contradiction.  We have therefore shown that $\l = a\l_1$. 
 Suppose $a \ge 4.$ We have $V_Y(\l) = S^a(\l_1)-S^{a-2}(\l_1)$, and  Lemma \ref{sympowers}  shows that $S^a(101) \supseteq ((a-2)0(a-2))^3$.  As this is the highest weight of $S^{a-2}(\l_1)$ we see that $V_Y(\l)\downarrow X \supseteq ((a-2)0(a-2))^2$ and so the restriction is not MF, a  contradiction.  Therefore $\l = a\l_1$ with $a \le 3$, as in the conclusion.  This completes the proof of  Proposition \ref{l=2}. \hal

The next result, whose proof relies on Lemmas \ref{5.5}--\ref{5.8}, completes the proof of Theorem \ref{b0...0b} for $b=1$.

\begin{prop} \label{b=1,l>2} Assume $X = A_{l+1}$ with $l \ge 3$ and $\d = \om_1 + \om_{l+1}$ and $X$ is embedded in $Y = SO(W) = SO_d$, where $W = V_X(\d)$ and $d = (l+2)^2-1.$ Let $ \l \ne 0,\l_1$ be a dominant weight for $Y.$  Then $V_Y(\l) \downarrow X$ is MF if and only if $\l =  \l_2$, $2\l_1$, $3\l_1$, or $\l_3.$
\end{prop}

Before starting the proof, we first observe that $V_Y(\l) \downarrow X$ is indeed MF for each of the weights $\l$ in the conclusion of the proposition: this follows from \cite{MF} for all except $3\l_1$, and in this case it follows from Lemma \ref{sympowers}(ix), noting that $V_Y(3\l_1) = S^3(\l_1)-\l_1$. 

We work towards the proof of Proposition \ref{b=1,l>2}. Assume that $V_Y(\l) \downarrow X$ is MF, and 
write $\mu^1 = a_1\l_1^1 + \cdots + a_t\l_t^1.$ It follows from the discussion at the beginning of this subsection that $V_{C^1}(\mu^1) \downarrow D$
has a composition factor of highest weight $a_1\nu_1+ \cdots + a_{t-1}\nu_{t-1} + (a_{t-1}+a_t)\nu_t$ or $a_1\nu_1 + \cdots + a_{t-2}\nu_{t-2} + (a_{t-1}+a_t)\nu_{t-1},$ according to whether $C^1 = B_t$ or $D_t.$ The irreducible $C^1$-module with this highest weight must be MF on restriction to $L_X'$,
so by induction we have $\mu^1 = 0$,  $\l_1^1$, $\l_2^1$, $2\l_1^1$, $3\l_1^1$, or $\l_3^1.$

\begin{lem}\label{5.5} $\mu^1 = 0$ or $\l_1^1.$
\end{lem}

\pf  By way of contradiction assume $\mu^1 = \l_2^1$, $2\l_1^1$, $\l_3^1,$ or $3\l_1^1.$  It follows that
$V_{C^1}(\mu^1) = \wedge^2(V_{C^1}(\l_1^1)),$ $S^2(V_{C^1}(\l_1^1)) - 0,$ $\wedge^3(V_{C^1}(\l_1))$, or $S^3(V_{C^1}(\l_1^1)) - V_{C^1}(\l_1^1)$, respectively.  We claim that in each case the restriction contains
$(\om_1 + \om_l)^2,$ which will establish the lemma.

Suppose $\mu^1 = \l_2^1.$ Then $V_{C^1}(\mu^1) \downarrow L_X' = \wedge^2((\om_1 + \om_l) +0) =\wedge^2(\om_1 + \om_l)  + (\om_1 + \om_l).$    Lemma  \ref{sympowers}(iii) shows  that the first summand contains $\om_1 + \om_l,$ so this gives the claim in this case.  Similarly for $\mu^1 = 2\l_1^1,$ but this time we use  Lemma \ref{sympowers}(ii). 

If $\mu^1 = \l_3^1$ then $\wedge^3(V_{C^1}(\l_1)) \downarrow  L_X'  = \wedge^3(\om_1+\om_l) + 
\wedge^2(\om_1+ \om_l)$ and parts (iii) and (iv) of Lemma \ref{sympowers} show that   $\om_1+\om_l$ occurs in each summand. So again the claim holds. 

 Finally assume that $\mu^1 = 3\l_1^1.$  From the first paragraph we see that the restriction to $L_X'$ contains
 $S^3((\om_1+\om_l) + 0) - ((\om_1+\om_l) + 0).$    
 Now $S^3((\om_1+\om_l) + 0) = S^3(\om_1+\om_l) + S^2(\om_1+\om_l) + (\om_1 + \om_l) + 0.$   Therefore $V_{C^1}(3\l_1) \downarrow L_X' = S^3(\om_1+\om_l) + S^2(\om_1+\om_l)$ and Lemma 
\ref{sympowers}(vii) - (ix) show that $S^3(\om_1+\om_l) \supseteq (\om_1 + \om_l)^2$, a contradiction.   \hal
 
 \begin{lem}\label{5.6} $\mu^1 = 0.$
  \end{lem} 
 
 \pf Suppose $\mu^1 = \l_1^1,$ so that $V_{C^1}(\mu^1) \downarrow L_X' = (\om_1 + \om_l) + 0.$
 As $V^1$ is MF we can apply Lemma \ref{stemb} to get $\mu^0 = a\l_i^0$ for some $i$.  If $a = 0$, 
 then $V^2 \supseteq (0 \cdots 01) \otimes \wedge^2((10\cdots 01) + 0) = (0 \cdots 01) \otimes ((10\cdots 01)^2 + (010 \cdots 02) + (20 \cdots 01)) \supseteq  (10 \cdots 02)^3$.  This is a contradiction as $V^1 = (10 \cdots 01) + 0.$ Therefore
 $a > 0.$
 
 We have $V^1 = a\om_i \otimes  ((10\cdots01) + 0)).$ An argument using the Littlewood-Richardson rule (as in \cite[\S 4.1]{MF})   shows that $a\om_i \otimes  (10\cdots01) \supseteq a\om_i.$ Therefore $V^1 \supseteq a\om_i^2$, a contradiction.  \hal

 \begin{lem}\label{5.7}  Either $\mu^0 = a_1\l_1^0$ with $a_1 \ge 0$, or $\mu^0 = \l_i^0$ with $i>1$.
  \end{lem}
 
 \pf   Write $\mu^0 = (a_1 \cdots a_l)$ and suppose the result is false. Then $\mu^0 \ne a_1\l_1^0$ and there exists $i >1$ with $a_i \ne 0.$ Take $i$ minimal for this.  Then
 $V^2$ contains a composition factor with highest weight  $(a_1 \cdots (a_{i-1}+1)(a_i-1) \cdots a_l) \otimes ((10\cdots 01) + 0). $  Hence \cite[Lemma 7.1.5(i)]{MF} implies that $V^2$ contains $(a_1 \cdots (a_{i-1}+1)(a_i-1) \cdots a_l)^3$ unless $(a_1 \cdots (a_{i-1}+1)(a_i-1) \cdots a_l)$ contains only a single nonzero term. So this would be  a contradiction  since only one such term can arise from $V^1.$ Therefore assume the sequence $(a_1 \cdots (a_{i-1}+1)(a_i-1) \cdots a_l)$ has only one nonzero term, namely the entry in the $(i-1)$ spot.
 Then  $a_i =1$ and  $a_j = 0$ for $j \ne i,i-1. $
 
Note that $a_{i-1} \ne 0$, since $\mu^0\ne \l_i^0$ by assumption. 
The minimality of $i$ implies that $i = 2$ and $\mu^0 = (a_110\cdots0).$ Since we are assuming that the result is false
 we have $a_1 \ne 0.$ Then $V^2$ has summands $((a_1+1)0\cdots0)\otimes((10\cdots01)+0)$ and $((a_1-1)10\cdots0)\otimes((10\cdots01)+0)$ which sum to $(a_10\cdots0)\otimes(10\cdots0)\otimes(10\cdots 0)
 \otimes (0\cdots01).$
 Expanding we get terms 
 $(a_10\cdots0)\otimes(010\cdots0) \otimes (0\cdots01)$
 and $(a_10\cdots0)\otimes(20\cdots0) \otimes (0\cdots01).$
 The first contains $(a_110\cdots0)\otimes (0\cdots01)$ and hence by Cor. 5.1.5 of \cite{MF} it covers summands arising from $V^1$ while the second contains $((a_1+1)0\cdots 0)^2.$
This is a contradiction.    \hal
 
\begin{lem}\label{5.8} $\la \l, \g_1 \ra = 0.$
\end{lem}

\pf Set $\la \l, \g_1 \ra = d.$  By way of contradiction suppose $d \ne 0.$  Then $V^2 \supseteq (\mu^0\downarrow L_X' + \om_l) \otimes ((\om_1 + \om_l) + 0).$ If $\mu^0 = a_1\l_1^0$ with $a_1 \ne 0$ or if $\mu^0 = \l_i^0$ with $1 < i < l$ then we find that $V^2 \supseteq (\mu^0\downarrow L_X'  + \om_l)^3$ and
only one such term can arise from $V^1.$ So this is a contradiction.

Now suppose that $\mu^0 = \l_l^0$. Then $V^2$ contains
$2\om_l \otimes (\om_1 + \om_l)$ as well as $\om_{l-1} \otimes (\om_1 + \om_l).$ These each contain $(\om_1 + \om_{l-1} + \om_l)$, and this again yields a contradiction.

To establish the lemma we still need to settle the case $\mu^0 = 0$.
If $d = 1$, then $V_Y(\l) = \wedge^{l+1}(\l_1)$, whereas \cite{MF} shows that $ \wedge^{l+1}(\l_1)\downarrow X$ is not MF.  So suppose $d \ge 2.$ Then $V^3 \supseteq
2\om_l \otimes V_{C^1}(2\l_1^1) \downarrow L_X'$.  The second tensor factor restricts to $L_X'$ as $S^2((\om_1 + \om_l) +0) - 0$ and
this contains $(2\om_1 + 2\om_l) \oplus (\om_1+\om_l)^2.$ Therefore tensoring with $2\om_l$ we have $(\om_1 +3\om_l)^3$.  However,
$V^2 = \om_l \otimes ((\om_1 + \om_l) + 0)$ so that only one
summand $(\om_1 +3\om_l)$ can arise from $V^2.$ This contradiction establishes the lemma.  \hal 

\vspace{2mm}
We now complete the proof of  Proposition \ref{b=1,l>2}.  First suppose that 
$\mu^0 = \l_i^0$ for $i > 1.$ Then $V_Y(\l) \downarrow X = \wedge^i(\om_1 + \om_{l+1})$.  Then \cite{MF} implies that $i = 2$ or $3$  and as already observed (just after the statement of \ref{b=1,l>2}), these representations are indeed MF.  So the proposition holds in this case.

Finally consider the case $\mu^0 = a_1\l_1^0$ so that $\l = a_1\l_1.$
We can assume $a=a_1 > 1.$ Then $V = S^a(\l_1) - S^{a-2}(\l_1)$ so restricting to $X$ we have $S^a(\om_1 + \om_{l+1}) - S^{a-2}(\om_1+\om_{l+1}).$  We apply Lemma \ref{sympowers}.  Parts (i),(viii), and (ix) of the lemma settle the cases  $a = 2$ and $a =3.$ 
So assume $a \ge 4.$ Then Lemma  \ref{sympowers}(v) shows that  $S^a(\om_1 + \om_{l+1})\supseteq ((a-2)\om_1 + (a-2)\om_{l+1})^3.$ Since $((a-2)\om_1 + (a-2)\om_{l+1})$ is the highest weight of $S^{a-2}(\om_1 + \om_{l+1})$,
we conclude that $S^a(\om_1 + \om_{l+1}) - S^{a-2}(\om_1 + \om_{l+1}) \supseteq ((a-2)\om_1 + (a-2)\om_{l+1})^2$,
and hence $V_Y(a\l_1) \downarrow X$ is not MF, completing the proof of  Proposition \ref{b=1,l>2}.

\subsection{The case $b\ge2$}

In this subsection we complete the proof of Theorem \ref{b0...0b} by proving 

\begin{prop} \label{b>1} Assume $X = A_{l+1}$ with $l \ge 1$, let $\d = b\om_1 +b\om_{l+1}$ with $b \ge 2$, and let $X<Y = SO(W)$, where $W = V_X(\d)$. Let $ \l \ne 0,\l_1$ be a dominant weight for $Y.$  Then $V_Y(\l) \downarrow X$ is MF if and only if $b=2$ and $\l =  \l_2$. 
\end{prop}

Observe first that \cite{MF} shows that if $b=2$, then $V_Y(\l_2)\downarrow X$ is indeed MF, giving the right to left implication in the proposition.

We will establish  Proposition \ref{b>1} using induction on $l$.  Assume the hypotheses, and note that Lemma \ref{setup} implies that $Y$ is an orthogonal group and  $k = b$. If $l \ge 2$  it is easy to see that $V_{C^k}(\l_1^k)\downarrow L_X' = (b0\cdots0b)+((b-1)0\cdots0(b-1)) + \cdots +(10\cdots 01) + 0$, while if $l = 1$ then $V_{C^k}(\l_1^k)\downarrow L_X' = 2b + (2b-2) + \cdots +0.$

\begin{lem}\label{a2} Proposition $\ref{b>1}$ holds if $X = A_2$.
\end{lem}

\pf Assume that $X = A_2$ and note that the embeddings of $L_X'$ in $C^j$ are given as in the diagram of Section 5.3 of \cite{MF}, although in the diagram we have $r = s = b$ and we only use the top portion of the diagram.
We will use  arguments from Chapter 8 of \cite{MF}.  While these arguments were for an overgroup of $X$ of type $A$,
 the arguments used are local in nature and apply to the current situation.

The proof of \cite[Lemma 8.1.5(a)]{MF}  shows that if $j \ne 0 $, then $\mu^j $ is either zero or the highest weight of the natural module or its dual (this requires also \cite[Prop. 6.1]{LST} for the case where $j=k$ and $C^k$ is orthogonal). Suppose there is a nonzero $\l$-label on a factor $C^j$ of $L_Y'$, with $j \ne 0$. Then the arguments of
\cite[Lemmas 8.1.6, 8.1.7]{MF} give a contradiction.  It follows that the only possible nonzero label is on $C^0$.  Then Lemma 8.1.9 of \cite{MF} shows that $\la \l, \g_i \ra = 0$ for
all $i$ and the arguments in Lemmas 8.1.12-8.1.18 of \cite{MF} (noting that $r=s$ here) establish the result. \hal

From now on we suppose that $l \ge 2.$ Assume that $V_Y(\l) \downarrow X$ is MF.

\begin{lem}\label{mukposs} $\mu^k = 0$ or $\l_1^k$.
\end{lem}

\pf Write $\mu^k = \sum c_i\l_i^k.$ We see from the above that $L_X'$ is embedded in a subgroup $F \le  C^k$ where $F = F_b \times \cdots \times F_1$
is a product of orthogonal groups and $L_X'$ is embedded in $F_j$ via the representation with highest weight $(j0\cdots0j)$.  Now $V_{C^k}(\mu^k) \downarrow L_X'$  is MF, so each composition factor of $F$ on $V_{C^k}(\mu^k)$ is MF when restricted to $L_X'$ and hence each irreducible summand of $V_{C^k}(\mu^k) \downarrow F_b$ must restrict to an MF representation of $L_X'$.
Inductively we see that each composition factor of $F_b$ on $V_{C^k}(\mu^k)$ is either trivial, a natural module, or if
$b = 2$  the alternating square of the natural module.
Set $Z = F_b$ and note that  $Z$ has rank at least 13. 

 First suppose that $Z$ and $C^k$ are  of the same type (i.e. both of type $D$ or of type $B$).  Then we can regard $Z$ as  a standard Levi group of $C^k.$ Say $Z$ has fundamental system $\la \b_s^k, \cdots, \b_{r_k}^k\ra$ where $s \ge 6.$  Let $\tilde Z$ denote the conjugate   of $Z$ with fundamental system $\la \b_1^k, \b_2^k,\b_3^k+ \cdots + \b_{s+2}^k,\b_{s+3}^k, \cdots, \b_{r_k}^k\ra$.  Then $\tilde Z$ has a composition factor of highest weight $(c_1, c_2, (c_3+ \cdots + c_{s+2}), c_{s+3}, \cdots, c_{r_k})$.  Inductively we have $c_3 = \cdots = c_{r_k} = 0$ and $(c_1,c_2) =  (0,0)$, $(1,0)$, or $(0,1),$ where the third
possibility only occurs if $b = 2$. Suppose that $b = 2$ and $(c_1,c_2) = (0,1)$, so that $V_{C^k}(\mu^k)$ is the wedge square of the orthogonal module. Then there is a composition factor for $F = F_2  \times F_1$ which is the tensor product
of natural modules. This restricts to $L_X'$ as $(20\cdots 02) \otimes (10 \cdots 01)$ and this contains $(20\cdots 02)^2$ by \cite[Lemma 7.1.5]{MF}, a contradiction.  Therefore the conclusion holds in this case.
 
 Now suppose  $Z$ and $Y$ are not of the same type.  Then $Z$ is contained in a proper Levi subgroup of $C^k$ and this Levi subgroup has rank equal to or one greater than the rank of $Z$, according to whether $Z$ has type $D$ or type $B$.  We can choose fundamental systems of roots and corresponding dominant weights for $Z$ as in the second paragraph of Section \ref{setupp}.  At this point the  above argument   again gives the result. \hal

\begin{lem} $\mu^i = 0$ for $i \ne 0.$
\end{lem}

\pf The proofs of Theorems 14.1 and 14.2 of \cite{MF} imply that $\mu^i = 0$ for $0 < i < k-1.$ 
The proof of Theorem 14.1  also shows that $\mu^{k-1} = 0, \l_1^{k-1}$ or $\l_{r_{k-1}}^{k-1}$,  and the argument of Theorem 14.2 shows that $\mu^{k-1} \ne \l_1^{k-1}$.  

Suppose that $\mu^{k-1} =\l_{r_{k-1}}^{k-1}.$ 
One of the composition factors of $L_X'$ on
the natural module for $C^{k-1}$ is $(b0\cdots0(b-1)).$ Therefore $V^2 \supseteq (\mu^0 \downarrow L_X') \otimes \wedge^2((b-1)0\cdots 0b) \otimes (b0\cdots0b)$
or
 $V^2 \supseteq (\mu^0 \downarrow L_X') \otimes \wedge^2((b-1)0\cdots 0b) \otimes ((2b)0\cdots0(2b))$, according to whether $\mu^k = 0$ or $\l_1^k$.  Now $\wedge^2((b-1)0\cdots b)  \supseteq ((2b-2)0 \cdots 01(2b-2))$.
Hence $V^2$ contains $ (\mu^0 \downarrow L_X') \otimes ((3b-3)0\cdots 01(3b-3))^2$  or $ (\mu^0 \downarrow L_X') \otimes ((4b-3)0\cdots 01(4b-3))^2$  and hence a repeated composition factor of $S$-value $S(\mu^0) +6b-5$ or $S(\mu^0)+8b-5$, respectively.  As $S(V^1) = S(\mu^0) + 2b-1$ or $S(\mu^0) + 4b-1$, this is a contradiction by Lemma \ref{prop38}.  Therefore, $\mu^{k-1} = 0$.
 
 Suppose $\mu^k = \l_1^k$.
Then $V^2 \supseteq (\mu^0 \downarrow L_X') \otimes (\l_{r_{k-1}}^{k-1} \downarrow L_X') \otimes \wedge^2((b0\cdots0b)+((b-1)0\cdots0(b-1)) + \cdots + (10 \cdots 01)).$  Therefore  \cite[Lemma 7.1.7(ii)]{MF} implies that 
\[
\begin{array}{ll}
V^2 & \supseteq  (\mu^0 \downarrow L_X') \otimes ((b-1)0\cdots0b) \otimes ((2b-1)0\cdots0(2b-1))^2 \\
       & \supseteq  (\mu^0 \downarrow L_X') \otimes ((3b-2)0\cdots0(3b-1))^2.
\end{array}
\]
We are assuming $b \ge 2$, so an $S$-value argument gives a contradiction by Lemma \ref{prop38}.   Therefore, Lemma \ref{mukposs}  implies that $\mu^k = 0$. \hal

\begin{lem} $\la \l, \gamma_i \ra = 0$ for all $i$.
\end{lem}

\pf If $b \ge 3$, then this follows from the proof of Theorem 16.1 of \cite{MF}. So suppose $b = 2.$ 
If $\la \l, \gamma_2 \ra \ne 0$, then
arguing as in the previous lemma we see that $V^2 \supseteq (\mu^0 \downarrow L_X') \otimes ((b-1)0\cdots0b) \otimes ((b0\cdots0b)+((b-1)0\cdots0(b-1)) + \cdots)$
and this contains $(\mu^0 \downarrow L_X') \otimes ((2b-2)0\cdots 0(2b-1))^2$ by \cite[Lemma 7.1.5(i)]{MF}.
   Then an $S$-value argument gives a contradiction.

Now suppose $\la \l, \gamma_1 \ra \ne 0.$ Then \cite[Prop. 5.4.1]{MF} shows that
$V^2_{\gamma_1}(Q_Y)  \supseteq V_{C^0}(\mu^0) \otimes V_{C^0}(\l_{r_0})\otimes V_{C^1}(\l_1^1).$ Restricting to $L_X'$ this is
$\mu^0\downarrow L_X' \otimes (0\cdots02)\otimes((20\cdots01)+(10\cdots0))$ which contains
$(\mu^0\downarrow L_X') \otimes (10\cdots02)^2$ and again an $S$-value argument gives a contradiction.  \hal

At this point it remains to determine $\mu^0$ as that will determine $\l.$ The main result of  \cite{MF}  gives a list of  possibilities for $\mu^0$ and this list must be reduced to those in the conclusion of Proposition \ref{b>1}.  
For $b=2$ we can use  \cite[Lemmas 17.3.4, 17.3.5]{MF} to rule out nearly all possibilities since the proofs of those lemmas are local in the sense that they only involve $\mu^0,$ $C^0,$ and $C^1$ rather than the structure of the group $Y$, so they  apply equally well in this situation. Hence we need to rule out cases resulting from Lemma 17.3.5 of \cite{MF}, together with inductive cases for  $b \ge 3.$ 
That is,  for $b = 2$ we assume that  $\mu^0$ or its dual is in $\{ \l_1^0, \l_2^0, 2\l_1^0 \}$, and for $b \ge 3$ we assume that  $\mu^0$ or its dual is in  $\{ c\l_1^0\,(2\le c \le 5),\, \l_i^0\,(2\le i\le 5),\, \l_1^0 + \l_{r_0}^0,\, \l_1^0 + \l_2^0 \}.$

\begin{lem} $\mu^0 \ne c\l_1^0$ or $c\l_{r_0}^0$ for $c \ge 2.$
\end{lem}

\pf  Suppose $\mu^0 = c\l_1^0$ with $c \ge 2.$   It follows that $V_Y(\l) = V_Y(c\l_1) = S^c(V_Y(c\l_1)) - S^{c-2}(V_Y(c\l_1)).$  Restricting to $X$ we get $S^c(b0\cdots0b)-S^{c-2}(b0\cdots0b).$ 

If $c = 2$, then Lemma 7.1.9(ii) of \cite{MF} shows that $S^2(b0\cdots0b) \supseteq ((2b-2)0\cdots0(2b-2))^2,$ a contradiction.  Now assume $c \ge 3.$
We claim that $S^c(b0\cdots0b)$ has a composition factor with highest weight $\gamma = c\d - 2\a_1-2\a_2$ with multiplicity at least 2.  Indeed, in view of the simple nature of $\gamma$ we can count weights to see  that $S^c(b0\cdots0b) \supseteq \gamma^2.$   Then the $S$-value of this composition factor is $2cb-4$.  On the other hand the highest weight of $S^{c-2}(b0\cdots0b)$ has $S$-value $(c-2)b + (c-2)b = 2cb-4b,$  so this is a contradiction by Lemma \ref{prop38}.

Now suppose $\mu^0 = c\l_{r_0}^0$.  If $b \ge 3$, then Theorem 15.1 of \cite{MF} gives a contradiction.  So assume $b = 2$.  
Then
$V^2 \supseteq ((0 \cdots 01(c-1))\downarrow L_X') \otimes (20 \cdots 01)$.
The first tensor factor contains a composition factor of highest weight
$(0 \cdots 01(2c))$ and hence $V^2 \supseteq (10 \cdots 01(2c))^2$.  An $S$-value argument now gives a contradiction.
 \hal

\begin{lem}  $\mu^0 \ne \l_1^0 + \l_{r_0}^0$ if $b \ge 3.$
\end{lem}

\pf Assume false.  Then $V^1$ has $S$-value $2b$.  On the other hand $V^2 \supseteq
(((\l_1+\l_{r_0-1}) \oplus \l_{r_0})\downarrow L_X') \otimes ((b0\cdots01)+((b-1)0\cdots0)) = ((\l_1\otimes \l_{r_0-1})\downarrow L_X') \otimes
(b0\cdots0)\otimes(0\cdots01).$ This contains 
$(b0\cdots0) \otimes (0\cdots01(2b-2))\otimes (b0\cdots0)\otimes(0\cdots01).$

From the first and third tensor factors there is a summand $((2b-2)10\cdots0)$
and hence Lemma 7.1.6 of \cite{MF} implies that  the full tensor product contains $((2b-2)0\cdots 0(2b-1))^2.$ 
Now an $S$-value argument gives a contradiction. \hal

\begin{lem} If $b = 3,$ then $\mu^0 \ne \l_1^0 + \l_2^0$ or its dual.
\end{lem}

\pf Theorem 15.1 of \cite{MF} rules out the dual, and the argument of Lemma 7.2.34(ii) of \cite{MF} together with 
Lemma \ref{tensorl1} shows that $\l \ne \l_1 + \l_2$.  \hal

\begin{lem} If $b \ge 3,$ then $\mu^0 \ne \l_i^0\,(2\le i\le 5)$ or its dual.
\end{lem}

\pf  Suppose $b \ge 3$ and $\mu^0 = \l_i^0$ with $2\le i\le 5$.  Then $V_Y(\l) \downarrow X = \wedge^i(b0\cdots0b)$ and 
 \cite{MF} shows that this is not MF. And Theorem 15.1 of \cite{MF} rules out the dual. \hal

At this point the remaining possibilities are $\mu^0 = \l_1^0$ or  $\mu^0 = \l_2^0$ (for $b = 2$). But then
$\l = \l_1$ or $\l_2$, completing the proof of Proposition \ref{b>1}.

\section{The case $X = A_{2s}$, $\d = \om_s+\om_{s+1}$}
 
 In this section we assume that $X=A_{2s}$ and $\d =  \om_s+\om_{s+1}$. Here $Y$ is orthogonal, by Lemma \ref{setup}. The case $s=1$ is covered by Theorem \ref{b0...0b}, so assume $s\ge 2$. 
In our proof we will adopt what we call the Inductive Hypothesis -- that is, the hypothesis that the conclusion of Theorem \ref{mainthm} holds for groups $X$ of type $A$ and rank less than $2s$.

  \begin{thm}\label{oms+om(s+1)} Let $X=A_{2s}$ with $s\ge 2$, $\d = \o_s+\o_{s+1}$ and $X<Y=SO(W)$ with $W = V_X(\d)$. Assume the Inductive Hypothesis and let $\l \ne 0,\l_1$ be a dominant weight for $Y$.  Then $V_Y(\l) \downarrow X$ is MF if and only if $\l = \l_2$ or $2\l_1$.
  \end{thm}

The fact that $V_Y(\l) \downarrow X$ is MF for the weights $\l$ in the conclusion follows directly from \cite{MF}. 

Now assume the hypotheses of the theorem, so that $X = A_{2s}$ with $s\ge 2$ and $\d = \o_s+\o_{s+1}$. 
Then   $L_Y' = C^0 \times C^1$, and   $L_X' = A_{2s-1}$ is embedded in $C^1$ via $2\om_s \oplus (\om_{s-1}+\om_{s+1})$. 
  Notice that each summand here is self-dual and so the image of $L_X'$ is contained in a product of two orthogonal subgroups of $C^1$.
  
   We will prove the result in a series of lemmas.  Assume  $V_Y(\l) \downarrow X$ is MF.    We first consider the possibilities
 for $\mu^0$ and   $\mu^1$ applying  \cite{MF}  and the Inductive Hypothesis to the two weights.
 
 \begin{lem}\label{possmu0mu1} {\rm  (i)} $\mu^0$ or its dual is in $\{ 0, \l_1^0, \l_2^0, 2\l_1^0 \}$.
 
 {\rm (ii)} $\mu^1 \in \{ 0,\, \l_1^1, \,\l_2^1, \l_3^1\,(s=2),\, 2\l_1^1,\, 3\l_1^1\,(s=2) \}$.
  \end{lem}

 \pf Part (i) is immediate from Theorem 1 of \cite{MF}.  Now consider $\mu^1$ and assume
 the image of $L_X'$ in $C^1$ is contained in $R\times S$ where these are orthogonal subgroups corresponding to the
 representations $2\om_s, (\om_{s-1}+\om_{s+1})$, respectively.  The factor $S$ could have type $B$ or $D$, and the natural orthogonal module for $S$ has dimension at least 15.    Embed $S$ in an orthogonal Levi subgroup of $C^1$, say $T$, of either the same type as $S$ or possibly of rank one greater.
 
 First assume that $S$ and $T$ are of the same type with base $\b_x^1, \ldots, \b_y^1$, $x < y$.  This group is conjugate to an orthogonal subgroup of $C^1$ with base $\b_1^1+ \cdots +\b_x^1, \b_{x+1}^1, \ldots ,\b_y^1$.  If we now consider the restriction of  $\mu^1$ to this subgroup and use the Inductive Hypothesis, we find that $\mu^1$ is as in (ii).
 
 The remaining cases are similar.  Suppose $S = B_c$ and $T = D_{c+1}$. Taking the base $\b_x^1, \cdots, \b_y^1$ for $ T$ as above, we can use the base  $\b_x^1, \ldots, \b_{y-2}^1, \b_y^1+\b_{y-1}^1$ for $S$ and take a conjugate with base $\b_1^1+ \cdots +\b_x^1, \b_{x+1}^1,  \ldots, \b_{y-2}^1, \b_y^1+\b_{y-1}^1$.  At this point we restrict the highest weight and again get (ii).  We leave the remaining case to the reader.  \hal
 
 \begin{lem}\label{possmu1}  $\mu^1 =0$.
  \end{lem}
 
 \pf  Suppose false. We first rule out all cases of (ii) in Lemma \ref{possmu0mu1} other than  $\l_1^1$.  Note that  
$\wedge^2(2\om_s \oplus (\om_{s-1}+\om_{s+1})) \supseteq (\om_{s-1}+2\om_s+\om_{s+1})^2$
 and $S^2(2\om_s \oplus (\om_{s-1}+\om_{s+1}) )\supseteq (2\om_{s-1}+2\om_{s+1})^2$. It follows that 
if $\mu^1 = \l_2^1$ or $2\l_1^1$, then  $V^1$ is not MF, a contradiction.  
Also, $\wedge^3(2\om_2 \oplus (\om_1+\om_3)) \supseteq (\om_1+2\om_2+\om_3)^2$ and $S^3(2\om_2 \oplus (\om_1+\om_3)) \supseteq (3\om_1+3\om_2)^2$. Hence also $\mu^1 \ne \l_3^1$ or $3\l_1^1$. 

Suppose finally that $\mu^1 = \l_1^1$.  Then $V^2(Q_Y) \supseteq (\mu^0 + \l_{r_0}^0) \otimes \l_2^1$.
The restriction of $\l_2^1$ to $L_X'$  contains $(\om_{s-1}+ 2\om_s+ \om_{s+1})^2$ which has
$S$-value 4, and Lemma 3.9 of \cite{MF} shows that 
$S(\mu^0 + \l_{r_0}^0) > S(\mu^0)$.  Therefore $V^2$ has a repeated composition factor with $S$-value at least $S(\mu^0)+5$ and this is a contradiction by Lemma \ref{prop38}.  \hal

\begin{lem} $\mu^0 \ne \l_{r_0}^0, 2\l_{r_0}^0$, or $\l_{r_0-1}^0$.
\end{lem}

\pf  First assume $\mu^0 = \l_{r_0}^0$.  Then $V^1 = \om_{s-1} + \om_s$ and
$V^2 \supseteq \wedge^2(\om_{s-1} + \om_s) \otimes (2\om_s \oplus (\om_{s-1}+\om_{s+1}))$.
The first tensor factor contains $3\om_{s-1} +\om_{s+1}$ and it follows using \cite[Lemma 7.1.7]{MF} 
that $V^2 \supseteq (3\om_{s-1} +\om_{s+1})^2$. Now comparing $S$-values gives a contradiction.

Next suppose $\mu^0 = 2\l_{r_0}^0$ so that $V^2(Q_Y) \supseteq ((\l_{r_0}^0 \otimes \l_{r_0-1}^0)- \l_{r_0-2}^0) \otimes \l_1^1$.  The restriction of the first tensor factor to $L_X'$ contains
$(3\om_{s-1} + 3\om_s) - \a_s = (4\om_{s-1} +\om_s +\om_{s+1})$ and tensoring with
the restriction of the second tensor factor we see that $V^2 \supseteq (\om_{s-2}+4\om_{s-1} +\om_s +\om_{s+1}+ \om_{s+2})^2$ where the first and last terms do not exist if $s = 2$.  This has
$S$-value at least 6 whereas $S(V^1)  = 4$, again a contradiction.

Finally assume $\mu^0 = \l_{r_0-1}^0$ so that $V^1 = \wedge^2( \om_{s-1} + \om_s)$ and
$V^2 \supseteq \wedge^3(\om_{s-1} + \om_s) \otimes (2\om_s \oplus (\om_{s-1}+\om_{s+1}))$.
 The first tensor factor contains $(\om_{s-2}+2\om_{s-1} + 2\om_s+\om_{s+1})$ where the first  term does not occur if $s = 2$.  Therefore $V^2 \supseteq  (2\om_{s-2}+2\om_{s-1} + 2\om_s+\om_{s+1}+\om_{s+2})^2$ or $(221)^2$ according as $s > 2$ or $s = 2$.  In the former case we have an $S$-value contradiction since $S(V^1) = 4$.  And in the latter case
 $V^1 = \wedge^2(110) =030+010+301+111$ and this cannot contribute a summand 221 to $V^2$.
 This establishes the result.  \hal
 
 \begin{lem}\label{mu0=l1,2l1,l2} If $\mu^0 = 0, \l_1^0, 2\l_1^0$ or $\l_2^0$, then $\l = 2\l_1$ or $\l_2$.
 \end{lem}
 
 \pf   Assume  $\mu^0 = 0, \l_1^0, 2\l_1^0$ or $\l_2^0$.   Suppose $\la \l, \g_1 \ra > 0$, 
which must be the case if $\mu^0 = 0$ or $\lambda_1^0$. 
 Then Lemma 5.4.1 of \cite{MF} shows that  $V^2(Q_Y) \supseteq (\mu^0 \otimes \l_{r_0}^0) \otimes \l_1^1$  and
 restricting to $L_X'$ this is $(\mu^0 \downarrow L_X') \otimes (\om_{s-1}+\om_s) \otimes (2\om_s \oplus (\om_{s-1}+\om_{s+1}))$ and this contains $(\mu^0 \downarrow L_X') \otimes (\om_{s-2}+\om_{s-1}+\om_s+\om_{s+2})^2$ where as usual we omit  terms if $s= 2$.  Therefore
 $V^2$ has a repeated composition factor of $S$-value at least $S(\mu^0) +2$ and this is a
 contradiction by Lemma \ref{prop38}.  Therefore $\l$ is as indicated in the assertion.  \hal
 
 This completes the proof of  Theorem \ref{oms+om(s+1)}.

\section{The case $\d = \o_i+\o_{l+2-i}$}

 In this section we assume that $X=A_{l+1}$ and $\d =  \om_i+\om_{l+2-i}$. Here $Y$ is orthogonal, by Lemma \ref{setup}. The cases $i=1$ and $i=\frac{l+1}{2}$ are covered by Theorems \ref{b0...0b} and \ref{oms+om(s+1)}, so we assume that $1 < i < \frac{l+1}{2}$. 

Here is the main result of the section. As in Theorem \ref{oms+om(s+1)} we shall assume the Inductive Hypothesis in the proof.

  \begin{thm}\label{omi+istar} Let $X=A_{l+1}$, let $\d = \o_i+\o_{l+2-i}$ with $1 < i < \frac{l+1}{2}$, and let $X<Y=SO(W)$ with $W = V_X(\d)$. Assume the Inductive Hypothesis and let $\l \ne 0,\l_1$ be a dominant weight for $Y$.  Then $V_Y(\l) \downarrow X$ is MF if and only if $i=2$ and $\l = \l_2$.
  \end{thm}

The fact that $V_Y(\l) \downarrow X$ is MF for $i=2$ and $\l = \l_2$ follows directly from \cite{MF}. 

Now assume the hypotheses of the theorem, so that $X = A_{l+1}$ and $\d = \o_i+\o_{l+2-i}$ with 
$1 < i < \frac{l+1}{2}$. By Lemma \ref{setup}, $Y$ is orthogonal, $k=1$ and 
 $L_Y' = C^0 \times C^1$, with $L_X' = A_{l}$ is embedded in $C^0$ via $(\om_i + \o_{l+2-i})$ and in 
$C^1$ via $(\o_i+\o_{l+1-i}) \oplus (\o_{i-1}+\o_{l+2-i})$. 
  Note that $C^1$ is orthogonal, and the image of $L_X'$ is contained in a product of two orthogonal subgroups of $C^1$.
  
   We will prove the result in a series of lemmas.  Assume  $V_Y(\l) \downarrow X$ is MF.    We first consider the possibilities
 for $\mu^0$ and   $\mu^1$.
 
 \begin{lem}\label{possmu0mu1caseistar} {\rm  (i)} $\mu^0$ or its dual is one of the following:
\[
\begin{array}{l}
 0,\, \l_1^0, \\
 \l_2^0,\,2\l_1^0\,(\hbox{both with }i=2), \\
\l_2^0\,(i=3,l=6).
\end{array}
\]
 
 {\rm (ii)} $\mu^1$ is one of:
\[
\begin{array}{l}
 0,\, \l_1^1, \\
 \l_2^1\,(i=2 \hbox{ or }(i,l)=(3,6)), \\
2\l_1^1\,(i=2,l=4).
\end{array}
\]
\end{lem}

 \pf Part (i) follows from \cite{MF}, while (ii) follows from the Inductive Hypothesis together with the argument of Lemma 
\ref{possmu0mu1}. \hal

\begin{lem}\label{mu10star} We have $\mu^1=0$.
\end{lem}

\pf Suppose $\mu^1\ne 0$, so that it is one of the possibilities given in Lemma \ref{possmu0mu1caseistar}(ii). 

If $\mu^1 = 2\l_1^1$ with $(i,l)=(2,4)$, then $V^1$ contains $S^2(0110 + 1001)-0000$ (for $L_X'=A_4$), and this is not MF, a contradiction.

Now suppose $\mu^1=\l_2^1$ with $i=2$. Then $V^1$ contains the wedge-square of 
$(\o_2+\o_{l-1}) \oplus (\o_{1}+\o_{l})$, and this is not MF by \cite[Lemma 14.2.1]{MF}. We reach a similar contradiction if $\mu^1 = \l_2^1$ with $(i,l)=(3,6)$.

It remains to consider the case where $\mu^1=\l_1^1$. Here $S(V^1) = S(\mu^0 \downarrow L_X')+2$, and
\[
V^2 \supseteq (\mu^0+\l_{r_0})\downarrow L_X' \otimes \wedge^2((\o_i+\o_{l+1-i}) \oplus (\o_{i-1}+\o_{l+2-i})).
\]
The second tensor factor has a repeated summand of highest weight $\o_{i-1}+\o_i+\o_{l+1-i}+\o_{l+2-i}$ (also by  \cite[Lemma 14.2.1]{MF}), and also $S((\mu^0+\l_{r_0})\downarrow L_X')  \ge S(\mu^0\downarrow L_X')$ by \cite[Lemma 3.9]{MF}. It follows that  $V^2$ has a repeated summand of $S$-value at least $S(\mu^0\downarrow L_X')+4$, hence at least $S(V^1)+2$, which is a contradiction by Lemma \ref{prop38}. \hal

\begin{lem}\label{g1star} We have $\la \l, \g_1 \ra =0$.
\end{lem}

\pf Suppose $\la \l, \g_1 \ra \ne 0$. Then by \cite[Prop. 5.4.1]{MF}, we have
\[
\begin{array}{ll}
V^2 & \supseteq (\mu^0\otimes \l_{r_0}^0\otimes \l_1^1)\downarrow L_X' \\
       & \supseteq (\mu^0\downarrow L_X') \otimes (\o_{i-1}+\o_{l+1-i}) \otimes 
((\o_i+\o_{l+1-i}) \oplus (\o_{i-1}+\o_{l+2-i})).
\end{array}
\]
By \cite[Lemma 7.1.7]{MF}, $(\l_{r_0}^0\otimes \l_1^1)\downarrow L_X'$ has a repeated summand of $S$-value at least 2. Hence $V^2$ has a repeated summand of $S$-value at least $S(\mu^0\downarrow L_X')+2$, contradicting Lemma \ref{prop38}. \hal

At this point we have $\l = \mu^0$, which is given by Lemma \ref{possmu0mu1caseistar}(i).

\begin{lem}\label{dualstars} $\l$ is not $\l_{r_0}$, $2\l_{r_0}$ or $\l_{r_0-1}$.
\end{lem}

\pf Observe first that $\l_{r_0}$ and  $\l_{r_0-1}$ are not possible, as in these cases $V\downarrow X = \wedge^{r_0}W$ or $\wedge^{r_0-1}W$, which is not MF by \cite{MF}.

Now suppose  $\l = 2\l_{r_0}$, so $i=2$. Here $V^1 = S^2(\o_1+\o_{l-1})$ has $S$-value 4 and 
$V^2$ contains the restriction to $L_X'$ of $(\l_{r_0-1}^0 + \l_{r_0}^0) \otimes \l_1^1$. The restriction of the first tensor factor is $(\wedge^2\nu \otimes \nu)-\wedge^3\nu$ where $\nu = \o_1+\o_{l-1}$, and this contains a summand of highest weight $3\nu - \a_{l-1} = 30\cdots 0111$. By \cite[Lemma 7.1.4]{MF}, the tensor product of this with $\o_2+\o_{l-1}$ has a repeated summand $210\cdots 0112$ if $l\ge 5$ (and $4112$ if $l=4$). This is contained in $V^2$ and has $S$-value 
at least 6, contradicting Lemma \ref{prop38}. 
 \hal

The remaining possibilities for $\l$ given by Lemma \ref{possmu0mu1caseistar}(i) are $2\l_1$ and $\l_2$, with $i=2$ in the first case and either $i=2$ or $(i,l)=(3,6)$ in the second. If $\l=2\l_1$ then $V_Y(\l) \downarrow X = S^2(\o_2+\o_l) - 0$ (and $l\ge 4$), and this is not MF by \cite{MF}, noting that 0 only appears with multiplicity 1 in $S^2(\o_2+\o_l)$.
Finally, if $\l = \l_2$ with $(i,l)=(3,6)$ then $V_Y(\l) \downarrow X = \wedge^2(\o_3+\o_5)$ which is not MF by \cite{MF}. This leaves $\l = \l_2$ with $i=2$, as in the conclusion of Theorem \ref{omi+istar}, completing the proof of this theorem.


\section{The case $X = A_5 < Y = C_{10}$ with $\d = \o_3$}\label{caseA5o3}

In the next three sections we prove Theorem \ref{mainthm} in the case where $X = A_{2s+1}$ with $\d = \o_{s+1}$.
There are many exceptional examples in the theorem in the smallest case where $s=2$, and it is convenient to deal separately with this case, which we do in this section.

Here we prove Theorem \ref{mainthm} in the case where $X=A_5$ and $\d = \o_3$. In this case, by Lemma \ref{setup}  we have $Y = Sp(W) = C_{10}$, where $W = V_X(\d)$. We prove 

\begin{thm}\label{A5<C10} Assume $X = A_5$ is embedded in $Y = C_{10}$ via $\d = \o_3$,  and let $V = V_Y(\l)$
where $l$ is a dominant weight for $Y$ and 
$\l \ne 0,\l_1.$  Then $V \downarrow X$ is MF if and only if $\l$ is one of the following:
\[
\begin{array}{l}
\l_i, \\
c\l_1\,(c\le 5), \\
\l_1+\l_i\,(i=2,3,4,10), \\
2\l_2, \\
2\l_1+\l_2.
\end{array}
\]
\end{thm}

We begin by showing that $V_Y(\l) \downarrow X$ is indeed MF for each weight $\l$ listed in the theorem.

\begin{lem}\label{MFA5o3} If $X = A_5 < Y=C_{10}$ via $\d=\o_3$, and $\l$ is one of the weights listed in the conclusion of Theorem $\ref{A5<C10}$, then  $V_Y(\l) \downarrow X$ is MF.
\end{lem}

\pf This follows from Magma computations, using the following constructions of the modules $V_Y(\l)$ based on taking tensor products of symmetric and exterior powers of the natural module $W = V_Y(\l_1)$; these are consequences of Lemmas \ref{tensorl1} and \ref{symwed} apart from the last two cases, which can be verified using Magma:
\[
\begin{array}{l}
V_Y(\l_i) = \wedge^iW - \wedge^{i-2}W, \\
V_Y(c\l_1) = S^cW, \\
V_Y(\l_1+\l_2) = W\otimes \wedge^2W - \wedge^3W-W,\\
V_Y(\l_1+\l_3) = W\otimes \wedge^3W - W\otimes W + 0 - \wedge^4W, \\
V_Y(\l_1+\l_4) = W\otimes \wedge^4W - W\otimes \wedge^2W +W - \wedge^5W, \\
V_Y(\l_1+\l_{10}) = W\otimes \wedge^{10}W - W\otimes \wedge^8W +\wedge^7W - \wedge^9W, \\
V_Y(2\l_2) = S^2(\wedge^2W)-\wedge^4W-\wedge^2W,\\
V_Y(2\l_1+\l_2) = W \otimes S^3W - S^4W - S^2W. 
\end{array}
\]
 \hal

We now embark on the proof of the other implication in the statement of Theorem \ref{A5<C10}. Let $X,Y$ be as in the hypothesis, and suppose $\l$ is a dominant weight for $Y$ such that  $V_Y(\l) \downarrow X$ is MF. 

We have $L_X' = A_4$ embedded in $L_Y' = A_9$ via the representation $\d' = 0010.$  It will be convenient to let $\d'' = 0100$.  The main result of \cite{MF} shows that $\mu^0$ or its dual is one of  the following: 
\begin{equation}\label{possa4}
\begin{array}{l}
0, \\
a\l_1^0, \\
a\l_1^0+\l_2^0, \\
a\l_1^0+\l_9^0, \\
\l_i^0, \\
\l_1^0+\l_i^0, \\
a\l_2^0\,(a\le 5), \\
2\l_3^0,\,2\l_4^0, \\
\l_1^0+2\l_2^0,\,\l_2^0+\l_3^0,\,\l_2^0+\l_8^0.
\end{array}
\end{equation}

At this point we work through the various possibilities for $\mu^0$ in (\ref{possa4}) doing the easiest cases first. Note that if $\mu^0 = 0$ then $\l = a\l_{10}$, which is covered by Lemma \ref{al10}. As usual we simplify notation by just writing the highest weight $\mu$ of an irreducible module rather than $V(\mu)$.


\begin{lem}\label{l2+l8} We have $\mu^0 \ne \l_2^0 +\l_8^0.$  
\end{lem}

\pf  Assume that  $\mu^0  = \l_2^0 +\l_8^0,$ so that $\l =  \l_2 +\l_8 + x\l_{10}$.   We first consider the case where   $x = 0$.  Magma shows that the $C_{10}$-module $\l_2 \otimes \l_8 = (\l_2 + \l_8) \oplus (\l_1+\l_9) \oplus (\l_1+\l_7) \oplus \l_{10} \oplus \l_8 \oplus \l_6.$  Using Lemma \ref{tensorl1} we see that $\l_2+\l_8 = (\l_2 \otimes \l_8) -(\l_1 \otimes \l_9) -(\l_1 \otimes \l_7) + \l_8.$ Now restrict to $X = A_5$ and recall that $\l_i \downarrow X =
\wedge^i(\d) - \wedge^{i-2}(\d)$ for $i \ge 2.$ It follows that
\[ 
\begin{array}{ll}
(\l_2 + \l_8)\downarrow X = & \left((\wedge^2(\d)-0) \otimes (\wedge^8(\d)-\wedge^6(\d))\right) -
\left(\d\otimes(\wedge^9(\d)-\wedge^7(\d))\right)\\
  & -\left(\d\otimes(\wedge^7(\d)-\wedge^5(\d))\right)+(\wedge^8(\d)-\wedge^6(\d)),
\end{array}
\]
which simplifies to
\[
(\l_2 + \l_8)\downarrow X = \left(\wedge^2(\d) \otimes (\wedge^8(\d)-\wedge^6(\d))\right) -
\left(\d \otimes(\wedge^9(\d)-\wedge^5(\d))\right).
\]
A Magma check shows that this contains $(21012)^6$, which establishes the result when $x = 0.$

Now suppose $x \ne 0.$ Then $V^2(Q_Y) \supseteq (\l_1^0+\l_8^0+\l_9^0) \oplus (\l_2^0+\l_7^0+\l_9^0) \oplus (\l_2^0+\l_8^0+2\l_9^0),$ where these summands are afforded by $\l - \b_2 - \cdots -\b_{10}$,  $\l - \b_8 - \b_9 -\b_{10}$, $\l - \b_{10}$, respectively.  Therefore $V^2(Q_Y) \supseteq ((\l_2^0+\l_8^0) \otimes 2\l_9^0) - (\l_1^0+\l_7^0)$ and restricting to $L_X'$ we have $V^1 \otimes (0200 + 0001) - 0010 \otimes \wedge^3(0100)+ \wedge^2(0100)$.  By Corollary 5.1.5 of \cite{MF}, $V^1 \otimes 0001$ covers all composition factors arising from $V^1$, so it suffices to show that $(V^1 \otimes 0200)  - (0010 \otimes \wedge^3(0100))+ \wedge^2(0100)$ is not MF.

Now $V^1 = \wedge^2(0010) \otimes \wedge^8(0010) - (0010 \otimes \wedge^9(0010)) $ and we use Magma to see that $V^1 \otimes 0200 \supseteq (1201)^8$
and $1201$ does not occur in $0010 \otimes \wedge^3(0100)$.  This completes the proof of the lemma. \hal

 \begin{lem}\label{2l3} $\mu^0 \ne 2\l_3^0$ or its dual.
\end{lem} 

\pf  First assume that $\l = 2\l_3$.  Then $V^2(Q_Y) \supseteq (\l_2^0+\l_3^0+\l_9^0)  \oplus 2\l_2^0$ where the summands are afforded by $\l- \b_ 3 -\cdots -\b_{10}$ and $\l- 2\b_3 - \cdots -2\b_9 -\b_{10}$, respectively.
To see that the second summand occurs we note that this weight occurs with multiplicity 8 for $C_{10}$ but only multiplicity $7$ in the first summand.  It follows that $V^2(Q_Y) \supseteq ((\l_2^0+\l_3^0) \otimes \l_9^0)   - (\l_1^0+\l_3^0).$  Then a Magma check shows that $V^2 \supseteq (1020)^4$. Now $1020$ can only arise from composition factors in $V^1$ of highest weights $2020$, $0120$, and $1011$ (see \cite[Cor. 5.1.5]{MF}). On the other hand $V^1(Q_Y) = S^2(\l_3^0) - (\l_1^0+\l_5^0)$ so that $V^1 = S^2(\wedge^3(0010)) - (0010 \otimes \wedge^5(0010)) + \wedge^6(0010)$ and a Magma check shows that
$V^1$ can only contribute $(1020)^2$ to $V^2$ since of the above composition factors it only contains $2020$ and $1011$. Therefore $V \downarrow X$ is not MF in this case, a contradiction.

Next assume $\l  = 2\l_3 + x\l_{10}$ for $x >0$.  Here we see that $V^2(Q_Y) \supseteq (\l_2^0+\l_3^0+\l_9^0) \oplus (2\l_3^0+2\l_9^0)$ where the summands are afforded by $\l- \b_ 3 -\cdots -\b_{10}$ and $\l-\b_{10}$, respectively.  Therefore $V^2(Q_Y) \supseteq (2\l_3^0 \otimes 2\l_9^0) - 2\l_2^0$.  Restricting to $L_X'$ we see that $V^2 \supseteq V^1 \otimes (0001 + 0200) - (2\l_2^0 \downarrow L_X')$ and so by \cite[5.1.5]{MF}, it suffices to show that $(V^1 \otimes 0200)) - (2\l_2^0 \downarrow L_X')$ is not MF.  Now $ 2\l_2^0 \downarrow L_X' = S^2(\wedge^2(0010)) - \wedge^4(0010) = 0011 + 1020+ 0202 + 2000$.
However, from the end of the last paragraph we see that $V^1 \otimes 0200 \supseteq (2020 + 1011) \otimes (0200)$ and
this contains $(1101)^3$.  Therefore, $V\downarrow X$ is not MF in this case.

Now we must deal with the duals.  First assume $\l = 2\l_7$.   Arguing as in the first paragraph we see that
$V^2(Q_Y) \supseteq  (\l_6^0+\l_7^0+\l_9^0) \oplus 2\l_6^0$ where the summands are afforded by $\l - \b_7 - \b_8 - \b_9 - \b_{10}$
and $\l - 2\b_7 - 2\b_8 - 2\b_9 - \b_{10},$ respectively.  Therefore $V^2(Q_Y) \supseteq ((\l_6^0+\l_7^0) \otimes \l_9^0) - (\l_5^0+\l_7^0)$
and restricting to $L_X'$ we see that 
$$V^2 \supseteq \wedge^4(\d'') \otimes \wedge^3(\d'')\otimes \d'' - \wedge^5(\d'') \otimes \wedge^2(\d'') \otimes \d'' - \wedge^3(\d'') \otimes \wedge^5(\d'') + \wedge^2(\d'') \otimes \wedge^6(\d'')$$ 
and a Magma check shows that this contains $(2111)^6$.  On the other hand
$V^1 \subseteq  S^2(\wedge^3(0100))$ and another check shows that this can only contribute 3 copies of 2111 to $V^2$.
Hence we again see that $V\downarrow X$ is not MF.

Finally, assume $\l  = 2\l_7 + x\l_{10}$ for $x >0$.   Here we find that $V^2(Q_Y) \supseteq (\l_6^0+\l_7^0+\l_9^0) \oplus  (2\l_7^0 +2\l_9^0) =
(2\l_7^0  \otimes 2\l_9^0) - 2\l_6^0$ and restricting to $L_X'$ this is $V^1 \otimes (0001 + 0200) - (2\l_6^0\downarrow L_X')$.  Therefore it suffices to show that 
$(V^1 \otimes 0200) - (2\l_6^0\downarrow L_X')$
is  not MF.  By taking duals of the $V^1$ analysis of the first paragraph we see that $V^1 \otimes 0200 \supseteq (1101+0202) \otimes 0200$ and this contains $(1110)^3$.  On the other hand $2\l_6^0 \downarrow L_X'$ is MF so the
difference is not MF.  This completes the proof of the lemma.  \hal

 \begin{lem}\label{2l4} $\mu^0  \ne 2\l_4^0$ or its dual.   
\end{lem} 

\pf  The analysis here is very similar to that of the previous lemma. First assume that $\l = 2\l_4.$ This time we argue that $V^2 \supseteq (\l_3^0+\l_4^0+\l_9^0) \oplus 2\l_3^0,$  where  the summands are afforded by $\l - \b_4 - \cdots - \b_{10}$ and 
by $\l - 2\b_4 - \cdots - 2\b_9 - \b_{10}$, respectively.   Therefore
$V^2(Q_Y) \supseteq ((\l_3^0+\l_4^0) \otimes \l_9^0) - (\l_2^0+\l_4^0)$   and  so 
$$V^2 \supseteq ((\wedge^3(\d') \otimes \wedge^4(\d')) - (\wedge^2(\d') \otimes \wedge^5(\d'))) \otimes \d'' - ((\wedge^2(\d') \otimes \wedge^4(\d')) + (\wedge^1(\d')) \otimes \wedge^5(\d'))$$
 and a Magma check shows that this contains $(0112)^6$.
Now $V^1(Q_Y) = 2\l_4^0 = S^2(\l_4^0) - (\l_2^0+\l_6^0) - \l_8^0$ and another Magma check shows
that this contributes only 3 copies of $0112$ to $V^2$.   Therefore $V \downarrow X$ is not MF.

Next assume $\l = 2\l_4 + x\l_{10}$ for $x >0$.  Here we see that $V^2(Q_Y) \supseteq (\l_3^0+\l_4^0+\l_9^0) \oplus (2\l_4^0 + 2\l_9^0),$ where the summands are afforded by $\l- \b_ 4 -\cdots -\b_{10}$ and $\l-\b_{10}$, respectively.   Therefore $V^2(Q_Y) \supseteq (2\l_4^0 \otimes 2\l_9^0) - 2\l_3^0$. Restricting to $L_X'$ we see that $V^2 \supseteq V^1 \otimes (0001 + 0200) - (2\l_3^0 \downarrow L_X')$ and it suffices to show that $(V^1 \otimes 0200)) - (2\l_3^0 \downarrow L_X')$ is not MF.  A Magma check shows that $V^1 \supseteq 1112+1104+2202$ and another check shows that $V^1 \otimes 0200 \supseteq (1113)^5$.  On the other hand $(2\l_3^0 \downarrow L_X')$ is MF and we conclude that $V\downarrow X$ is not MF.

Now assume $\l = 2\l_6$.  Arguing as above we see that $V^2(Q_Y) \supseteq (\l_5^0+\l_6^0+\l_9^0) \oplus 2\l_5^0 =
((\l_5^0+\l_6^0) \otimes \l_9^0) - (\l_4^0+\l_6^0)$. Therefore 
$$V^2 \supseteq \wedge^5(\d') \otimes \wedge^6(\d')\otimes \d'' - \wedge^4(\d') \otimes \wedge^7(\d') \otimes \d'' - \wedge^4(\d') \otimes \wedge^6(\d') + \wedge^3(\d')) \otimes \wedge^7(\d').$$  
 A Magma check shows that this contains $(1111)^{10}$.  On the other hand $V^1$ 
can be computed as the dual of what appears in the first paragraph and we find that $V^1$ can only contribute
$3$ copies of $1111$ to $V^2.$ Hence here too we find that  $V \downarrow X$ is not MF.

Finally assume $\l = 2\l_6 + x\l_{10}$ for $x >0$. Here we see that $V^2(Q_Y) \supseteq (\l_5^0+\l_6^0+\l_9^0) \oplus (2\l_6^0+2\l_9^0)$ where the summands are afforded by $\l- \b_ 6 -\cdots -\b_{10}$ and $\l-\b_{10}$, respectively.   Therefore $V^2(Q_Y) \supseteq (2\l_6^0 \otimes 2\l_9^0) - 2\l_5^0 $. Restricting to $L_X'$ we see that $V^2 \supseteq V^1 \otimes (0001 + 0200) - (2\l_5^0 \downarrow L_X')$ and it suffices to show that $(V^1 \otimes 0200) - (2\l_5^0) \downarrow L_X'$ is not MF.  Taking duals of the weights in the second paragraph we see that $V^1 \otimes 0200 \supseteq (2111+4011+2022) \otimes (0200)$ and this contains $(2112)^5$.  On the other hand $2\l_5^0 \downarrow L_X' \subset
S^2(\wedge^5(0010))$  which only contains one copy of $2112$.  This completes the proof of the lemma. \hal

\begin{lem}\label{l2+l3} $\mu^0  \ne \l_2^0 +\l_3^0$ or its dual.
\end{lem}

\pf  Assume $\l = \l_2 +\l_3$.  A Magma check shows that  the $C_{10}$-module $\l_2 \otimes \l_3 = (\l_2+\l_3) \oplus (\l_1+\l_2) \oplus (\l_1+\l_4) \oplus \l_1 \oplus \l_3 \oplus \l_5$ and Lemma \ref{tensorl1} shows that $\l_1 \otimes \l_4 = (\l_1+\l_4) \oplus \l_3\oplus\l_5$ and $\l_1 \otimes \l_2 = (\l_1+\l_2) \oplus \l_1\oplus \l_3$.  Therefore $\l_2+\l_3 = (\l_2 \otimes \l_3) - (\l_1 \otimes \l_4) - (\l_1 \otimes \l_2) + \l_3$. Restricting
this to $X = A_5$ we have 
$$(\wedge^2(\d) - 0) \otimes (\wedge^3(\d) - \d) - (\d \otimes (\wedge^4(\d) - \wedge^2(\d))) - (\d \otimes (\wedge^2(\d)-0)) + (\wedge^3(\d)-\d).$$
A Magma check shows that this contains $(10101)^2$ and hence $V\downarrow X$ is not MF.

Next suppose $\l = \l_2 +\l_3+x\l_{10}$ with $x \ne 0$.  Then $V^2(Q_Y) \supseteq (\l_1^0+\l_3^0+\l_9^0) \oplus (2\l_2^0+\l_9^0)\oplus (\l_2^0+\l_3^0+2\l_9^0)
= ((\l_2^0+\l_3^0) \otimes2\l_9^0) - (\l_1^0+\l_2^0)$. Therefore $V^2 \supseteq (V^1 \otimes (0001)) + (V^1 \otimes 0200) - (0010 \otimes \wedge^2(0010) - \wedge^3(0010))$.  By Cor. 5.1.5 of \cite{MF} the  first term covers the contribution of $V^1$ to $V^2$ so it suffices
to show that the remaining terms fail to be MF. 
Now $V^1 = (\wedge^2(0010) \otimes \wedge^3(0010)) - (0010 \otimes \wedge^4(0010))$, and using this we see
that the remaining terms contain $(1112)^7$.  Therefore the result holds in this case.

It remains to consider the cases where $\l = \l_7+\l_8 + x\l_{10}$.  First assume $x \ne 0$. Then
$V^2(Q_Y) \supseteq (\l_6^0+\l_8^0+\l_9^0) \oplus (2\l_7^0+\l_9^0)\oplus (\l_7^0+\l_8^0+2\l_9^0) = ((\l_7^0+\l_8^0) \otimes 2\l_9^0) - (\l_6^0+\l_7^0)$.
Restricting to $L_X'$ this becomes $V^1 \otimes (0001 + 0200) - (\wedge^6(0010) \otimes \wedge^7(0010)) + (\wedge^5(0010) \otimes \wedge^8(0010))$. By Corollary 5.1.5 of \cite{MF}, $V^1 \otimes (0001)$ covers the contribution of $V^1$ to $V^2$
and $V^1 = (\wedge^7(0010) \otimes \wedge^8(0010)) - (\wedge^6(0010) \otimes \wedge^9(0010))$.  We then
check that the terms other than $V^1 \otimes 0001$ contain $(1112)^4$ and hence $V\downarrow X$ is not MF.

The final case here is where $\l = \l_7+\l_8$.  A Magma check shows that the $C_{10}$-module
\[
\begin{array}{ll}
\l_7+\l_8 = & (\l_7 \otimes \l_8) - (\l_6 \otimes \l_9) -(\l_6 + \l_7) -(\l_4+\l_5) \\ 
  & - (\l_3+\l_4) - (\l_5+\l_6) - (\l_2+\l_3) - (\l_1+\l_2) - \l_1.
\end{array}
\]
Consider the multiplicity of $21223$ in $V \downarrow X.$
Now $(\l_7 \otimes \l_8) \downarrow X = (\wedge^7(\d)- \wedge^5(\d)) \otimes (\wedge^8(\d)- \wedge^6(\d)).$
A Magma check shows that this contains $(21223)^9$.  Similarly we see that the multiplicity of
$21223$ in $(\l_6 \otimes \l_9) \downarrow X$  is $4$.  For the remaining terms we note that
$\l_6+\l_7 \subseteq  \l_6\otimes \l_7$, $\l_4+\l_5 \subseteq  \l_4\otimes \l_5$ and so on, and we conclude
that the subtracted terms contribute at most $1$ copy of $21223.$  It follows that $V \downarrow X$ is not MF.  \hal

\begin{lem}\label{l1+2l2} $\mu^0  \ne \l_1^0 +2\l_2^0$ or its dual. 
\end{lem}

\pf First assume that $\l = \l_1 + 2\l_2$.  Using Magma for $C_{10}$ gives  
$$\l_1 + 2\l_2 = (\l_1 \otimes 2\l_2) -  (\l_2+\l_3) - (\l_1+ \l_2).$$
Also $2\l_2 = S^2(\l_2) - (\l_4 + \l_2 + 0) = S^2(\l_2) - \wedge^4(\l_1)$ and $(\l_2+\l_3)\oplus (\l_1+\l_2) = \l_2\otimes \l_3 - \l_1\otimes \l_4 -\l_1$. Putting all this together and using Magma, we compute that $(\l_1+2\l_2)\downarrow X \supseteq (10101)^3$, a contradiction. 

Next suppose $\l =  \l_1 + 2\l_2 + x\l_{10}$ with $x \ne 0$. Then 
$$V^2(Q_Y) \supseteq (2\l_1^0+\l_2^0+\l_9^0) \oplus (2\l_2^0+\l_9^0) \oplus (\l_1^0+2\l_2^0+2\l_9^0)$$
$$ =((\l_1^0+2\l_2^0) \otimes 2\l_9^0) - (\l_1^0+\l_2^0) - 3\l_1^0.$$
Restricting this to $A_4$ we obtain $V^1 \otimes (0001 + 0200) - (\d' \otimes \wedge^2(\d' )) + \wedge^3(\d' ) - S^3(\d' )$. By Corollary 5.1.5 of \cite{MF},  $V^1 \otimes 0001$ covers all composition factors in $V^2$ arising from $V^1$.  So it will suffice
to show that $V^1 \otimes 0200 - (\d'  \otimes \wedge^2(\d' )) + \wedge^3(\d' ) - S^3(\d' )$ is not MF. 

Towards this end note that as $A_9$-modules,  $\l_1^0 + 2\l_2^0 = (\lambda_1^0\otimes 2\l_2^0) - (\l_2^0+\l_3^0)$ and $2\l_2^0 = S^2(\l_2^0) - \wedge^4(\l_1^0)$.  Therefore
$$V^1 = (\d'  \otimes (S^2(\wedge^2(\d' )) - \wedge^4(\d'))) - (\wedge^2(\d' ) \otimes \wedge^3(\d' )) + (\d'  \otimes \wedge^4(\d' )).$$
Using Magma we see that $V^1 \otimes 0200$ contains $(1112)^7$.  Moreover  neither of  $(\d'  \otimes \wedge^2(\d' ))$  and  $S^3(\d' )$ contains $1112$, so this shows that $V\downarrow X$ fails to be MF.

To complete the proof of the lemma we must consider the case $\l = 2\l_8 +\l_9 + x\l_{10}$.  First assume $x = 0$.  We consider certain dominant weights for $L_Y'=A_9$ as follows:

\vspace{2mm}
$$\begin{tabular}{|l|l|l|}
\hline
$(1)\ \  000000031$ & $\l-\b_9-\b_{10}$ & $1$  \\
\hline
$(2)\ \ 000000112$ & $\l-\b_8-\b_9-\b_{10}$ & $2$   \\
 \hline
$(3)\ \ 000000120$ &$\l-\b_8-2\b_9-\b_{10}$ & $4$  \\
\hline
$(4)\ \ 000000201$ & $\l-2\b_8-2\b_9-\b_{10}$ & $6$   \\
\hline
\end{tabular}$$
\vspace{2mm}

\noindent In this table the third column gives the multiplicity of the weight in $C_{10}$.  These multiplicities are clear
for the first two weights listed and one can use Magma for the others.  We claim that each of these weights
is the highest weight of a composition factor of $V^2(Q_Y)$ appearing with multiplicity $1$.  The claim is immediate
for $(1)$  since this is the largest dominant weight in $V^2(Q_Y)$.  Weight  $(2)$ is the next largest dominant weight and as this weight appears with multiplicity $1$ in the irreducible module of highest weight $(1)$, we have a unique composition factor with this highest weight.  Weight $(3)$ occurs
with multiplicity $2$ and $1$ in the irreducibles with highest weight  $(1)$ and $(2)$, respectively.  So the claim holds here.  Finally, the dominant weights above the weight $(4)$ are $(1)$, $(2)$ and $(3)$, while $(4)$ appears in these irreducibles with multiplicity $2,2,1$, respectively.  Therefore the claim holds for $(4)$.  

Next we note that $(2\l_8^0+\l_9^0)\otimes 2\l_9^0  =(2\l_8^0+3\l_9^0)\oplus  (3\l_8^0+\l_9^0)\oplus (\l_7^0 + \l_8^0+2\l_9^0) \oplus (\l_7^0 + 2\l_8^0) \oplus (2\l_7^0 + \l_9^0) $, so that
$V^2(Q_Y) \supseteq ((2\l_8^0+\l_9^0)\otimes 2\l_9^0) - (2\l_8^0+3\l_9^0)$.
Restricting to $L_X'$ this becomes $(V^1 \otimes (0001 + 0200)) - (2\l_8^0+3\l_9^0) \downarrow L_X'$ so, as in other lemmas, it will suffice to show that $(V^1 \otimes 0200) - (2\l_8^0+3\l_9^0) \downarrow L_X'$ is not MF. Now 
$$V^1 = (0100 \otimes S^2(\wedge^2(0100))) - (\wedge^2(0100) \otimes \wedge^3(0100))$$
and  also  
$$2\l_8^0+3\l_9^0  = (5\l_9^0 \otimes 2\l_9^0) - (6\l_9^0 \otimes \l_9^0)$$
A Magma check shows that $(V^1 \otimes 0200) - (2\l_8^0+3\l_9^0) \downarrow L_X'$ contains $(1112)^4$ and we have the result here.

Finally assume $\l = 2\l_8 +\l_9 + x\l_{10}$ with $x > 0$. This is very similar to the above  case only easier.  As $x \ne 0$, $\l -\b_{10}$ affords an irreducible of highest weight $2\l_8^0+3\l_9^0$ which we label as weight $(0)$.  We then check that weights $(1)$, $(2)$, $(3)$ and $(4)$ occur with multiplicities $2$, $4$, $6$, and $9$ respectively.  We now argue as above that $V^2(Q_Y)$ contains composition factors of highest weights $(0)$, $(1)$, $(2)$,  $(3)$, and $(4)$ each with multiplicity $1$.  Moreover these sum
to $000000021 \otimes 000000002$ and restricting to $L_X'$ this becomes $(V^1 \otimes 0001) + (V^1 \otimes 0200)$. By Cor. 5.1.5 of \cite{MF} the first summand covers the composition factors arising from $V^1$ and the argument above shows that the second summand is not MF.  This completes the proof of the lemma.  \hal

 \begin{lem}\label{al1+l9} $\mu^0  \ne a\l_1^0+\l_9^0\,(a>0)$ or its dual.
\end{lem} 

\pf First assume $\l = a\l_1+\l_9+x\l_{10}$.  Then 
$V^2(Q_Y) \supseteq ((a-1)\l_1^0+2\l_9^0) \oplus (a\l_1^0+\l_8^0 +\l_9^0) =
((a\l_1^0+\l_9^0) \otimes 2\l_9^0) - (a\l_1^0+3\l_9^0) - ((a-1)\l_1^0+\l_8^0) - ((a-2)\l_1^0+\l_9^0)$, where
the last term is only present if $a \ge 2$.  The usual argument shows that the assertion holds
provided we can show that $D = (V^1 \otimes 0200) - ((a\l_1^0+3\l_9^0) + ((a-1)\l_1^0+\l_8^0) + ((a-2)\l_1^0+\l_9^0)) \downarrow L_X'$ is not MF.

Now $V^1 = S^a(0010) \otimes 0100 - S^{a-1}(0010) $ and so 
\[
\begin{array}{ll}
V^1 \otimes 0200  & = (S^a(0010) \otimes 0100 - S^{a-1}(0010)) \otimes 0200 \\
  & = S^a(0010) \otimes (0300 + 1110 + 0101) - S^{a-1}(0010) \otimes 0200.
\end{array}
\]
Also $ (a\l_1^0+3\l_9^0) \downarrow L_X' = S^a(0010) \otimes S^3(0100) -  (((a-1)\l_1^0+2\l_9^0) + ((a-2)\l_1^0+\l_9^0) + (a-3)\l_1^0)\downarrow L_X' $ where the last two terms occur only if $a \ge 2, a\ge 3$, respectively.
Using this and noting that $S^a(0010) \otimes S^3(0100) = S^a(0010) \otimes  (0300+0101)$ we combine terms and obtain
\[
\begin{array}{ll}
D =&  (S^a(0010) \otimes 1110 - S^{a-1}(0010)) \otimes 0200 \, + \\
 & (((a-1)\l_1^0+2\l_9^0) + (a-3)\l_1^0 - ((a-1)\l_1^0+ \l_8^0)) \downarrow L_X' 
\end{array}
\]
where we delete one term if $a \le 2.$
Now $(a-1)\l_1^0+2\l_9^0 = S^{a-1}(\l_1^0) \otimes S^2(\l_9^0) - ((a-2)\l_1^0+\l_9^0) -(a-3)\l_1^0$,
so the above reduces to
$$ D = S^a(0010)\otimes (1110)\otimes(0200)+S^{a-1}(0010)\otimes
(0001)-((a-2)\lambda^0_1+\lambda^0_9 +(a-1)\lambda^0_1+\lambda^0_8)\downarrow
L_X'.$$
Note as well that $(a-1)\lambda^0_1\otimes \lambda^0_8 =
((a-1)\lambda^0_1+\lambda^0_8) \oplus ((a-2)\lambda^0_1+\lambda^0_9)$. Now we find that $$D = S^a(0010)\otimes (1110)\otimes(0200)+S^{a-1}(0010)\otimes (0001) - S^{a-1}(0010)\otimes\wedge^2(0100).$$
Note that $(00a0)\otimes (0200)\supseteq (11(a-1)1)$ and by \cite[7.1.5]{MF}, $(11(a+1)1)^2\subseteq (11(a-1)1\otimes 1110)$. 
Now an $S$-value argument shows that this weight does not appear as a summand in $S^{(a-1)}(0010)\otimes \wedge^2(0100)$, completing the consideration of this case.

It remains to consider the case $\l = \l_1+a\l_9+x\l_{10},$ for $a \ge 2$.  This
time $V^2(Q_Y) \supseteq (a+1)\l_9^0 \oplus (\l_1^0 + \l_8^0+a\l_9^0) = ((\l_1^0 +a\l_9^0) \otimes 2\l_9^0) - (\l_1^0 +(a+2)\l_9^0) - (\l_8^0 +(a-1)\l_9^0)- (\l_1^0+2\l_8^0+(a-2)\l_9^0)$.  Restricting to $L_X'$ we see that it suffices to show
that  the following module is not MF:
\[
\begin{array}{l}
\left((S^a(0100) \otimes 0010 - S^{a-1}(0010)) \otimes 0200\right) - \left(S^{a+2}(0100) \otimes 0010\right) + \\
S^{a+1}(0100) - (\l_8^0 +(a-1)\l_9^0)\downarrow L_X' - (\l_1^0+2\l_8^0+(a-2)\l_9^0)\downarrow L_X'.
\end{array}
\]
 Now $S^a(0100) \otimes 0010  \otimes 0200 \supseteq 0a00 \otimes (0210 + 1101)\supset
(2(a-1)11)^2$ and this has $S$-value $a+3$.  The only subtracted term with this $S$-value is contained in 
$S^{a+2}(0100) \otimes 0010 $ which by \cite[3.8.1]{Howe} is equal to $(0(a+2)00 + 0a01 +\cdots) \otimes 0010$. This does not 
contain $(2(a-1)11)$ (see \cite[Cor. 4.1.3]{MF}).  Therefore $V\downarrow X$ is not MF, completing the proof of the lemma.  \hal

 \begin{lem}\label{al1} Assume $\mu^0  = a\l_1^0\,(a>0)$ or its dual.   Then $\l = a\l_1\,(a\le 5)$, $\l_9$, 
 $\l_1 + \l_{10}$ or $\l_9+\l_{10}$.
\end{lem}

Note that all the possibilities for $\l$ except $\l_9+\l_{10}$ appear in the conclusion of Theorem \ref{A5<C10}. The latter possibility will be excluded in Lemma \ref{li813}.

\vspace{2mm}
\pf First assume $\l = a\l_1$. Then $V = S^a(\l_1)$ and it follows from \cite{MF} together with Lemma \ref{symwed} that $a \le 5$. 
 Now assume that $\l = a\l_1 + x\l_{10}$, with
$x > 0$.  Then $V^2(Q_Y) \supseteq ((a-1)\l_1^0+\l_9^0) \oplus (a\l_1^0+2\l_9^0) = (a\l_1^0 \otimes 2\l_9^0)- (a-2)\l_1^0$, where the last term does not appear if $a < 2$.
Restricting to $L_X'$ this becomes $(V^1 \otimes 0001) + (V^1 \otimes 0200) - S^{a-2}(0010)$.  By Cor. 5.1.5 of  \cite{MF} the first summand
covers all composition factors of $V^2$ arising from $V^1$, so we consider the second summand.
If $a \ge 2$ this contains $(00a0 + 10(a-2)0) \otimes 0200$ (see \cite[3.8.1]{Howe}), and using Lemma 7.1.2 of \cite{MF} we see that this contains $(01(a-1)0)^2$. A consideration of $S$-values shows that  $S^{a-2}(0010)$ contains no such composition factor, so  we conclude that
$V\downarrow X$ is not MF if $a\ge 2$.

Now suppose $a = 1$.  If $x = 1$ then $\l=\l_1+\l_{10}$ is as in the conclusion, so suppose $x \ge 2$.  Here we
check that $V^2(Q_Y) = \l_9^0 \oplus (\l_1^0+2\l_9^0) = \l_1^0 \otimes 2\l_9^0$.  Then restricting to $L_X'$ we have $V^2 = 0010 \otimes (0001 + 0200) = (0100)^2 + 0210+0011+1101$.
Next we note that $V^3(Q_Y) \supseteq (\l_1^0+4\l_9^0) \oplus (\l_1^0+\l_7^0+\l_9^0) $ where the summands are afforded by
$\l - 2\b_{10}$ and $\l - \b_8-2\b_9-2\b_{10}$, respectively.  Indeed it is clear that the first summand appears and one checks that the second weight has multiplicity $2$ in the $C_{10}$-module of highest weight $\l$ but only multiplicity $1$ in $\l_1^0+4\l_9^0.$

We will show that $V\downarrow X$ is not MF by considering the multiplicity of $1110$ in $V^3$.
Now $(\l_1^0+4\l_9^0) \downarrow L_X' = (0010 \otimes S^4(0100)) - S^3(0100)$ and Magma shows
that this contains one copy of $1110$.  The next summand is more complicated.
We have $(\l_1^0+\l_7^0+\l_9^0) = (\l_1^0 \otimes (\l_7^0+\l_9^0)) - (\l_8^0+\l_9^0) - \l_7^0$. In addition,
$(\l_7^0+\l_9^0) = (\l_7^0 \otimes \l_9^0) - \l_6^0$ and $(\l_8^0 +\l_9^0)= (\l_8^0 \otimes \l_9^0) - \l_7^0$.  Therefore the contribution of the second summand to $V^3$ is
$$(0010 \otimes ((\wedge^3(0100) \otimes 0100) - \wedge^4(0100)) -(\wedge^2(0100) \otimes 0100)$$
and we find that this contains $(1110)^3$.  Therefore $V^3 \supseteq (1110)^4$.  However, from the expression for $V^2$ given in the last paragraph we see that only two copies of 1110 can arise from $V^2$ and we conclude that $V\downarrow X$ is not MF.

Next we consider the case where $\l = a\l_9 + x\l_{10}$.  First assume $x > 0$.  Note that the case $a = x = 1$ is in the conclusion of the lemma, so assume this is not the case.  Then
$V^2(Q_Y) \supseteq (a+2)\l_9^0 \oplus (\l_8^0+a\l_9^0) \oplus (2\l_8^0+(a-2)\l_9^0) = a\l_9^0 \otimes 2\l_9^0$, where the third term occurs only if $a \ge 2$.  These summands are afforded by $\l-\b_{10}$, 
$\l-\b_9 -\b_{10}$, and $\l-2\b_9-\b_{10}$  (assuming $a \ge 2$) which have multiplicities $1,2$, and $3$,
respectively.
Restricting to $L_X'$ this becomes
$(V^1 \otimes 0001) + (V^1 \otimes 0200) $  and it will suffice to show
that the last expression is not MF.

If $a \ge 2$, then $V^1 \otimes 0200 \supseteq (0a00 + 0(a-2)01 + \cdots) \otimes 0200$ and this contains $(0a01)^2$ as required.

Suppose $a = 1$, and so $x \ge 2$.  Then $V^2 = 0100 \otimes S^2(0100) = 1000 + 0300 + (0101)^2 + 1110$.  Now $V^3(Q_Y) \supseteq 5\l_9^0 \oplus (\l_8^0+3\l_9^0)  \oplus (2\l_8^0+\l_9^0)$ where the summands are afforded by $\l - 2\b_{10},$  $\l - \b_9 - 2\b_{10},$ and  $\l - 2\b_9 - 2\b_{10},$ respectively.  This follows from the fact that these weights have multiplicities $1,2,$ and $3$ for $C_{10}.$  This is obvious
for the second weight and Magma gives the result for the other two weights using \cite[Prop. A]{cavallin}.
In the following we list the weights of $V^3(Q_Y)$ (using Magma for the expressions on the right hand sides):
\[
\begin{array}{rl}
5\l_9^0 = & S^5(\l_9^0) \\
\l_8^0+3\l_9^0 = & (3\l_9^0 \otimes \l_8^0) - (2\l_9^0 \otimes \l_7^0) + (\l_9^0 \otimes \l_6^0) - \l_5^0 \\
2\l_8^0+\l_9^0 = & (S^2(\l_8^0)\otimes \l_9^0) -   (\l_7^0\otimes \l_8^0).
\end{array}
\]
 Restricting to $L_X'$, Magma shows that $V^3 \supseteq (0301)^3$.  However $V^2$ can contribute only one such composition factor so we see that $V\downarrow X$ is not MF.
 
 We are left with the case where $\l = a\l_9$.  If $a = 1$, then $\l$ is as in the conclusion, so we assume $a \ge 2$.  Then
 $V^2(Q_Y) \supseteq (\l_8^0+a\l_9^0) \oplus (2\l_8^0+(a-2)\l_9^0)$ where the summands are afforded by
 $\l-\b_9 -\b_{10}$, and $\l-2\b_9-\b_{10}$.  This is  $(a\l_9^0 \otimes 2\l_9^0) - (a+2)\l_9^0$, and so
 $$V^2 \supseteq (S^a(0100) \otimes (0001 + 0200)) - S^{a+2}(0100)$$
and it suffices to show that $ (S^a(0100) \otimes 0200) - S^{a+2}(0100)$ is not MF.

If $a \ge 4$, then $S^a(0100) \otimes 0200 \supseteq (0a00+0(a-2)01+0(a-4)02) \otimes (0200)$ and
this contains $(0(a-2)02)^3$, whereas $S^{a+2}(0100)$ is MF.  And if $a = 3$  then $S^3(0100) \otimes (0200)$ contains $(1111)^2$ and $1111$ is not contained in $S^5(0100)$.

The final case is where $a = 2$.  Here $V^2(Q_Y) = (\l_8^0+2\l_9^0) \oplus 2\l_8^0\oplus (\l_7^0+\l_9^0)$ where the summands are afforded by $\l - \b_9 -\b_{10},$  $\l - 2\b_9 -\b_{10}$, and  $\l - \b_8 -2\b_9 -\b_{10}$.
Therefore $V^2(Q_Y) = (\l_8^0 \otimes 2\l_9^0) \oplus 2\l_8^0$ and we find that
$V^2 \supseteq (1011)^2$.
On the other hand, $V^1 = S^2(0100) = 0200 + 0001$ and this does not contribute a term $1011$ to
$V^2$.  Therefore $V\downarrow X$ is not MF in this case.  \hal

\begin{lem}\label{al10} If $\l = a\l_{10}\,(a>0)$, then $a=1$.
\end{lem}

\pf  Suppose $a \ge 2$. We first check that $V^2(Q_Y) = 
2\l_9^0$ which is afforded by $\l - \b_{10}$.  To see this consider weights of the form $\l - a_1\b_1-\cdots - a_9\b_9-  \b_{10}$ whose restriction to $L_Y'$ is dominant, and for which the weight space is not contained in $2\l_9^0$.

Next we claim that $V^3(Q_Y) = 4\l_9^0 \oplus 2\l_8^0$ where the summands are
afforded by weights $\l- 2\b_{10}$ and  $\l- 2\b_9 -2\b_{10}$. The corresponding weight spaces for $C_{10}$  have dimensions $1$ and $2$, respectively.  The only other possible
weights to consider that afford dominant weights for $L_Y'$ are $\l- \b_9 -2\b_{10}$, $\l- \b_8 -2\b_9 -2\b_{10}$ , and $\l- \b_7-2\b_8 - 3\b_9 -2\b_{10}$ where the dimensions of the
corresponding weight spaces are $1$, $2$ and $3$.  The first weight space is contained in $4\l_9^0$, the second has one occurrence there and another in $2\l_8^0,$ while
the third occurs once in $4\l_9^0$ and twice in $2\l_8^0$, establishing the claim. 
Therefore
\[
\begin{array}{ll}
V^3  &= S^4(0100) + S^2(\wedge^2(0100)) - \wedge^4(0100) \\
 &= 0400+(0201)^2+(0002)^2+2020+1100.
\end{array}
\]
Assume $a \ge 3$.  We claim $V^4(Q_Y) \supseteq 6\l_9^0 \oplus (2\l_8^0+2\l_9^0) \oplus 2\l_7^0$ where the summands are afforded by  weights $\l- 3\b_{10}$,   $\l- 2\b_9- 3\b_{10}$, and  $\l-2\b_8- 4\b_9- 3\b_{10}$ respectively. Using \cite[Prop. A]{cavallin} and a Magma check, we see that the corresponding weight spaces have dimensions $1$,  $2$, and $5$, respectively. 
Call these weights $(1)$,  $(2)$, and $(3)$. Clearly $(1)$ occurs as the highest weight of a composition factor. To make sure that $(2)$ occurs we first observe that this weight appears with multiplicity $1$ in $(1)$.  Next we consider the weight $\l- \b_9- 3\b_{10}$.  This weight occurs
with multiplicity $1$ in $V$ and it already appears in $(1)$.  Therefore $(2)$ does indeed appear
as a composition factor.  Weight (3) appears with multiplicity $5$ and it has multiplicities $1$ and $3$ in
$(1)$ and $(2)$.  Additional checks with Magma show that there are no additional composition factors for which $(3)$ appears as a weight.  This gives the claim.

We will study the multiplicity of $(0202)$ in $V^4$.  First note that
$$6\l_9^0 \oplus (2\l_8^0+2\l_9^0)   = (2\l_9^0 \otimes 4\l_9^0)  - (\l_8^0+4\l_9^0).$$ 
Magma shows that $S^2(0100)\otimes S^4(0100)\supseteq (0202)^4$.    Moreover
$(\l_8^0+4\l_9^0) \downarrow L_X'  \subseteq  \wedge^2(0100) \otimes S^4(0100)$ and $0202$ has multiplicity 1 in the tensor product.  Next consider $2\l_7^0 = S^2(\l_7^0) - (\l_5^0+\l_9^0)$. We check that $S^2(\wedge^3(0100)) \supseteq (0202)^2$ whereas the restriction of
$(\l_5^0+\l_9^0)$ to $L_X'$ is contained in $(0100) \otimes \wedge^5(0100)$ and here
the multiplicity of $(0202)$ is 1.  It follows that $V^4 \supseteq (0202)^4$.  But at most
two such terms can arise from $V^3$. Therefore $V\downarrow X$ is not MF. 

Finally assume $a =2$.  Here there are some changes in $V^4(Q_Y)$.  Indeed $6\l_9^0$ does not
occur as $\l-3\b_{10}$ is not a weight.  However $\l- 2\b_9- 3\b_{10} = 2\l_8+2\l_9 -2\l_{10}$ is conjugate in $C_{10}$ to $2\l_9$ which has multiplicity $1$ in $\l$, so $V^4(Q_Y) \supseteq 2\l_8^0+2\l_9^0 = \nu$.  In addition the summand $\mu = 2\l_7^0$, which is afforded by $\l - 2\b_8-4\b_9-3\b_{10}$, occurs.  Indeed this weight has multiplicity $4$ in $\l$ and $\mu$ occurs with multiplicity 3 in $\nu$, and occurs in no other summand of $V^4(Q_Y)$.  Therefore
$V^4(Q_Y) \supseteq (2\l_8^0+2\l_9^0) \oplus 2\l_7^0$.  A careful check  shows that this is actually an equality.  Now 
\[
\begin{array}{ll}
2\l_7^0\downarrow L_X'  &=  S^2(\wedge^3(0100)) - (\wedge^5(0100) \otimes (0100)) + \wedge^6(0100) \\
 & = 4002+0040+1020+2021+2000+0202+2110+1101+0003.
\end{array}
\]
From the $a\ge 3$ analysis  above we have $(2\l_8^0+2\l_9^0)   = (2\l_9^0 \otimes 4\l_9^0)  - 6\l_9^0 -(\l_8^0+4\l_9^0)$, and using Magma we check that 
\[
(\l_8^0+4\l_9^0) = (\l_8^0\otimes 4\l_9^0) - (\l_7^0\otimes 3\l_9^0) + (\l_6^0\otimes 2\l_9^0) - (\l_5^0\otimes \l_9^0) + \l_4^0.
\]
Hence, restricting to $L_X'$ we find that 
we find that
\[
\begin{array}{ll}
(2\l_8^0+2\l_9^0) \downarrow L_X' = & 1020+0210+2110+2021+ 
1101^2+2220+0202^2 \\
& +1211+0011+0401 +0003+1300+2000.
\end{array}
\]
 We next claim that 
\begin{equation}\label{v5eq}
V^5(Q_Y) \supseteq  4\l_8^0 \oplus (2\l_7^0+2\l_9^0)\oplus 2\l_6^0.
\end{equation}
  These weights arise from $\l -4\b_9 -4\b_{10}= 4\l_8 -2\l_{10}$,   $\l -2\b_8-4\b_9 -4\b_{10} = 2\l_7 +2\l_9-2\l_{10}$, and $\l -2\b_7-4\b_8-6\b_9 -4\b_{10} = 2\l_6$, respectively.  They are conjugate in $C_{10}$ to   $2\l_{10}$, $2\l_8$, and $2\l_6$ which have multiplicities $1$,  $2$, and $10$ for $C_{10}$, respectively.  Obviously $2\l_7^0+2\l_9^0$ appears with multiplicity $1$ in $4\l_8^0$ and Magma shows that  $2\l_6^0$ occurs with multiplicity $3$ and $6$ in
$4\l_8^0$ and $2\l_7^0+2\l_9^0$, respectively, so this gives the claim once we check that these weights are not subdominant to any other weights in $V^5(Q_Y)$.  This is a straightforward check.

We next use Magma to check that the right hand side of (\ref{v5eq}) is equal to $S^2(2\l_8^0) -(\l_6^0 + 2\l_8^0)$, and hence
\[
\begin{array}{ll}
V^5(Q_Y)  \supseteq & S^2(2\l_8^0) -(\l_6^0 + 2\l_8^0)  \\
 & = S^2(2\l_8^0) - (\l_6^0 \otimes 2\l_8^0) + (\l_5^0+\l_8^0+\l_9^0)+(\l_4^0+\l_8^0).
\end{array}
\]
Also
\[
\begin{array}{l}
(\l_5^0+\l_8^0+\l_9^0)\oplus (\l_4^0+\l_8^0) = ((\l_5^0+\l_8^0) \otimes \l_9^0) - (\l_5^0+\l_7^0), \hbox{ and} \\
2\l_8^0 = S^2(\l_8^0) - \l_6^0.
\end{array}
\]
At this point we can restrict terms  to $L_X'$ and use Magma to compute that $V^5 \supseteq (2022)^4 $. However $V^4$
can only account for two copies of $2022$.  This completes the proof.
 \hal

\begin{lem}\label{al1l2} Assume $\mu^0  = a\l_1^0 +\l_2^0\,(a\ge 2)$ or its dual.  Then $\l = 2\l_1 +\l_2$. 
\end{lem} 

\pf  First assume $\l =  a\l_1 +\l_2$ with $a\ge 3$.  Lemma \ref{tensorl1} shows that 
$$a\l_1 +\l_2 = (\l_1 \otimes (a+1)\l_1) - (a+2)\l_1 - a\l_1$$
so that 
$$V \downarrow X \supseteq (00100 \otimes S^{a+1}(00100)) - S^{a+2}(00100) - S^a(00100).$$
For $a = 3$ and $a = 4$ a Magma check shows that  $V \downarrow X$ contains $(01110)^2$ and $(10201)^3$, respectively, so these cases are not MF.

Now assume $a \ge 5.$ We claim that 
$$S^{a+1}(00100) \supseteq 10(a-1)01+ 00(a-1)00 + 20(a-3)02 + 01(a-3)10 + 10(a-3)01 .$$
To see this we first note that these weights have the form $00(a+1)00 - \mu$ where $\mu$
is $\a_2 + 2\a_3+\a_4$, $\a_1+2\a_2 + 3\a_3+2\a_4+\a_5$, $2\a_2 + 4\a_3+2\a_4$, $\a_1+2\a_2 + 4\a_3+2\a_4+\a_5$, and $\a_1+3\a_2 + 5\a_3+3\a_4+\a_5$, respectively.  It follows from
the form of $\mu$ and \cite[Prop. A]{cavallin} that the multiplicities of the above irreducibles in $S^{a+1}(00100)$ are the same as they are when $a = 5$ and one checks using Magma that 
they each occur in $S^6(00100)$.  This establishes the claim.

Next we use the  Littlewood-Richardson rule to see that tensoring each of the five modules displayed above 
with $00100$ yields a summand $10(a-2)01$.  Therefore 
$$00100 \otimes S^{a+1}(00100) \supseteq (10(a-2)01)^5.$$
Next we observe that $10(a-2)01 = (00(a+2)00) - (\a_1+3\a_2 + 5\a_3+3\a_4+\a_5)$ and
$10(a-2)01 = (00a00) - (\a_2 + 2\a_3+\a_4)$, so \cite[Prop. A]{cavallin} implies that the multiplicity of $10(a-2)01$ in $S^{a+2}(00100)$ and in $S^a(00100)$ is the same as it is when $a= 5$. In each case this multiplicity is $1$ and hence we conclude that $V\downarrow X \supseteq (10(a-2)01)^3$ and so this is not MF.

Now consider the case where $\l =  a\l_1 +\l_2 + x\l_{10}$ with $x > 0$, $a\ge 2$.   Here we aim to show
that $V\downarrow X$ is not MF.  We find $4$ composition factors in $V^2(Q_Y)$ afforded
by $\l-\b_{10}$, $\l - \b_2 - \cdots - \b_{10}$, $\l -\b_1 - \b_2 - \cdots - \b_{10}$, and $\l -\b_1 -2\b_2 - \cdots -2\b_9 - \b_{10}$.  To see that the last one occurs we use Magma to see that it has multiplicity $90$ in $\l$ and multiplicities $72, 9,$ and $8$  in the other composition factors. Therefore 
\[
\begin{array}{ll}
V^2(Q_Y) \supseteq & (a\l_1^0 +\l_2^0 +2\l_9^0) \oplus ((a+1)\l_1^0 +\l_9^0) \oplus ((a-1)\l_1^0 +\l_2^0 +\l_9^0) \oplus a\l_1^0 \\
 &  = (a\l_1^0 +\l_2^0) \otimes 2\l_9^0.
\end{array}
\]
Restricting to $L_X'$ this is $(V^1 \otimes 0001) + (V^1 \otimes 0200)$, 
so it suffices to show that $V^1 \otimes 0200$ is not MF.  Now $V^1(Q_Y) = ((a+1)\l_1^0 \otimes \l_1^0)- (a+2)\l_1^0$ so that 
$$V^1 = (S^{a+1}(0010) \otimes 0010) - S^{a+2}(0010).$$
 Now $S^{a+1}(0010) \supseteq 00(a+1)0 + 10(a-1)0$
and hence 
$$S^{a+1}(0010)  \otimes 0010 \supseteq (00(a+1)0 + 10(a-1)0) \otimes 0010 \supseteq (10a0)^2 +01a1.$$
Therefore $V^1 \supseteq 10a0 + 01a1$ (note that $01a1 \not \subseteq S^{a+2}(0010)$ by \cite[3.8.1]{Howe}), and so  $V^1\otimes 0200 \supseteq (12a0)^2,$ as required.  

We must still consider duals, so we now assume $\l = \l_8 +a\l_9 + x\l_{10}$ with $a\ge 2$. 
We claim that this fails to be MF for all values of $x$. Taking duals of the above we have
$$V^1 = S^{a+1}(0100) \otimes 0100 - S^{a+2}(0100).$$
Magma checks show that we get composition factors in $V^2(Q_Y)$ of highest weights $(2\l_8^0 + a\l_9^0)$, $(\l_7^0 + (a+1)\l_9^0) $, 
$ (3\l_8^0 + (a-2)\l_9^0)$, and $ (\l_7^0 + \l_8^0 + (a-1)\l_9^0)$ afforded by $\l-\b_9-\b_{10}$, $\l-\b_8 -\b_9-\b_{10}$, $\l -2\b_9-\b_{10}$ and $\l-\b_8 -2\b_9-\b_{10}$, respectively. (The counts are slightly different depending on whether $x = 0$ or $x > 0$ since in the latter case there is also a composition factor afforded by $\l - \b_{10}$ which we do not use in the argument to follow.)  Therefore in any case we have 
$$V^2(Q_Y) \supseteq  ((\l_8^0 + a\l_9^0)  \otimes 2\l_9^0) - (\l_8^0 + (a+2)\l_9^0)$$
and restricting to $L_X'$ this becomes
$$(V^1 \otimes 0001) + (V^1 \otimes 0200) - (\l_8^0 + (a+2)\l_9^0) \downarrow L_X'.$$
By Corollary 5.1.5 of \cite{MF} the first summand covers those composition factors arising from $V^1$ and the subtracted
term is MF (by \cite{MF}).  Therefore it will suffice to show that  $V^1 \otimes 0200$ contains a composition factor of multiplicity at least $3$.  Now $S^{a+1}(0100) \supseteq 0(a+1)00 + 0(a-1)01$ and tensoring
these summands with $0100$ we get terms $1a10$ and $1(a-2)11$, respectively.  Neither of these terms are summands of $S^{a+2}(0100)$, so from the above  expression for $V^1$ (and using \cite{cavallin} as usual), we see
that 
$$V^1 \otimes 0200 \supseteq (1a10 + 1(a-2)11) \otimes 0200 \supseteq (1a11)^3$$
which proves the claim and completes the proof of the lemma.  \hal

\begin{lem}\label{l1li} Assume $\mu^0  = \l_1^0 +\l_i^0$ with $2 \le i \le 9$, or its dual.  Then $\l = \l_1 +\l_2$, $\l_1 +\l_3,$ or $ \l_1 +\l_4$.  
\end{lem} 

\pf First assume that $\l = \l_1 + \l_i$ with $i>4$. Then Lemma \ref{tensorl1} shows that $\l_1 + \l_i = (\l_1 \otimes \l_i) - \l_{i+1} - \l_{i-1}$, so that
$$V\downarrow X = (\d \otimes (\wedge^i(\d) - \wedge^{i-2}(\d))) -\wedge^{i+1}(\d) + \wedge^{i-3}(\d).$$
 For $i = 5,6,7,8$ and $9$ we use Magma to see that $V \downarrow X$ contains $(11011)^2$, $(11111)^2$, $(21012)^2$, $(11111)^2$,
and $(20202)^2$, respectively.  Therefore none of these is MF.

Next assume $\l = \l_1 + \l_i+x\l_{10}$ with $x > 0$. We claim that none of these
cases is MF.  We have 
\[
\begin{array}{ll}
V^2(Q_Y) \supseteq & (\l_1^0 + \l_{i-1}^0 + \l_9^0) \oplus ( \l_i^0 + \l_9^0) \oplus (\l_1^0 + \l_i^0 + 2\l_9^0) \\
 & = (\l_1^0 + \l_i^0) \otimes (2\l_9^0) - \l_{i-1}^0.
\end{array}
\]
Restricting to $L_X'$ this becomes
$$(V^1 \otimes 0001) + (V^1 \otimes 0200) - \wedge^{i-1}(0010).$$
 So since $\wedge^{i-1}0010$ is MF, as usual it will suffice to show that $V^1 \otimes 0200$ contains a composition factor of multiplicity $3$.  Note that $V^1 = (0010 \otimes \wedge^i(0010)) - \wedge^{i+1}(0010)$ so a Magma check shows that $V^1 \otimes 0200$ contains $(0201)^3$, $(1110)^6$, $(1201)^7$, $(1101)^9$, $(1111)^{10}$, $(1011)^7$, $(1110)^5$, when  $i = 2,3,4,5,6,7,8$, respectively.  This leaves the case $i = 9$, where $V^1 \otimes 0200 \supseteq (1201)^2$ and this does not appear in $\wedge^8(0010)$; hence here as well $V^1\otimes 0200-\wedge^{i-1}(0010)$ contains repeated summands, giving  the usual contradiction.
 
 Now assume that $\l = \l_i + \l_9+x\l_{10}$.  The $i = 1$ case was settled above so we can assume that $2\le i \le 8$.  We have
\[
\begin{array}{ll}
V^2(Q_Y) \supseteq & ( \l_{i-1}^0 + 2\l_9^0) \oplus (\l_i^0 + \l_8^0 + \l_9^0) \\
 & =  (\l_i^0 + \l_9^0) \otimes (2\l_9^0) - (\l_i^0 + 3\l_9^0) - (\l_{i-1}^0+\l_8^0).
\end{array}
\]
Restricting to $L_X'$ this becomes
\begin{equation}\label{eq3}
(V^1 \otimes 0001) + (V^1 \otimes 0200)  - (\l_i^0 + 3\l_9^0)\downarrow L_X' - (\l_{i-1}^0+\l_8^0)\downarrow L_X'.
\end{equation}
By Corollary 5.1.5 of \cite{MF} the first summand covers all composition factors in $V^2$ that arise from $V^1$ so it remains
to show that the remaining terms fail to be MF.  One checks that
$$V^1 \otimes 0200 =  (\wedge^i(0010) \otimes (0300+1110+0101)) - (\wedge^{i-1}(0010) \otimes 0200).$$
 For the subtracted terms in (\ref{eq3}) we have
\[
\begin{array}{ll}
(\l_i^0 + 3\l_9^0)\downarrow L_X' = & (\wedge^i(0010) \otimes S^3(0100)) - (\wedge^{i-1}(0010) \otimes S^2(0100)) \\
& +  (\wedge^{i-2}(0010) \otimes 0100) - \wedge^{i-3}(0010),
\end{array}
\]
where the last term is omitted  if $i = 2.$  Similarly
$$(\l_{i-1}^0+\l_8^0)\downarrow L_X' = (\wedge^{i-1}(0010) \otimes \wedge^2(0100)) - (\wedge^{i-2}(0010) \otimes 0100) .$$
So the subtracted terms in (\ref{eq3}) sum to 
$$(\wedge^i(0010) \otimes S^3(0100)) - (\wedge^{i-1}(0010) \otimes S^2(0100))+(\wedge^{i-1}(0010) \otimes \wedge^2(0100)) -\wedge^{i-3}(0010).$$
Now $S^3(0100) = 0300 + 0101$ and $S^2(0100) = 0200 + 0001$ so combining the expression for $V^1 \otimes 0200$ with the subtracted terms we obtain
$$  (\wedge^i(0010) \otimes 1110) + (\wedge^{i-1}(0010) \otimes 0001)
 - (\wedge^{i-1}(0010) \otimes 1010) +\wedge^{i-3}(0010)$$
  and we must show that this is not MF.  A Magma check shows that for $i =2,3,4,5,6$ and $7$ there are repeated  summands of highest weights $2110,1111,1121,2111,1120,$ and $1211$, respectively.
So this leaves the case $i =8.$ 

For this case we have used the composition factors in $V^2(Q_Y)$ afforded by $\l -\b_9 -\b_{10}$ and $\l -\b_8-\b_9 -\b_{10}$.   We claim that $\l -\b_8-2\b_9 -\b_{10}$ also appears.  Indeed this has multiplicity $4$ or $6$ for $C_{10}$ according to whether or not $x = 0$, and multiplicity $2$ and $1$, respectively, in the other two factors.  And if $x > 0$ it has multiplicity $2$ in $\l_8^0 + 3\l_9^0$ which is afforded
by $\l - \b_{10}$.  This establishes the claim.
Using this we see that
$$V^2(Q_Y) \supseteq   (\l_8^0 + \l_9^0) \otimes (2\l_9^0) - (\l_8^0 + 3\l_9^0).$$
Arguing as above we need to show that
$$(\wedge^8(0010) \otimes 1110) + (\wedge^7(0010) \otimes 0001) - (\wedge^6(0010) \otimes 0100)
 +\wedge^5(0010)$$
is not MF.  But a Magma check shows that this contains repeated copies of $1111$, completing
the proof.   \hal

\begin{lem}\label{li813} Assume $\mu^0  =\l_i^0$.  Then $\l =\l_i$ or $\l_1+\l_{10}$.  
\end{lem} 

\pf Assume $\l = \l_i + x\l_{10}$ with $x > 0$.  Then
$V^2(Q_Y) \supseteq (\l_{i-1} +\l_9) \oplus (\l_i + 2\l_9) = \l_i \otimes 2\l_9$.  Restricting to $L_X'$ this becomes $(V^1 \otimes 0001) + (V^1 \otimes 0200)$, so $V\downarrow X$ is not MF provided
we can show that $V^1 \otimes 0200$ is not MF. For $i = 3,4,5,6$ and $7$ a Magma check shows that
$V^1 \otimes 0200$  contains $(0201)^2, (2001)^2, (2011)^2, (1020)^2$, and $(1111)^2$, respectively.
So none of these cases are MF.

So we must now consider the cases where $i = 1,2,8,$ and $9$.  The case $i = 1$ was settled in 
Lemma \ref{al1}.  We begin with the case $i = 2.$

If $i = 2$ and $x = 1$ then a Magma computation gives 
$(\l_2 +\l_{10}) = (\l_2 \otimes \l_{10}) - (\l_1 \otimes \l_9) +\l_{10}.$
Restricting to $L_X'$ we see that
$$V \downarrow X = ((\wedge^2(\d)-0) \otimes (\wedge^{10}(\d)-\wedge^8(\d)) - (\d \otimes (\wedge^9(\d) - \wedge^7(\d)) + (\wedge^{10}(\d) - \wedge^8(\d))$$
and this contains $(21012)^3$ and hence $V\downarrow X$ is not MF.

Next suppose $i = 2$ and $x \ge 2$.  A check shows that $V^2(Q_Y)  =  (\l_2^0 + 2\l_9^0)  \oplus (\l_1^0 + \l_9^0) = \l_2^0 \otimes 2\l_9^0,$ where the summands are afforded by $\l - \b_{10}$ and $ \l - \b_2 - \cdots - \b_{10}$, respectively.  Restricting to $L_X'$ we see that
$$V^2 = (1001)^2 +2010 + 1200 +(0110)^2 + 1111 + 0301 +(0102)^2.$$
 We have $V^3(Q_Y) \supseteq (\l_2^0 + 4\l_9^0) \oplus (\l_2^0 + 2\l_8^0)$ where the summands are afforded by 
$\l - 2\b_{10}$ and $ \l -  2\b_9- 2\b_{10}$.  To see that the second summand occurs we note that $ \l -2\b_9- 2\b_{10}$ has multiplicity $2$ for $C_{10}$ but it only has multiplicity $1$ in $\l_2^0 + 4\l_9^0$, and also that 
$\lambda-\beta_9-2\beta_{10}$ has multiplicity 1 for $C_{10}$.  Now
\begin{equation}\label{eq5}
\begin{array}{ll}
 (\l_2^0 + 4\l_9^0) \oplus (\l_2^0 + 2\l_8^0) = & ((\l_2^0 + 2\l_9^0) \otimes 2\l_9^0) - (\l_1^0 + 3\l_9^0) - \\
 & (\l_1^0 + \l_8^0 +\l_9^0) - (\l_2^0 +\l_8^0 + 2\l_9^0).
\end{array}
\end{equation}
Using the fact that $(\l_2^0 + 2\l_9^0) = (\l_2^0 \otimes 2\l_9^0) - (\l_1^0 + \l_9^0)$, we find that restriction to $L_X'$ of the tensor product in (\ref{eq5}) contains $(0302)^5$.  Of the subtracted terms only the last one can have a summand  of $S$-value
at least $5$. Moreover $(\l_2^0 +\l_8^0 + 2\l_9^0) \downarrow L_X' \subseteq  \wedge^2(0010) \otimes \wedge^2(0100) \otimes S^2(0100)$ and this only contains $(0302)^2$. Therefore $V^3 \supseteq (0302)^3$, whereas only one such summand can arise from $V^2$.

Consider the case  $i = 9$ and $x\ge 2$.    We argue as above
that $V^2(Q_Y) = 3\l_9^0 \oplus (\l_8^0 +\l_9^0) = \l_9^0 \otimes 2\l_9^0$ and restricting to $L_X'$ we have 
$$V^2 = 1000 + 0300 + (0101)^2 + 1110.$$
We next claim that 
$$V^3(Q_Y) \supseteq 5\l_9^0 \oplus (\l_8^0 +3\l_9^0) \oplus (2\l_8^0 +\l_9^0)$$
which can be expressed $3\l_9^0 \otimes 2\l_9^0$.
To see that the three summands occur we note that these weights are afforded by $\l - 2\b_{10}$, $\l - \b_9 -2\b_{10}$, and $\l - 2 \b_9 -2\b_{10}$, and the dimensions of the corresponding weight spaces of $V$ are $1,2,$ and $3$, respectively.  From here it is easy to see that the claim holds.
It follows that
$$V^3 \supseteq S^3(0100) \otimes S^2(0100) \supseteq (0301)^3.$$
Only one such composition factor can arise from $V^2$ and we conclude that $V\downarrow X$ is not MF in this case either.

Now consider the case $i = 9$ and $x=1$, that is $\l = \l_9 +\l_{10}$.  We first claim that $V^2(Q_Y) = 3\l_9^0 \oplus (\l_8^0 +\l_9^0)$, where the summands are afforded by $\l - \b_{10}$ and $\l - \b_9 - \b_{10}$, respectively.  Call these weights and their corresponding restrictions $(1)$ and $(2)$, respectively.  The summand $(1)$ clearly occurs.  So does $(2)$ because  $\l - \b_9 - \b_{10}$ has multiplicity $2$ in $\l$,  and $(2)$ only has multiplicity $1$ in $(1)$.  The only other possible dominant weight for $V^2(Q_Y)$  is afforded by $\l - \b_8 -2\b_9 - \b_{10}$, which has multiplicity $3$. But this weight has   multiplicity $1$ and $2$ in $(1)$ and $(2)$,  respectively, and  so  does not occur  as the highest weight of a composition factor. This  establishes the claim.  As above we have  $V^2(Q_Y) = 2\l_9^0 \otimes \l_9^0$ and $V^2 = 1110 + 0300 + (0101)^2+1000.$
We next claim that  $V^3(Q_Y) \supseteq (\l_8^0 +3\l_9^0) \oplus (2\l_8^0 +\l_9^0) \oplus (\l_7^0 +\l_8^0)$, where the summands are afforded
by $\l -\b_9 - 2\b_{10}$, $\l -2\b_9 - 2\b_{10}$, and $\l -\b_8- 3\b_9 - 2\b_{10}$, respectively.  We call these weights and their
corresponding restrictions $(3)$, $(4)$, and $(5)$.  The first two weights are conjugate in $C_{10}$ to $\l_9 + \l_{10}$ and $\l_8+\l_9$, respectively and these weights have respective multiplicities $1$ and $2$.  They both occur since
the second only appears once in $(3)$.  The only other dominant weight above $(5)$ is afforded by $\l -\b_8- 2\b_9 - 2\b_{10}$.
This weight is conjugate to $\l_7+\l_{10}$ which has multiplicity $3$.  However, the restriction has respective multiplicities $2$ and $1$ in $(3)$ and $(4)$, so we do not get a summand from this weight.  On the other hand $(5)$ has 
multiplicity $2$ in each of $(3)$ and $(4)$ and it has multiplicity $5$ in $V$.  Hence $(5)$ does occur as a summand.  It follows that
$$V^3(Q_Y) \supseteq ((\l_8^0 +\l_9^0) \otimes 2\l_9^0) - (\l_7^0 +2\l_9^0)$$
and restricting to $L_X'$ and writing $\d'' = 0100$, we have 
$$V^3 \supseteq (((\d'' \otimes \wedge^2(\d'')) - \wedge^3(\d'')) \otimes S^2(\d'')) - (\wedge^3(\d'')\otimes S^2(\d'') )+ (\wedge^4(\d'')\otimes \d'') - \wedge^5(\d'').$$
A Magma computation shows that $V^3 \supseteq (1111)^3$ but $V^2$ can contribute only one such composition factor.
Hence $V\downarrow X$ is not MF.

The final case is where $i = 8$. First assume $x \ge 2$. Arguing as in previous cases, we see that $V^2(Q_Y) = (\l_8^0 +2\l_9^0) \oplus (\l_7^0 +\l_9^0) = \l_8^0 \otimes 2\l_9^0$. Restricting to $L_X'$ we find that
\begin{equation}\label{randomeq}
V^2 = (0010)^2+0201+0120+(1100)^2+(1011)^2+2101+1210.
\end{equation}
Again we work with $V^3(Q_Y)$.  We claim that
$$V^3(Q_Y) \supseteq (\l_8^0 +4\l_9^0) \oplus 3\l_8^0\oplus (\l_7^0 +3\l_9^0)\oplus (\l_7^0 +\l_8^0 +\l_9^0)$$
where the highest weights are afforded by $\l - 2\b_{10},$ $\l - 2 \b_9 -2\b_{10}$,
 $\l - \b_8 - \b_9 -2\b_{10}$, and $\l - \b_8 -2 \b_9 -2\b_{10},$ respectively.  Call these
 weights and corresponding summands A,B,C, and D, respectively.  A Magma check and an application of \cite[Prop. A]{cavallin} shows that the dimensions of the corresponding weight spaces are $1,2,3$, and $5$, respectively.  Clearly A occurs.
 The second weight occurs with multiplicity 1 in A, so it will occur as a summand provided there is no
 summand corresponding to the weight $\l -  \b_9 -2\b_{10}$.  But this weight space has multiplicity $1$ and it appears in A.  Therefore B appears in $V^3(Q_Y)$.  The weight C appears twice in  A but not in B, so it also appears.  Finally, D appears twice in A, once in B, and once in C, so it also appears.
 
We note that  $(\l_8^0 +4\l_9^0) \oplus (\l_7^0 +3\l_9^0) = \l_8^0 \otimes 4\l_9^0$ and that
\begin{equation}\label{3l8}
3\l_8^0 = S^3(\l_8^0) - \l_4^0 - (\l_6^0 \otimes  \l_8^0) + (\l_5^0 \otimes \l_9^0).
\end{equation}
Next we check that
\begin{equation}\label{l7l8}
(\l_7^0 +\l_8^0 +\l_9^0) = ((\l_7^0\otimes \l_8^0) - (\l_6^0\otimes \l_9^0)) \otimes \l_9^0 ) - (\l_6^0 \otimes \l_8^0)+ (\l_5^0 \otimes \l_9^0) - 2\l_7^0
\end{equation}
and the last term is contained in $S^2(\l_7^0).$
Restricting all the terms to $L_X'$  a Magma check shows that $V^3$ contains $(0210)^6$.  On the other hand  only $3$ copies of 0210 can arise from $V^2$.  So $V\downarrow X$ is not MF, completing the proof.

Finally assume $x=1$ so that $\l = \l_8 + \l_{10}.$  We first claim that $V^2(Q_Y) = (\l_8^0 +2\l_9^0) + (\l_7^0 +\l_9^0)$.  To see this first note that the summands are afforded by $ \l -\b_{10}$ and $\l -\b_8-\b_9-\b_{10}$, respectively.  Call these summands $(1)$ and $(2)$, respectively.  The corresponding weights appear with multiplicities $1$  and $3$, respectively. Clearly $(1)$ occurs and the second weight occurs
with multiplicity $2$ in $(1)$.  Furthermore $\l-\b_9-\b_{10}$ is conjugate to $2\l_8$ which occurs with multiplicity $1$ and already appears in $(1)$.  So this shows that $V^2(Q_Y)$ contains $(1)$ and $(2)$.  The only other possible weight is $\l - \b_7-2\b_8-2\b_9-\b_{10}$.  But this appears with multiplicity $3$ in both $(1)$ and $(2)$ and it appears with multiplicity $6$ in $\l$.  Hence it does not occur as a composition factor and this establishes the claim.  Therefore  $V^2(Q_Y) = \l_8^0 \otimes 2\l_9^0$ and $V^2$ is as in (\ref{randomeq}). 
We next consider $V^3(Q_Y)$.  Neither  $\l-2\b_{10}$ nor $\l-\b_9-2\b_{10}$ is a weight of $V$.  This is
clear for the first weight and the second is conjugate to $2\l_9$ which is  $\l +\b_9$ and hence is not subdominant to $\l$.  On the other hand $\l-2\b_9-2\b_{10}$ and $\l-\b_8-\b_9-2\b_{10}$ are each
conjugate to $\l$ and these weights afford composition factors $3\l_8^0$ and $\l_7^0 + 3\l_9^0$, respectively.  The last weight to consider is $\l - \b_8-2\b_9-2\b_{10}$ which has multiplicity $3$ and which restricts to $\l_7^0 + \l_8^0 +\l_9^0$.  The restriction appears with multiplicity $1$ in each of the weights produced so far and there are no other dominant weights above it.  Therefore $V^3(Q_Y) \supseteq 
3\l_8^0 \oplus (\l_7^0 + 3\l_9^0) \oplus (\l_7^0 + \l_8^0 +\l_9^0)$.   We next work out what occurs in $V^3$.
We have 
\[
(\l_7^0 + 3\l_9^0) =  (\l_7^0 \otimes 3\l_9^0) - (\l_6^0 + 2\l_9^0), 
\]
and together with (\ref{3l8}) and (\ref{l7l8}), this enables us to use Magma to restrict to $L_X'$ and find that  $V^3 \supseteq (2110)^4$.  Only $(2110)^2$ can arise
from $V^2$ and hence $V\downarrow X$ is not MF.
 \hal

\begin{lem}\label{al2} Assume $\mu^0  =a\l_2^0$ with $2 \le a \le 5$, or its dual.  Then $\l =2\l_2$.
\end{lem}

\pf   First assume $\l = a\l_2 + x\l_{10}$ with $x > 0$.   We show that $V\downarrow X$ is not MF. 

We first claim that  $V^2(Q_Y) \supseteq (\l_1^0 + (a-1)\l_2^0 +\l_9^0) \oplus (a\l_2^0 +2\l_9^0) \oplus (2\l_1^0 + (a-2)\l_2^0) = a\l_2^0 \otimes 2\l_9^0 $.  
To see this, note that the first two composition factors clearly occur.  The third weight is afforded by $\l - 2\b_2 - \cdots -2\b_9 - \b_{10}$. To check that this occurs as a composition factor it follows from  \cite[Prop. A]{cavallin} that we can use use Magma with $a = 2$ and $x = 1$. We find that the third weight appears with multiplicity $45$ in $V$  but it has multiplicities $8$ and $36$ in the first two summands, respectively.  This implies that the third weight does appear as a composition factor of $V^2(Q_Y)$, as claimed. Restricting to $L_X'$ we therefore have
$$V^2 \supseteq (V^1 \otimes 0001) + (V^1 \otimes 0200).$$
As usual it will suffice to show that the second summand is not MF.   For each possible value $2\le a \le 5$,  Table 6.2 of \cite{MF} together with a Magma computation shows
 that the second summand is indeed not MF.

Next assume that  $\l = a\l_2$ with $a\ge 3$. If $\l = 3\l_2$ then a Magma computation gives 
\[
V_Y(3\l_2) = S^3(\l_2) - S^2(\l_2) -\l_6 - (\l_2\otimes \l_4) + (\l_1\otimes \l_5),
\]
and restricting to $X$, there are multiple copies of $01010$. Hence
this case is not MF and we can assume that $a=4$ or 5.

Set $b = a-1$ so that $b = 3$ or $4$. Then $V^2(Q_Y) \supseteq (\l_1^0 + b\l_2^0 +\l_9^0) = \xi$. There is another composition factor as well, but we will not need that summand.

We begin by considering the embedding of $L_X'$ in $A_9$ via the weight $\d' =0010. $ We can
take as basis vectors of weights $\d'$, $\d' - \a_3$, $\d' - \a_3-\a_4$, $\d' -\a_2- \a_3$, $\d' -\a_2- \a_3-\a_4$, $\d' -\a_1- \a_2-\a_3$, $\d' -\a_1-\a_2- \a_3-\a_4$, $\d' -\a_2- 2\a_3-\a_4$, $\d' -\a_1-\a_2- 2\a_3-\a_4$, $\d' -\a_1-2\a_2- 2\a_3-\a_4$.  The basis for the natural module of $A_9$ consists of weight vectors of weights $\l_1^0, \l_1^0-\b_1^0, \l_1^0 - \b_1^0-\b_2^0, \cdots$ and we find that  $\b_1^0, \cdots, \b_9^0$ restrict to $S_X$ respectively as:
\[
 \a_3,\, \a_4,\, \a_2-\a_4,\,\a_4,\,\a_1- \a_4,\,\a_4,\ \a_3-\a_1,\,\a_1,\, \a_2.
\]
 We know that $\l_1^0, \l_2^0,$ and $\l_9^0$ restrict to $L_X'$ as $0010, 0101,$ and $0100$, respectively.
 Therefore the maximal vector of the representation $(\l_1^0 + b\l_2^0 +\l_9^0)$ is a weight vector for $L_X'$ of  weight
 $\gamma = 0(b+1)1b$ and $S$-value $2b+2$. We first argue that this vector is also a maximal vector for $L_X'$.   Indeed $\l_1^0 + b\l_2^0 +\l_9^0 \subseteq  S^b(\l_2^0) \otimes (\l_1^0 +\l_9^0)= S^b(\l_2^0) \otimes (\l_1^0 \otimes \l_9^0 -0)$,  and restricting
 to $L_X'$ this is $S^b(0101) \otimes (0110 + 1001)$.  Now use Magma to see
 that any composition factor of $S $-value at least $2b+2$ must be contained in $(0b0b) \otimes (0110 + 1001)$ and the only possible highest weights are $\gamma$, $\gamma-\a_2$,  $\gamma-\a_2-\a_3$, and $\gamma-2\a_2-\a_3$.  Therefore a weight vector of weight $\gamma$ is necessarily a maximal vector.  
 
Consider the weight $\gamma-\a_2$, which is the restriction of $\xi-(\b_2^0+\b_3^0)$ and of $\xi-\b_9^0$, but only has  multiplicity $1$ in the summand afforded by $\gamma$.  It has $S$-value $2b+2$ so we see from the above that there must be an irreducible summand of highest weight $\gamma-\a_2 = 1(b-1)2b$.
 
 Next we claim that there is a composition factor of highest weight $\gamma- (\a_2+\a_3+\a_4)$ appearing with multiplicity $3$.
 
We begin by considering possibilities for highest weights $\eta$ of composition factors  which lie above
this weight.  From the above we see that $\eta=\gamma-\a_3$ is not possible.  Now $ \gamma- (\a_2+\a_3)$ can only arise from $\xi-(\b_1^0+\b_9^0)$,  $\xi-(\b_1^0+\b_2^0+\b_3^0)$, and $\xi-(\b_7^0 +\b_8^0 +\b_9^0)$ and hence has multiplicity $4$.  This weight appears  twice in $\gamma$ and once in $\gamma-\a_2$, hence there does exist a composition factor with highest weight $ \gamma - (\a_2+\a_3)$. Similarly, we check that the same assertion holds for $\g-\a_4$. 

 Similarly $ \gamma - (\a_2+\a_4)$ can only arise from $\xi-(\b_2^0+\b_9^0)$ and $\xi-(\b_2^0+\b_3^0+\b_4^0)$ and hence has multiplicity $2$ in the restriction of $\xi$.  But it already occurs in $\gamma$ and $\gamma-\a_2$. And $ \gamma - (\a_3+\a_4)$ only arises from $\xi-(\b_1^0+\b_2^0)$ and hence occurs with multiplicity $2$ which it also does in $\gamma$.
 
 Now $\nu =\gamma- (\a_2+\a_3+\a_4)$ arises from $\xi-(\b_1^0+\b_2^0+\b_3^0+\b_4^0)$, $\xi-(\b_1^0 +\b_2^0+\b_9^0)$, $\xi-(\b_1^0+2\b_2^0+\b_3^0)$,  $\xi-(\b_6^0+\b_7^0+\b_8^0+\b_9^0)$,  $\xi-(\b_2^0+\b_7^0+\b_8^0+\b_9^0)$, and $\xi-(\b_2^0+\b_3^0+\b_4^0+\b_5^0 +\b_6^0 +\b_7^0)$ and these have multiplicities $2$, $2$, $3$, $1$, $1$, and $1$, respectively.  So $\nu$ occurs with multiplicity at least $10$.  On the other hand
 $\nu$ appears with multiplicity $4$ in $\gamma$,  multiplicity 2 in $\gamma-\a_2$, and multiplicity $1$ in $ \gamma - (\a_2+\a_3)$ so this establishes the claim.  We complete the analysis of these cases by using \cite[Table 6.2, p. 56]{MF} to see that at most  one copy of $\nu$ can arise from $V^1$.
 
 Finally we consider the case where $\l = a\l_8 + x\l_{10}$.  Set $b = a-1$ as before, so that $b = 1,2,3$ or $4$. Then $V^2(Q_Y) \supseteq (\l_7^0 + b\l_8^0 + \l_9^0) = \xi$.     Then
 $$(\l_7^0 + b\l_8^0 + \l_9^0) \subseteq  \l_7^0  \otimes S^b(\l_8^0) \otimes \l_9^0.$$
We note that $\l_8^0 = 2\l_9^0 - \b_9^0$ and $\l_7^0 = 3\l_9^0 - 2\b_9^0-\b_8^0$ so the above restrictions imply that $\l_8^0$ and $\l_7^0$ restrict to $1010$ and $0020$ for $L_X'$, respectively, and the maximal vector of $\l_7^0 + b\l_8^0 + \l_9^0$ affords weight $\gamma = b1(b+2)0$ of $S$-value $2b+3$. 

Assume first that $a=2$. Here 
there is an  irreducible in $V^2(Q_Y)$ of highest weight $\xi = \l_7^0+\l_8^0+\l_9^0$. 
One checks that there is also an irreducible in $V^2(Q_Y)$ of highest weight $2\l_7^0$  afforded by 
$\l - 2\b_8 - 2\b_9 -\b_{10}$.  
Using Magma, we check that between them, the restrictions of these irreducibles to $L_X'$ contain
$(1101)^4$. However, $V^1 = 0002+0201+2020+1100$ (see \cite[Table 6.2]{MF}), so this gives the usual contradiction using Proposition \ref{induct}.
 
 Hence $a\ge 3$ (so that $b\ge 2$). 
 We will establish the lemma in this case by showing that there is a 
 repeated composition factor of highest weight 
 $\gamma-(\alpha_1+\alpha_2+\alpha_3)$ in
 $((\lambda_7^0+b\lambda_8^0+\lambda_9^0) \oplus
 (2\lambda_7^0+(b-1)\lambda_8^0))\downarrow L_X'$. Note that the second 
 summand appears in $V^2(Q_Y)$ by the same argument as given in the $a=2$ case.
 
 Using the form of $\xi= \l_7^0 + b\l_8^0 + \l_9^0$ and the restriction of the roots $\b_1^0, \ldots, \b_9^0$ one checks that 
$\nu =\gamma +\a_1-\a_3 = (b+2)1b1$ is a dominant weight  in the restriction of $\xi$ which is not subdominant to any other weight occurring in this restriction. Also $\nu$ 
 is afforded only by $\xi-\b_7$ and hence has multiplicity $1$. Since $\gamma$ is not subdominant to $\nu$ it follows that both $\gamma$ and $\nu$ occur as composition factors of  $\xi \downarrow L_X'$ of multiplicity 1.   
 
The weight $\gamma-\a_1$ only arises from the restriction of $\xi-\b_8^0$ and so the weight space has dimension $1$ and is contained in $\gamma$.

 We next note that $\gamma-\a_3$  only occurs as the restriction of $\xi - (\b_7^0 +\b_8^0)$ and  $\xi - (\b_5^0+\b_6^0+\b_7^0)$ and the corresponding weight spaces have dimensions $2$ and $1$, respectively.  Therefore $\gamma-\a_3$ appears with multiplicity $3$ in $\xi$.  It appears in both $\gamma$ and $\nu$ with multiplicity $1$ and so $\gamma -\a_3 = b2b1$ appears as a composition factor of $\xi$ with multiplicity $1$. 
 
 Consider $\gamma-(\a_1+\a_2)$.  This appears with multiplicity 2 in $\xi$ since it only arises from $\xi-(\b_8^0+\b_9^0)$.  It also occurs with  multiplicity $2$ in $\gamma$ so $\gamma-(\a_1+\a_2)$ does not occur as a composition factor.  The weight $\gamma-(\alpha_2+\alpha_3)$ arises from the weights 
 $\xi-(\beta_7^0+\beta_8^0+\beta_9^0)$,
 $\xi - (\beta_5^0+\beta_6^0+\beta_7^0+\beta_9^0)$ and $\xi -  (\beta_3^0+\beta_4^0+\beta_5^0+\beta_6^0+\beta_7^0)$, and hence has 
 multiplicity  6 in $\xi$. It occurs in $\gamma$, $\nu$ and $\gamma-\alpha_3$, with respective multiplicities $2,2,1$ and so affords a summand of $V^2$.

 There is a composition factor of highest weight $\gamma-(\a_1+\a_3)$.  Indeed this arises from $\xi -(\b_7^0 +2\b_8^0)$ and from $\xi - (\b_5^0+\b_6^0+\b_7^0+\b_8^0)$ and hence the weight space has multiplicity $4$ in $\xi$.  It occurs with multiplicity $1$ in each of $\gamma, \gamma-\a_3,$ and $\nu$ and so we conclude that there is a composition factor of multiplicity 1 and highest weight
 $\gamma-(\a_1+\a_3)$. 

Finally we consider the weight $\gamma-\alpha_1-\alpha_2-\alpha_3$. 
This weight occurs with multiplicity 4 in $\g$, multiplicity 1 in $\gamma-\alpha_1-\alpha_3$, multiplicity 2 in 
$\gamma-\alpha_3$,  multiplicity 1 in $\gamma-\alpha_2-\alpha_3$,  and multiplicity 2 in
$\gamma-\alpha_3+\alpha_1$, while we can obtain it from the weights $\xi-\beta_5-\cdots-\beta_9$, 
$\xi-\beta_7-2\beta_8-\beta_9$ and $\xi-\beta_5-\cdots-\beta_8$, and therefore there is  a summand of $\xi$ of this highest weight.
We now argue that  $\gamma-\alpha_1-\alpha_2-\alpha_3$ also occurs
 as a summand of $(2\lambda_7^0+(b-1)\lambda_8^0)\downarrow L_X'$.
 Set $\tau = 2\lambda_7^0+(b-1)\lambda_8^0$ and $\eta = \tau\downarrow L_X'$, so $\eta = ((b-1)0(b+3)0)$.
 Note that $\eta-\alpha_3+\alpha_1$ occurs as the restriction of $\tau - \beta_7$ and no other weight in $\tau$ restricts to 
$\eta-\alpha_3+\alpha_1$.
 Since this weight is not subdominant to $\eta$ and one checks that there is no
 other weight in $\tau$ whose restriction lies above 
 $\eta-\alpha_3+\alpha_1$, we have that $\eta-\alpha_3+\alpha_1$ also affords a summand. 
Then we check that $\eta-\alpha_3$ arises from the 
 restriction of $\tau-\beta_5^0 -\beta_6^0-\beta_7^0$ and from 
 $\tau-\beta_7^0-  \beta_8^0$, and so has multiplicity 3 in $\tau$, but only occurs with 
 multiplicity 1 in each of $\eta$ and $\eta-\alpha_3+\alpha_1$.
 Hence we have a summand $\eta-\alpha_3$, which is of the form 
$((b-1)1(b+1)1)$, that is, the same as the summand $\gamma-\alpha_1-\alpha_2-\alpha_3$ occurring in $\xi$. Thus we have two summands $((b-1)1(b+1)1)$ in $V^2$.
 Using Table 6.2 of \cite{MF}, and recalling that there the 
 embedding is given by $X$ acting on the natural $A_9$-module with highest weight $\omega_2$, while here we have the dual action, we see that the list of summands given is precisely the summands of $V^1$ in our case.
 It is straightforward to check that no summand of the form $((b-1)1(b+1)1)$ can arise from $V^1$, for $a=b-1$, and $a=3,4,5$. This final contradiction completes the proof. 
 \hal

\vspace{4mm}
We have now completed the proof of Theorem \ref{A5<C10}.

\section{The case $X = A_{2s+1} < Y=C_n$ with $\d = \o_{s+1}$, $s\ge 4$ even}

In this section we prove Theorem \ref{mainthm} in the case where $X = A_{2s+1}$ with $\d = \o_{s+1}$ with $s$ even. In this case Lemma \ref{setup} shows that $W = V_X(\d)$ is a symplectic module, so that $X < Y = Sp(W) = C_n$. Also the case where $s=2$ was covered by Theorem \ref{A5<C10}, so we assume that $s \ge 4$.
 
\begin{thm}\label{A(2s+1)<C} Let $X = A_{2s+1}$  with $s \ge 4$ even.  Assume $X$ is embedded in $Y$ of type $C_n$ via the representation with highest weight $\d = \omega_{s+1}$.  If $V = V_Y(\l)$
where $\l \ne 0,\l_1$,  then $V \downarrow X$ is MF if and only if $\l =\l_2$, $2\l_1$,  $\l_3 (s = 4,6)$, or
$3\l_1 (s = 4)$.
\end{thm}

Note that if $\l$ is as in the conclusion of the theorem, then $V_Y(\l)\downarrow X$ is indeed MF: this follows directly from Theorem 1 of \cite{MF}.

Now assume the hypotheses of Theorem \ref{A(2s+1)<C}, and suppose $V_Y(\l)\downarrow X$ is MF.
We have $L_X' = A_{2s}$ embedded in $L_Y'$ which is of type $A$ via the representation $\d' = \omega_{s+1}$. Also write $\d''=\o_s$ (for $L_X'$). 
 Theorem 1 of \cite{MF} shows that $\mu^0$ or its dual is one of the the following: 
\begin{equation}\label{listos1}
\begin{array}{l}
0, \l_1^0, \\
\l_1^0+\l_{r_0}^0, \\
2\l_1^0, \\
\l_2^0, \\
\l_3^0\, (s = 4,6), \\
\l_4^0\, (s = 4), \\
3\l_1^0\,(s = 4). 
\end{array}
\end{equation}

We shall work through the list of possibilities in (\ref{listos1}).

\begin{lem}\label{l1ln} We have $\mu^0  \ne 0, \l_1^0$, or $\l_{r_0}^0$. 
\end{lem}

\pf  Suppose $\mu^0 = 0$, then $\l = c\l_n$ with $c>0$. Then $V^2(Q_Y) = 0 \cdots 02$ (afforded by $\l-\b_n$), while  
 $V^3(Q_Y) \supseteq 0\cdots 020$ (afforded by $\l-2\b_{n-1}-2\b_n$). Hence $V^2 = S^2(\d'')$ of $S$-value 2, and 
$V^3 \supseteq S^2(\wedge^2(\d''))-\wedge^4(\d'')$ where $\d''=\o_s$. For $s=4$, Magma shows that $V^3$ contains $(10101010)^2$. This weight is equal to $2(2\d''-\a_s)-\a_{s-2}-2\a_{s-1}-2\a_s-2\a_{s+1}-\a_{s+2}$. Hence it has multiplicity at least 2 in $V^3$ for all values of $s\ge 4$. As it has $S$-value 4, this is a contradiction.

  Next suppose  $\mu^0 = \l_1^0$ so that $\l = \l_1 + x\l_n$.  By the hypothesis of Theorem \ref{A(2s+1)<C} we have $\l \ne \l_1$, so $x > 0$.  Then $V^2(Q_Y) \supseteq 0\cdots 01 + 10 \cdots 02 = (10\cdots 0) \otimes (0\cdots 02)$, and
  $V^2 \supseteq \d' \otimes S^2(\d'')$.  A Magma check for $s = 4$, shows that $V^2 \supseteq (\om_2 + \om_4 + \om_7)^2$.  For $s = 4 $, we have $(\om_2 + \om_4 + \om_7) = (\d' + 2\d'') - (\a_{s-1} + 2\a_s + 2\a_{s+1} + \a_{s+2})$, so it follows that $V^2 \supseteq  (\om_{s-2} + \om_s + \om_{s+3})^2$ for
  all even $s \ge 4$.  Now $V^1 = \om_{s+1}$ cannot contribute such a composition factor, so we conclude
  that $V\downarrow X$ is not MF.
  
  Finally assume  $\mu^0 = \l_{r_0}^0$ so that $\l = \l_{n-1} + x\l_n$.  Then $V^2(Q_Y) \supseteq 0 \cdots 011 =
  (0\cdots 010 \otimes 0 \cdots 01) - (0\cdots 0100)$ and $V^2 \supseteq (\wedge^2(\d'') \otimes \d'') - \wedge^3(\d'')$. Assume $s \ge 6$.  For $s = 6$, Magma shows that $\wedge^2(\d'') \otimes \d'' \supseteq (\om_1+\om_6+\om_{11})^3$ and that  this composition factor appears with multiplicity 1 in  $\wedge^3(\d'')$.
  Therefore for $s = 6$, $V^2 \supseteq (\om_1+\om_6+\om_{11})^2$ and  as $V^1$ cannot contribute such a composition factor we see that $V\downarrow X$ is not MF when $s =6$.
  
For $s > 6$ and $\mu = 3\om_s -  (\a_{s-4}+ 2\a_{s-3} +3\a_{s-2} +4\a_{s-1} +5\a_s+4\a_{s+1}+3\a_{s+2} + 2\a_{s+3}+ \a_{s+4}) $, we have $\mu = \om_{s-5}+\om_s +\om_{s+5}$ and $V^2 \supseteq \mu^2$.  Again we conclude that $V\downarrow X$ is not MF. 
 
 Finally assume $s = 4$.  If $x > 0$, then $V^2(Q_Y) \supseteq 0\cdots 03 + 0\cdots 011 = 0\cdots 02 \otimes 0 \cdots 01$. Then $V^2 \supseteq S^2(\om_4) \otimes \om_4 \supseteq (\om_4 + \om_8)^3$.  This is a contradiction as $V^1 = \om_4$.  Now assume $x = 0$.  Here $V^2(Q_Y) = 0\cdots 011$ so that $V^2 = (\wedge^2(\om_4) \otimes \om_4) -\wedge^3(\om_4)$.  One checks that $V^3(Q_Y) \supseteq 0 \cdots 0110$
 which is afforded by $\l - \b_{n-2} - 3\b_{n-1} - 2\b_n$.  Then $V^3 \supseteq (\wedge^3(\om_4) \otimes \wedge^2(\om_4)) - (\wedge^4(\om_4) \otimes \om_4)$ and this contains $(2\om_1 + \om_5 + \om_6 + \om_7)^2$.  As $S(V^2) = 3$, this is a contradiction. \hal
 
 \begin{lem}\label{l1ln} We have $\mu^0  \ne \l_1^0+\l_{r_0}^0$. 
\end{lem}

\pf  Suppose $\mu^0  = \l_1^0+\l_{r_0}^0$. Here $V^1 = \o_{s+1}\otimes \o_s-0$, of $S$-value 2. Also
$V^2(Q_Y) \supseteq (10\cdots 011) + (0\cdots  02)$, which is equal to $(10\cdots 0 \otimes 0\cdots 011) - 0\cdots 010$, and so 
\[
V^2 \supseteq (\d'\otimes \wedge^2\d'' \otimes \d'') - (\d'\otimes \wedge^3\d'')-\wedge^2\d''.
\]
This contains $(\o_{s-2}+2\o_s+\o_{s+3})^2$ (obtained for arbitrary rank as $(\o_{s-1}+\o_s+2\o_{s+1}-\a_{s-1}-\a_s-2\a_{s+1}-\a_{s+2})^2$). This has $S$-value 4, a contradiction. \hal

\begin{lem}\label{l2m} If $\mu^0  = 2\l_1^0$ or $2\l_{r_0}^0$, then $\l = 2\l_1$. 
\end{lem}

\pf Suppose first that $\mu^0  = 2\l_{r_0}^0$, so that $\l = 2\l_{n-1}+x\l_n$. Here $V^1 = S^2(\d'')$, while $V^2(Q_Y)$ contains $0\cdots 012 = (0\cdots 010)\otimes (0\cdots 02)-(0\cdots 0101)$. Hence 
\[
\begin{array}{ll}
V^2 & \supseteq (\wedge^2\d'' \otimes S^2\d'') - (\wedge^3\d''\otimes \d'') + \wedge^4\d'' \\
       & \supseteq (\o_{s-1}+2\o_s+\o_{s+1}-\a_{s-2}-2\a_{s-1}-2\a_s-2\a_{s+1}-\a_{s+2})^2 \\
       & = (\o_{s-3}+2\o_s+\o_{s+3})^2,
\end{array}
\]
of $S$-value 4, a contradiction.

Now suppose $\mu^0  = 2\l_1^0$, so $\l = 2\l_1+x\l_n$. If $x=0$ then $\l=2\l_1$ as in the conclusion, so assume that $x>0$. Then $V^1=S^2\d'$, while $V^2(Q_Y) \supseteq (10\cdots 01)+(20\cdots 02)$. Hence $V^2 \supseteq (S^2\d'\otimes S^2\d'')-0$. For $s=4$, Magma shows that this contains $(11000011)^2$, and so the usual parabolic argument shows that for any $s\ge 4$, $V^2$ contains a repeated composition factor of $S$-value at least 4, a contradiction. \hal

\begin{lem}\label{l3m} If $\mu^0  = \l_2^0$ or $\l_{r_0-1}^0$, then $\l = \l_2$. 
\end{lem}

\pf Suppose $\mu^0  = \l_{r_0-1}^0$. Then $V^1=\wedge^2\d''$. Also $V^2(Q_Y) \supseteq 0\cdots 0101$, so $V^2$ contains $(\wedge^3\d''\otimes \d'')-\wedge^4\d''$. For $s=4$ this contains $(20000101)^2$, and so the parabolic argument shows that $V^2$ contains a repeated composition factor of $S$-value at least 4, for any $s\ge 4$, a contradiction. 

Now let $\mu^0  = \l_2^0$, so $\l = \l_2+x\l_n$. If $x=0$ then $\l=\l_2$ as in the conclusion, so assume that $x>0$. Then $V^1=\wedge^2\d'$. Also $V^2(Q_Y)$ contains $(010\cdots 02)+(10\cdots 01) = (010\cdots 0)\otimes (0\cdots 02)$,  so $V^2 \supseteq \wedge^2\d'\otimes S^2\d''$. For $s=4$ Magma shows that this contains $(10011001)^2$, and now we obtain the usual $S$-value contradiction using a parabolic argument. \hal 

\begin{lem}\label{l4m} Suppose $s \in \{4,6\}$ and  $\mu^0  = \l_3^0$ or $\l_{r_0-2}^0$. Then $\l = \l_3$. 
\end{lem}

\pf Assume $\mu^0  = \l_{r_0-2}^0$. Then $V^1 = \wedge^3\d''$ and $V^2(Q_Y) \supseteq 0\cdots 01001$, so that 
$V^2 \supseteq (\wedge^4\d''\otimes \d'')-\wedge^5\d''$. But Magma shows that this has a repeated composition factor of $S$-value at least 5, a contradiction. 

Now assume $\mu^0  = \l_3^0$, so $\l = \l_3+x\l_n$. If $x=0$ then the conclusion holds, so take $x>0$. Then $V^1 = \wedge^3\d'$, while we see in the usual way that $V^2 \supseteq \wedge^3\d'\otimes S^2\d''$, which Magma shows to contain a repeated composition factor of $S$-value 5, a contradiction. \hal 

\begin{lem}\label{l4m} For $s=4$ we have $\mu^0  \ne \l_4^0$ or $\l_{r_0-3}^0$. 
\end{lem}

\pf Suppose $s=4$. If $\mu^0  = \l_{r_0-3}^0$ then $V^1=\wedge^4\d''$ and $V^2(Q_Y) \supseteq 0\cdots 010001$. Hence $V^2 \supseteq (\wedge^5\d'' \otimes \d'')-\wedge^6\d''$, and Magma shows that this contains $(30000111)^2$, of $S$-value 6, a contradiction. 

Now suppose $\mu^0  = \l_4^0$, so that $\l = \l_4+x\l_n$. If $x=0$ then $\l = \l_4$ and $V_Y(\l)\downarrow X = \wedge^4\o_5 - \wedge^2\o_5$ (for $X = A_9$), which Magma shows to be non-MF. Hence $x>0$. 
Then $V^1 = \wedge^4\d'$ while $V^2 \supseteq \wedge^4\d' \otimes S^2\d''$, which contains $(11002011)^2$, of $S$-value 6, a contradiction. \hal

\begin{lem}\label{l5m} If $s=4$ and $\mu^0  = 3\l_1^0$ or $3\l_{r_0}^0$, then $\l = 3\l_1$. 
\end{lem}

\pf Suppose $\mu^0  = 3\l_{r_0}^0$. Then $V^1 = S^3\d''$ while $V^2(Q_Y) \supseteq 0\cdots 013$ which is equal to $(0\cdots 01)\otimes (0\cdots 04)-(0\cdots 05)$. Hence $V^2 \supseteq (\d''\otimes S^4\d'')-S^5\d''$, and this contains $(11100101)^2$, of $S$-value 5, a contradiction. 

Now assume $\mu^0  = 3\l_1^0$, so $\l = 3\l_1+x\l_n$. If $x=0$ then the conclusion holds, so take $x>0$. Then $V^1 = S^3\d'$, while we see in the usual way that $V^2 \supseteq (S^3\d'\otimes S^2\d'')-\d'$, which Magma shows to contain $(10100210)^2$ of $S$-value 5, a contradiction. \hal 

\vspace{4mm}
We have now covered all the possibilities for $\mu^0$ given by (\ref{listos1}), so the proof of Theorem \ref{A(2s+1)<C} is complete.

\section{The case $X = A_{2s+1} < Y=D_n$ with $\d = \o_{s+1}$, $s\ge 3$ odd}

 In this section we complete the proof of Theorem \ref{mainthm} for the case where $X=A_{2s+1}$ with $\d=\o_{s+1}$ by proving

\begin{thm}\label{A_{2s+1}<Dn} Assume $X = A_{2s+1}$, with $s \ge 3$ odd, is embedded in $Y = D_n$ via $\d = \om_{s+1}$,  and let $V = V_Y(\l)$
where $\l \ne 0,\l_1$.  Then $V \downarrow X$ is MF if and only if $\l$ is one of the following weights:
\[
\begin{array}{l}
\l_2,\,2\l_1\,(\hbox{both with any }s), \\
\l_3,\,\,3\l_1 \,(\hbox{both with }s=3,5), \\
\l_4,\,4\l_1,\,\l_{34},\,\l_{35}\,(\hbox{all with }s=3).
\end{array}
\]
\end{thm}

Note that when $s=3$ we have $Y = D_{35}$ and $\l_{34},\,\l_{35}$ are spin modules for $Y$.

We begin by proving that $V_Y(\l)\downarrow X$ is indeed MF for the weights $\l$ in the conclusion. 

\begin{lem}\label{MFAD}  If $X = A_{2s+1} < Y=D_n$ via $\d=\o_{s+1}$ with $s\ge 3$ odd, and $\l$ is one of the weights listed in the conclusion of Theorem $\ref{A_{2s+1}<Dn}$, then  $V_Y(\l) \downarrow X$ is MF.
\end{lem}

\pf For the cases with arbitrary $s$, this follows from Theorem 1 of \cite{MF}. The cases $\l = \l_3,\,3\l_1$ $(s=3,5)$ 
and $\l = \l_4,\,4\l_1$ ($s=3$) are Magma computations (using Lemma \ref{symwed}): for example,  
for $\l = 4\l_1$, Lemma \ref{symwed} gives that for $Y=D_{35}$ we have $4\l_1 = S^4(\l_1)-S^2(\l_1)$, and a Magma check shows that the restriction to $X = A_7$, namely $S^4(\o_4)-S^2(\o_4)$, is MF.

This leaves the case where $s=3$ and $\l = \l_{34}$ or $\l_{35}$, a spin module for $Y=D_{35}$.
This case requires a great deal more argument. Since $V_Y(\l_{34}) = V_Y(\l_{35})^*$, we need to consider just one of these two weights.

Let $X = A_7$ and let $X<Y = D_{35}$ via $\d = \o_4$. As in  Section 2, we have $L_X' = A_6$ embedded in the Levi factor $L_Y' = A_{34}$ of the parabolic subgroup $P_Y = P_{35}$ of $Y$, and centralized by the 1-dimensional torus $T = \{T(c) : c \in K^*\}$  of $X$, where $T(c) = \prod_1^7 h_{\a_i}(c^i)$. Then $\a_7(T(c)) = c^8$ and $\a_i(T(c)) = 1$ for $i\le 6$.

Take $\l = \l_{35}$ and let $V = V_Y(\l)$. We shall find all the composition factors of the restriction $V\downarrow X$. 

If $J$ is the natural module for $X$, then $T(c)$ acts as $(c,c,c,c,c,c,c,c^{-7})$ on $J$ and hence it acts on $\wedge^4(J)$ with two diagonal blocks  $c^4I_{35}$ and $ c^{-4}I_{35}$.  Let $A$ and $B$ be the respective singular subspaces. Choosing dual bases for $A$ and $B$ we can choose a system of fundamental roots $ \nu_1, \cdots, \nu_{35} $  for $Y$ such that  
$\nu_i(T(c)) = 1$ for $1\le i \le 34$ and $\nu_{35}(T(c)) = c^8$.    Now $\l = \l_{35} = \frac{1}{2}(\nu_1 + 2\nu_2 + \cdots  + 33\nu_{33} + \frac{33}{2}\nu_{34} +\frac{35}{2}\nu_{35})$ so that $\l(T(c)) = c^{70}$. 

The levels of $V$ for $L_X'$ are as follows: $V^1 = 000000$, and $V^i = \wedge^{2i-2}\om_3$ for $i>0$. We shall produce composition factors of $V\downarrow X$ corresponding to these levels for $1\le i\le 8$. Start at $i=1$: since $V^1 = 000000$, there must be a composition factor of $V \downarrow X$ of the form $000000x$ for some $x$. Observe that 
$\l(T(c)) = c^{7x}$ so from the above we have $x = 10$ and $V\downarrow X$ has a composition factor of highest weight 
\[
\mu_1 = 000000(10).
\]

We use similar arguments to produce further composition factors corresponding to succeeding levels. Consider $V^2 = \wedge^2\o_3$. This has a summand $010100$ that does not arise from $V^1$; hence $V\downarrow X$ has a composition factor $\mu_2=010100z$. The weights of $\mu_2$ have the form $\l-\a_7-\psi$, so we see that $\mu_2(T(c)) = c^{7x-8} = c^{6+7z}$. As $x=10$ this gives $z=8$, so $V\downarrow X$ has a  further composition factor 
\[
\mu_2 = 0101008.
\]
Next consider $V^3 = \wedge^4\o_3$. Magma gives the composition factors of $V^3$, and those not arising from $\mu_1$ or $\mu_2$ are $030001$, $000300$ and $110110$. Arguing as above, we see that $V\downarrow X$ has  further composition factors
\[
\begin{array}{l}
\mu_3 = 0300016, \\
\mu_4 = 1101106, \\
\mu_5 = 0003006.
\end{array}
\]
We now repeat this procedure for the levels $V^4,\ldots, V^8$, and obtain further composition factors $\mu_i$ ($6\le i \le 36$) of $V\downarrow X$, as listed in Table \ref{comps}.

\begin{table}[h!] 
\caption{}\label{comps}
\[
\begin{array}{|c|c|}
\hline
\hbox{level} & \hbox{composition factors} \\
\hline
V^4 & \mu_6 = 2110114 \\
       & \mu_7 = 1301005 \\
       & \mu_8 = 3000304 \\
       & \mu_9 = 1011204 \\
       & \mu_{10} = 0202014 \\
\hline
V^5 & \mu_{11} =1211103 \\
       & \mu_{12} =1010402 \\
       & \mu_{13} =4010122 \\
       & \mu_{14} =3120013 \\
       & \mu_{15} =2101212 \\
       & \mu_{16} =2021022 \\
       & \mu_{17} =0120212 \\
       & \mu_{18} =3200203 \\
\hline
V^6 & \mu_{19} =0040030 \\
       & \mu_{20} =1210301 \\
       & \mu_{21} =1130111 \\
       & \mu_{22} =0200410 \\
       & \mu_{23} =6000040 \\
       & \mu_{24} =3031002 \\
       & \mu_{25} =0000600 \\
       & \mu_{26} =3111111 \\
       & \mu_{27} =5101021 \\
       & \mu_{28} =5020102 \\
       & \mu_{29} =2020220 \\
       & \mu_{30} =4002030 \\
\hline
V^7 & \mu_{31} =7010020 \\
       & \mu_{32} =7002001 \\
       & \mu_{33} =5012010 \\
       & \mu_{34} =3030200 \\
       & \mu_{35} =1050010 \\
\hline
V^8 & \mu_{36} =9000100 \\
\hline
\end{array}
\]
\end{table}

Finally, we use Magma once more to compute and add the dimensions of the composition factors $\mu_i$ listed so far, and find  that 
\[
\sum_{i=1}^{36} \dim V_X(\mu_i) = 2^{34} = \dim V.
\]
It follows that  $V\downarrow X = \sum_1^{36} V_X(\mu_i)$ is MF, completing the proof. \hal

Now assume the hypotheses of Theorem \ref{A_{2s+1}<Dn}, and assume that  $V_Y(\l)\downarrow X$ is MF. 
We have $L_X' = A_{2s}$ embedded in $L_Y' = A_{n-1}$ via the representation $\d' = \om_{s+1}.$  It will be convenient to write $\d'' = \om_s$.  Note that $r_0 = n-1$ and also that $n = \frac{1}{2}{{2s+2}\choose {s+1}}$.

The main result
of \cite{MF} shows that $\mu^0$ or its dual is one of the the following weights:
\begin{equation}\label{listad}
\begin{array}{l}
0,\,\l_1^0,\,2\l_1^0,\,\l_2^0,\,\l_1^0+\l_{r_0}^0, \\
3\l_1^0,\,\l_3^0\,(s=3,5 \hbox{ in each case}), \\
4\l_1^0,\,5\l_1^0,\,\l_4^0,\,\l_5^0,\,\l_6^0,\,\l_1^0+\l_2^0\,(s=3 \hbox{ in each case}).
\end{array}
\end{equation}
We shall deal with each of these possibilities in the following lemmas.

\begin{lem}\label{ad1} Suppose $\mu^0=0$. Then $s=3$, $Y=D_{35}$ and $\l$ is a spin module $\l_{35}$.
\end{lem}

\pf Here $\l = x\l_n$ with $x\ge 1$. We have $V^2(Q_Y) = \l_{r_0-1}^0$ (afforded by $\l-\g_1$), so $V^2 = \wedge^2\d'' = \wedge^2\o_s$, of $S$-value 2. 

Suppose $x>1$. Then $V^3(Q_Y)$ contains $0\cdots 020$ (afforded by $\l-2\g_1$), so $V^3 \supseteq S^2(\wedge^2\o_s)-\wedge^4\o_s$. For $s=5$, Magma shows that this contains $(1010001010)^2$, and the usual parabolic argument then implies that for any $s\ge 5$, $V^3$ has a repeated summand of $S$-value at least 4, a contradiction. Hence $s=3$ and $V^2 = \wedge^2\o_3 = 010100 + 000001$. As $x>1$, one checks that $V^3(Q_Y)$ has  a further summand afforded by $\l-\b_{n-3}-2\b_{n-2}-\b_{n-1}-2\g_1$, and the restriction of this to $L_X'$ is equal to $\wedge^4\o_s$. Hence 
$V^3 \supseteq S^2(\wedge^2\o_s) \supseteq (010101)^3$. Only one of these summands can come from $V^2$, a contradiction. 

Hence $x=1$ and $\l=\l_n$, a spin module for $Y=D_n$. Here $V^3(Q_Y)$ has  a  summand afforded by 
$\l-\b_{n-3}-2\b_{n-2}-\b_{n-1}-2\g_1$, and the restriction of this to $L_X'$ contains $\wedge^4\o_s$.
For $s=7$, Magma shows that this contains $(01010001000001)^3$, and so the parabolic argument implies that for any $s\ge 7$, $V^3$ has a repeated summand of $S$-value at least 4, a contradiction. For $s=5$, Magma gives
\[
\begin{array}{l}
V^2 = \wedge^2\o_5 = 0001010000+0100000100+0000000001, \\
V^3 \supseteq \wedge^4\o_5 \supseteq (0001010001)^3,
\end{array}
\]
and only one of the summands $0001010001$ in $V^3$ can come from $V^2$. Hence $s=3$ and $\l = \l_{35}$, as in the conclusion. \hal 

\begin{lem}\label{ad2} We have $\mu^0 \ne \l_1^0$.
\end{lem}

\pf Suppose $\mu^0 = \l_1^0$, so that $\l = \l_1+x\l_n$ with $x>0$ (as $\l \ne \l_1$). Here $V^1 = \o_{s+1}$ of $S$-value 1, while $V^2(Q_Y) = 10\cdots 010 + 0\cdots 01 = (10\cdots 0)\otimes (0\cdots 010)$, so that $V^2 = \o_{s+1}\otimes \wedge^2\o_s$. For $s=5$ this contains $(0100010100)^2$, and so the parabolic argument shows that for any $s\ge 5$, $V^2$ has a repeated summand of $S$-value at least 3, a contradiction. 

Finally assume $s=3$. Now $V^3(Q_Y)$ contains a summand $10\cdots 01000$ afforded by 
 $\l-\b_{n-3}-2\b_{n-2}-\b_{n-1}-2\g_1$, and so 
\[
V^3 \supseteq (\o_4\otimes \wedge^4\o_3)-\wedge^3\o_3 \supseteq (120011)^2,
\]
of $S$-value 5, a contradiction (as $S(V^2) = 3$). \hal

\begin{lem}\label{ad3} Suppose $\mu^0 = (\l_1^0)^* = \l_{r_0}^0$. Then $s=3$ and $\l=\l_{34}$.
\end{lem}

\pf Here $\l = \l_{n-1} + x\l_n$.  First assume $x = 0$.  If $s = 3$, then
$Y = D_{35}$  and $\l = \l_{34}$, as in the conclusion.   Now suppose $s \ge 5$. Then $V^2(Q_Y) = 0 \cdots 0100$ and $V^3(Q_Y) = 0 \cdots 010000$. Therefore $V^3 \supseteq \wedge^5(\om_5)$.  For $s = 5$ a Magma computation shows that $V^3$ contains $(1100100110)^2$ which is afforded by $5\omega_5-(\alpha_2+3\alpha_3+5\alpha_4+7\alpha_5+5\alpha_6+3\alpha_7+\alpha_8)$.  This is a contradiction as $S(V^2) = 3$.  A parabolic argument now shows that we also obtain a contradiction for $s > 5$.  Hence we may assume that $x > 0$.

If $x=1$ then $\l = \l_{n-1}+\l_n$ and $V_{D_n}(\l) = \wedge^{n-1}(\l_1)$. However, $\wedge^{n-1}(\o_{s+1})$ is not MF for $X = A_{2s+1}$ by \cite{MF}.

Hence $x>1$. Now $V^1 = \o_s$ and $V^2(Q_Y) = 0\cdots 011 + 0\cdots 0100 = (0\cdots 010)\otimes (0\cdots 01)$. Hence $V^2 = \wedge^2\o_s\otimes \o_s$. For  $s\ge 5$ this has a repeated summand 
$3\omega_s-\alpha_{s-2}-2\alpha_{s-1}-3\alpha_s-2\alpha_{s+1}-\alpha_{s+2}$
of $S$-value at least 3.
Hence $s=3$. Here $V^3(Q_Y)$ has summands $0\cdots 021$ and $0\cdots 0110$ afforded by $\l-2\g_1$ and $\l-\b_{n-2}-\b_{n-1}-2\g_1$, and the sum of these is $(0\cdots 020)\otimes (0\cdots 01)$. Consequently $V^3 \supseteq (S^2(\wedge^2\o_3)-\wedge^4\o_3)\otimes \o_3$, and this contains $(201020)^2$, of $S$-value 5. This is a contradiction, as $S(V^2)=3$. \hal

\begin{lem}\label{ad4} If $\mu^0 = 2\l_1^0$, then $\l=2\l_1$.
\end{lem}

\pf Suppose $\mu^0 = 2\l_1^0$, so that $\l = 2\l_1+x\l_n$. Assume $x>0$. Then $V^1 = S^2\o_{s+1}$, while 
$V^2(Q_Y) = (20\cdots 0) \otimes (0\cdots 010)$, so that $V^2 = S^2\o_{s+1}\otimes \wedge^2\o_s$. For $s=3$ this contains a repeated summand of $S$-value at least 4 (namely $(110011)^2$), and so the parabolic argument shows that the same holds for any $s\ge 3$, a contradiction. \hal

\begin{lem}\label{ad5} We have $\mu^0 \ne (2\l_1^0)^*$.
\end{lem}

\pf Suppose $\mu^0 = (2\l_1^0)^*$, so that $\l = 2\l_{n-1}+x\l_n$. Here $V^1 = S^2\o_s$, while $V^2(Q_Y) 
\supseteq  0\cdots 0101$, so that $V^2 \supseteq (\wedge^3\o_s \otimes \o_s) - \wedge^4\o_s$. For $s=5$ this contains a repeated summand of $S$-value at least 4 (namely $(1001100001)^3$), and so the parabolic argument shows that 
a similar conclusion holds for any $s\ge 5$, a contradiction. Finally, for $s=3$ we have $V^1 = S^2\o_3 = 002000+100010$, while $V^2 \supseteq (010101)^2$, and $010101$ cannot arise from $V^1$. \hal

\begin{lem}\label{ad6} If $\mu^0 = \l_2^0$, then $\l=\l_2$.
\end{lem}

\pf Suppose $\mu^0 = \l_2^0$, so that $\l = \l_2+x\l_n$. Assume $x>0$. Then $V^2 = \wedge^2\o_{s+1}$ and 
$V^2(Q_Y) \supseteq 10\cdots 01 + 010\cdots 010$, so that $V^2 \supseteq (\wedge^2\o_{s+1}\otimes \wedge^2\o_s) - 0$. 
For $s=3$ this contains a repeated summand of $S$-value at least 4 (namely $(110011)^2$), and so the parabolic argument shows that the same holds for any $s\ge 3$, a contradiction. \hal

\begin{lem}\label{ad7} We have $\mu^0 \ne (\l_2^0)^*$.
\end{lem}

\pf Suppose $\mu^0 = (\l_2^0)^*$. Here $V^1 = \wedge^2\o_s$ and $V^2(Q_Y) \supseteq 
0\cdots 0101 + 0\cdots 01000$ which is equal to $0\cdots 0100 \otimes 0\cdots 01$. Hence $V^2 \supseteq \o_s\otimes \wedge^3\o_s$. 
For $s=3$ this contains a repeated summand of $S$-value at least 4 (namely $(110110)^2$), and so the parabolic argument shows that the same holds for any $s\ge 3$, a contradiction. \hal

\begin{lem}\label{ad8} We have $\mu^0 \ne \l_1^0+\l_{r_0}^0$.
\end{lem}

\pf Suppose $\mu^0 = \l_1^0+\l_{r_0}^0$. Then $V^1 = \o_{s+1}\otimes \o_s-0$, while $V^2(Q_Y)$ contains
$10\cdots 0100 + 0\cdots 02 + 0\cdots 010$. Hence $V^2$ contains $(\o_{s+1}\otimes \wedge^3\o_s)+S^2\o_s$. 
For $s=3$ this contains a repeated summand of $S$-value at least 4 (namely $(110101)^2$), and so the parabolic argument shows that the same holds for any $s\ge 3$, a contradiction. \hal

\begin{lem}\label{ad9} Suppose $\mu^0$ or its dual is one of the following weights, with $s=3$ or $5$:
\[
\begin{array}{l}
s=3,5:\,3\l_1^0,\,\l_3^0 \\
s=3: \,4\l_1^0,\,5\l_1^0,\,\l_4^0,\,\l_5^0,\,\l_6^0,\,\l_1^0+\l_2^0.
\end{array}
\]
Then either $s=3$ and $\l \in \{4\l_1,\,\l_4\}$, or $s\in \{3, 5\}$ and $\l\in \{3\l_1,\,\l_3\}$.
\end{lem}

\pf Set $x = \langle \lambda,\beta_n\rangle$.
If $\mu^0$ is one of the weights listed in the statement of the lemma (rather than the dual of one of these), and $x=0$, then either $\l$ is as in the conclusion, or $\l = 5\l_1$, $\l_5$, $\l_6$ or $\l_1+\l_2$ (with $s=3$). In each of the latter cases a Magma computation  shows that $V_Y(\l)\downarrow X$ is not MF. Hence $x>0$ for such $\mu^0$.

In Table \ref{remna}, for each possible $\mu^0$ we list $V^1$, some summands of $V^2(Q_Y)$, and a repeated 
summand $Z$ of $V^2$; in the second column of the last two rows of the table, the module 
$D = \omega_4\otimes \wedge^2\omega_4 - \wedge^3\omega_4$. 
    In all cases other than $\mu^0 = (\l_3^0)^*$ (for $s=3$) and $(\l_5^0)^*$ we have $S(Z) \ge S(V^1)+2$, so this gives a contradiction.  In the exceptional cases $V^2 \supseteq (011101)^2$ and $(111111)^3$, respectively, and one checks that $V^1$ cannot contribute such  summands.  This completes the proof of the lemma.  \hal

\begin{table}[h!]
\caption{}\label{remna}
\[
\begin{array}{|cccc|}
\hline
\mu^0 & V^1 & V^2(Q_Y) \supseteq & V^2 \supseteq \\
\hline
3\l_1^0  & S^3\o_{s+1} & 30\cdots 0\otimes 0\cdots 010 & (111002)^2,\,s=3\\
   &&& (0110100101)^{13},\,s=5 \\
(3\l_1^0)^* & S^3\o_s & 0\cdots 0102 & (201101)^2,\,s=3 \\
   &&& (0111001010)^4, s=5 \\
\l_3^0 & \wedge^3\o_{s+1} &  0010\cdots 010+010\cdots 01 & (021011)^2,\,s=3 \\
   &&& (0111000011)^9, s=5 \\
(\l_3^0)^* & \wedge^3\o_s &  0\cdots 01001 & (011101)^2,\,s=3 \\
   &&&(0111000011)^9, s=5 \\
4\l_1^0 & S^4\o_4 & 40\cdots 0 \otimes 0\cdots 010 & (101301)^2 \\
(4\l_1^0)^* & S^4\o_3 & 0\cdots 0103 & (301002)^2 \\
5\l_1^0 & S^5\o_4 & 50\cdots 0 \otimes 0\cdots 010 & (101401)^2 \\
(5\l_1^0)^* & S^5\o_3 & 0\cdots 0104 & (203101)^2 \\
\l_4^0 & \wedge^4\o_4 & 00010\cdots 010 + 0010\cdots 01 & (111021)^2 \\ 
(\l_4^0)^* & \wedge^4\o_3 &  0\cdots 010001 & (101201)^2 \\
\l_5^0 & \wedge^5\o_4 & 000010\cdots 010 + 00010\cdots 01 & (112012)^2 \\
(\l_5^0)^* & \wedge^5\o_3 &  0\cdots 0100001 & (111111)^2 \\
\l_6^0 & \wedge^6\o_4 & 0000010\cdots 010 + 000010\cdots 01 & (211013)^2 \\
(\l_6^0)^* & \wedge^6\o_3 &  0\cdots 01000001 & (401012)^2 \\
\l_1^0+\l_2^0 & D & 110\cdots 0\otimes 0\cdots 010 - 10\cdots 0 & (200021)^2 \\
(\l_1^0+\l_2^0)^* & D^*  & 0\cdots 0101\otimes 0\cdots 01 - 0\cdots 01001 & (111110)^2 \\
\hline
\end{array}
\]
\end{table}


\section{The case $X = A_{2s+1}$ with $\d = b\o_{s+1}$ and $b\ge 2$}

In this section we  prove Theorem \ref{mainthm} in the case where  $X = A_{2s+1}$, $\d = b\om_{s+1}$ with $s \ge 1$,  $b \ge 2$, and $Y = SO(W)$ or $Sp(W)$ in accordance with Lemma \ref{setup}, where $W = V_X(\d)$.  As in Theorem \ref{oms+om(s+1)} we shall assume the Inductive Hypothesis in the proofs.

  We will prove the following result.

\begin{thm}\label{0...0b0...0} Let $X = A_{2s+1}$ with $s\ge 1$, and assume $X$ is embedded in $Y =C_n$ or $D_n$ via $\d = b\om_{s+1}$, where $b\ge 2$. 
Assume the Inductive Hypothesis and let $V  = V_Y(\l)$ with $\l \ne 0,\l_1$.  Then $V\downarrow X$ is MF if and only if one of the following holds:
\begin{itemize}
\item[{\rm (i)}]  $\l = \l_2$ or  $2\l_1$,
\item[{\rm (ii)}]  $X = A_3$, $b= 2$ and $\l = \l_3,\, 3\l_1,\,\l_9$ or $\l_{10}$.
\end{itemize}
\end{thm}

We begin by verifying that $V_Y(\l)\downarrow X$ is indeed MF for $\l$ as in the conclusion.

\begin{lem}\label{MFa2s1} Let $X<Y$ be as in the hypothesis of Theorem $\ref{0...0b0...0}$, and suppose $\l$ is as in conclusion (i) or (ii). Then  $V_Y(\l) \downarrow X$ is MF.
\end{lem}

\pf  For $\l = \l_2$ or $2\l_1$ this follows directly from Theorem 1 of \cite{MF}.  Now suppose $X = A_3$ with $\d = 2\om_2$.  If $\l = \l_3$ the assertion also follows from Theorem 1 of \cite{MF} while if $\l = 3\l_1$ it follows from a Magma check that $S^3(2\om_2) - 2\om_2$ is MF.

Finally, let $X=A_3$, $\d=2\o_2$ and let $\l$ be a spin module for $Y = D_{10}$, so that $\l = \l_9$ or $\l_{10}$.
 We claim that  
\begin{equation}\label{spinres}
V_Y(\l) \downarrow X = 311 + 113.
\end{equation}
  As $V^1 = 11$, it must be the case that $V_Y(\l)\downarrow X$ has 
 a composition factor of highest weight $11x$.   Moreover the dual of this composition
 factor must also occur as the embedding of $X$ in $Y = SO_{20}$ is invariant under the graph automorphism of $X$.  A dimension argument forces $x \le 3$. 
 
 Now $V^2 \supseteq 20 \otimes 11 \supseteq 31$ and this summand cannot arise from $V^1$. Therefore $V_Y(\l) \downarrow X \supseteq 31y$ for some $y$ and the usual argument with the central torus $T$ shows that $3+2+3y = 1+2 +3x-4$ so that $y = x-2$.  In particular $x \ge 2$.
 
  Let $z$ generate the center of $SL_4$  so that $z$ acts as $\zeta$, a fourth root of unity, on the natural module.  Then $z$ acts as $\zeta^4$ on $020$ and hence $X$ is embedded as an adjoint group in $SO_{20}$ and has center of order at most 2 in its action on the spin module.  But for $x = 2$,  $z$ acts via $\zeta^{1+2-2}$, a contradiction.  Therefore $x = 3$  and a dimension check  gives the assertion (\ref{spinres}).   \hal

We now embark on the proof of the `only if' part of Theorem \ref{0...0b0...0}. Adopt the hypotheses of the theorem, and assume that 
$V_Y(\l)\downarrow X$ is  MF.

We have  $L_X' = A_{2s}$ and we write $b = 2k$ or $2k+1$.   Then
$L_Y' = C^0 \times \cdots \times C^{k-1} \times C^k$  where $C^i$ has type $A$ for $0 \le i \le k-1$. For each $i$  with $0 \le i \le k-1$, $L_X'$  is embedded in $C^i$ via the representation $i\om_s + (b-i)\om_{s+1}$.  

If $b = 2k$,  then $C^k$ has type $B$ or $D$ and $L_X'$ is embedded in $C^k$ via the representation  $(\frac{b}{2}\om_s +\frac{b}{2}\om_{s+1})$.

If $b = 2k+1$, then $Y$ could be symplectic or orthogonal, $C^k$ has type $A$ and $L_X'$ is embedded via the representation $k\om_s + (k+1)\om_{s+1}$.

Throughout the proofs to follow we assume the Inductive Hypothesis. There are minor differences between the cases $s >1$ and $s = 1$. We will often produce weights of the form $\nu =c\om_{s-1}+a\om_s +e\om_{s+1} +d\om_{s+2}$, but when  $s = 1$ the terms $\om_{s-1}$ and $\om_{s+2}$ do not occur.  Rather than repeatedly calling attention to these differences we will simply ignore them and leave the small changes to the reader.  

In several places to follow we will consider Magma computations involving small wedge or symmetric powers of representations of $L_X'$ with highest weights of the form $x\om_s +y\om_{s+1}$.  When $s = 1$ we  find composition factors of the form $\xi = e\om_1 +f\om_2 = e\om_s + f\om_{s+1}$ with certain multiplicities.  These are obtained by subtracting linear combinations of roots of form $x\a_s + y\a_{s+1}$ from a suitable dominant weight above $\xi$.  It then follows that for  $s > 1$ there are composition factors of highest weight  $\nu$ as above with multiplicity at least as large as that as the $s = 1$ case. Consequently we can use Magma computations for  a group of type  $A_4$ to obtain information, again with the convention of the preceding paragraph.

\begin{lem}\label{mui} {\rm (i) }  $\mu^i \ne 0$ for at most one value of $i  >0$.

{\rm  (ii)}  $\mu^k = 0$ or $\l_1^k$ if either  $b = 2k > 4$ or $b = 4$ and $s > 1$.

{\rm  (iii)} If $2 \le i \le k-1$ then $\mu^i = 0,\,  \l_1^i$, or $\l_{r_i}^i$.
\end{lem}

\pf   Part (ii) follows from the Inductive Hypothesis, and (iii) follows from  \cite{MF}.  Finally,  (i) follows from (ii) and
Lemma 7.3.1 of \cite{MF} and Lemma \ref{stemb}, since $L(i\om_s + (b-i)\om_{s+1}) \ge 2$ for each $0 < i \le k-1$ and similarly for $i = k$.  \hal

At this point we divide the proof into subcases, according to whether $b\ge 4$, $b=3$ or $b=2$.

\subsection{The case where $b \ge 4$}

We continue the proof of Theorem \ref{0...0b0...0}, assuming now that $b\ge 4$.

\begin{lem}\label{muk=0}  Assume $b\ge 4$ is even.  Then $\mu^k = 0$.
\end{lem}

\pf Assume $\mu^k\ne 0$.  The embedding of $L_X'$ in $C^k$ is given by  the weight $\frac{b}{2}\om_s + \frac{b}{2}\om_{s+1}$ so the inductive hypothesis implies that either $\mu^k = \l_1^k$ or $b = 4, s = 1$ and
$\mu^k = \l_2^k$.  

Assume $\mu^k = \l_1^k$.  Then Lemma \ref{mui} shows that $\mu^{i} = 0$ for $1\le i\le k-1$.
Using the composition factor of $V^2(Q_Y)$ afforded by $\l -\gamma_k-\b_1^k$ we see
that $V^2 \supseteq (\mu^0 \downarrow L_X') \otimes (\frac{b+2}{2}\om_s +\frac{b-2}{2}\om_{s+1})
\otimes \wedge^2(\frac{b}{2}\om_s + \frac{b}{2}\om_{s+1})$. The last tensor factor contains a composition factor of highest weight $\om_{s-1} +(b-2)\om_s +(b+1)\om_{s+1}$ and it follows using \cite[7.1.5]{MF} that $V^2$ contains a repeated composition factor of $S$-value at
least $S(\mu^0 \downarrow L_X' ) + 3b-3$.  As $S(V^1) = S(\mu^0 \downarrow L_X' ) +b$
this is a contradiction.  

Now suppose $\mu^k = \l_2^k$.  Then $b = 4$ and $X = A_3$. The embeddings of $L_X'$ in $C^0, C^1, $ and $C^2$ are given by $(04),(13)$, and $(22)$, respectively.  As above $\mu^1 = 0$ and using the representation afforded by $\l -\gamma_2-\b_1^2-\b_2^2$ we see that $V^2 \supseteq  (\mu^0 \downarrow L_X') \otimes (31) \otimes \wedge^3(22) \supseteq (\mu^0 \downarrow L_X') \otimes (31) \otimes (25)^4 $ which has a repeated summand of $S$-value at least
$S(\mu^0 \downarrow L_X' ) + 11$.  As  $S(V^1)  = S(\mu^0 \downarrow L_X' ) + 7$ this is a contradiction.  \hal

\begin{lem}\label{mu(k-1)=0even}  Assume $b \ge 4$ is even.  Then $\mu^{k-1} = 0$ and $\la \l, \g_k \ra = 0$.
\end{lem}

\pf  Suppose $\mu^{k-1} \ne 0$  and apply \cite{MF}.  If $b \ge 6$ then  $\mu^{k-1}  = \l_1^{k-1}$ or 
 $\l_{r_{k-1}}^{k-1}$.  While if $b = 4$ then $\mu^{k-1} = \mu^1$, and either it or its dual is $\l_1^1, \l_2^1,$ or $2\l_1^1.$ By Lemma \ref{muk=0} we know that $\mu^k = 0$,  and by Lemma \ref{mui} we have $\mu^i = 0$ if $i \ne 0,\, k-1$. 
Also the possibilities for $\mu^0$ are given by \cite{MF}, and it follows using \cite[Lemma 7.3.1]{MF} and Lemma \ref{stemb} that $(\mu^0\otimes \mu^{k-1})\downarrow L_X'$ is not MF unless  $\mu^0 = 0, \l_1^0$ or $\l_{r_0}^0$.   So $\mu^0$ must be one of these.
 
 First assume 
 $\mu^{k-1} =\l_{r_{k-1}}^{k-1}$.  Then $V^2(Q_Y) \supseteq \mu^0 \otimes \l_{r_{k-1}-1}^{k-1} \otimes \l_1^k$.  Restricting to $L_X'$ we see that $V^2$ contains $(\mu^0 \downarrow L_X')
 \otimes (\om_{s-1}+b\om_s+(b-1)\om_{s+1}) \otimes (\frac{b}{2}\om_s +\frac{b}{2}\om_{s+1})$ and
 by \cite[7.1.5]{MF} this contains $(\mu^0 \downarrow L_X')
 \otimes (2\om_{s-1}+(\frac{3b}{2}-1)\om_s+(\frac{3b}{2}-2)\om_{s+1}+\om_{s+2})^2$, giving a contradiction as 
$S(V^1) = S(\mu^0)+b$.
 
 
We now complete the treatment of the case $b\ge 6$, where the only remaining case is $\mu^{k-1} =\l_1^{k-1}$. Then $V^2 \supseteq (\mu^0 \downarrow L_X') \otimes  (\l_{r_{k-2}}^{k-2}\downarrow L_X') \otimes \wedge^2((\frac{b-2}{2})\om_s + (\frac{b+2}{2})\om_{s+1})$.
If $s \ge 2$ this contains  $(\mu^0 \downarrow L_X') \otimes  ((\frac{b+4}{2})\om_s + (\frac{b-4}{2})\om_{s+1}) )\otimes ((b-1)\om_s + b\om_{s+1} +\om_{s+2})$  and this contains  $(\mu^0 \downarrow L_X') \otimes  (\om_{s-1}+(\frac{3b}{2})\om_s + (\frac{3b}{2}-3)\om_{s+1} + 2\om_{s+2})^2$  of $S$-value at least $3b +S(\mu^0 )$, again a contradiction.  If $s  = 1$ the same argument yields a repeated summand with $S$-value $3b-3 +S(\mu^0)$, again
a contradiction.

 Now assume $b = 4$  so that $k = 2$ and the embeddings in $C^0$ and $C^1$ are given by $4\om_{s+1}$ and $\om_s+3\om_{s+1}$, respectively.  There are several remaining possibilities for $\mu^1$. Recall that $\mu^0 = 0,\,\l_1^0$ or $\l_{r_0}^0$. 
 
 First assume $\mu^1 = \l_1^1$.  Then $S(V^1) = S(\mu^0) + 4$.  On the other hand
 $V^2(Q_Y) \supseteq ((\mu^0 + \l_{r_0}^0) \downarrow L_X') \otimes \wedge^2(\om_s + 3\om_{s+1})$. The second tensor factor contains $(\om_{s-1}+\om_s + 5\om_{s+1} +\om_{s+2})$ and $ (3\om_s + 4\om_{s+1}+\om_{s+2})$.
  
   If $\mu^0 = 0$, then $S(V^1) = 4$ while $V^2 \supseteq 4\om_s \otimes \wedge^2(\om_s + 3\om_{s+1})$.  By the above this contains  $(\om_{s-1}+5\om_s + 5\om_{s+1}+\om_{s+2})^2$ and we obtain a contradiction using  $S$-values. 
   
  If $\mu^0 = \l_1^0$, then   $S(V^1) = 8$ and $\mu^0 + \l_{r_0}^0 =  \l_1^0 + \l_{r_0}^0 $, so that  $V^2 \supseteq (4\om_s + 4\om_{s+1}) \otimes \wedge^2(\om_s + 3\om_{s+1}) \supseteq (4\om_s + 4\om_{s+1}) \otimes  (3\om_s + 4\om_{s+1}+\om_{s+2})$ and $V^2$ contains $(\om_{s-1}+6\om_s + 7\om_{s+1}+2\om_{s+2})^2$.  Once again $S$-values yield a contradiction.
  
  Now assume $\mu^0 = \l_{r_0}^0$.  Here $\mu^0 + \l_{r_0}^0 = 2 \l_{r_0}^0$ and   $S(V^1) = 4+4$ while $V^2 \supseteq (2\om_{s-1}+4\om_s +2\om_{s+1})  \otimes \wedge^2(\om_s + 3\om_{s+1}) \supseteq (2\om_{s-1}+4\om_s +2\om_{s+1})  \otimes  (3\om_s + 4\om_{s+1}+\om_{s+2})$.  Hence  $V^2$ contains $(3\om_{s-1}+6\om_s + 5\om_{s+1}+2\om_{s+2})^2$.  Once again $S$-values yield a contradiction.  This completes the analysis of the case $\mu^1 = \l_1^1$.
  
  Next suppose that $\mu^1 = 2\l_{r_1}^1$ (still with $b=4$).  Then $S(V^1) = S(\mu^0) + 8$ and $V^2 \supseteq 
 (\mu^0 \downarrow L_X') \otimes ((\l_{r_1-1}^1+\l_{r_1}^1)\downarrow L_X') \otimes (2\o_s + 2\o_{s+1})$ and this contains $(\mu^0 \downarrow L_X') \otimes (10\om_s+\om_{s+1} +\om_{s+2}) \otimes  (2\om_s + 2\om_{s+1})$.  So $V^2 \supseteq  (\mu^0 \downarrow L_X') \otimes (\om_{s-1}+11\om_s +2\om_{s+1} +2\om_{s+2})^2$
 which has $S$-value at least $S(\mu^0) + 13$, a contradiction.
 
 Suppose $\mu^1 = 2\l_1^1$.  Then $S(V^1) = 8,12$, or $12$ according to  $\mu^0 = 0, \l_1^0$, or $\l_{r_0}^0$.  Also $V^2(Q_Y) \supseteq (\mu^0 + \l_{r_0}^0) \otimes (\l_1^1+\l_2^1)$ and one checks that the restriction to $L_X'$ of the second tensor factor contains $(\om_{s-1}+2\om_s + 8\om_{s+1}+\om_{s+2})^2$.  Therefore $V^2$ contains a repeated composition factor of $S$-value at
 least $4+10, 8+10$, or $8+10$, according to whether $\mu^0 = 0, \l_1^0$, or $\l_{r_0}^0$.
 An $S$-value consideration gives a contradiction.
 
 Next assume $\mu^1 = \l_{r_1 -1}^1$.  Then $S(V^1) = S(\mu^0) + 7$ (or $S(\mu^0) + 8$ if $s>1$) and
 $V^2(Q_Y) \supseteq \mu^0 \otimes \l_{r_1 -2}^1 \otimes \l_1^2$.  Therefore $V^2 \supset
( \mu^0 \downarrow L_X') \otimes (4\om_{s-1} +5\om_s+\om_{s+1}) \otimes (2\om_s +2\om_{s+1}) \supseteq ( \mu^0 \downarrow L_X') \otimes (6\om_{s-1} +5\om_s +\om_{s+1}+2\om_{s+2})^2$ which has $S$-value at least $S(\mu^0) + 11$, a contradiction.
 
 The final case is $\mu^1 = \l_2^1$.  Then $S(\mu^1 \downarrow L_X') = 7$ (respectively 8 if $s > 1$) so that  $S(V^1) = 7,11$, or $11$ (respectively $ 8,12,12$) according as  $\mu^0 = 0, \l_1^0$, or $\l_{r_0}^0$.   Then $V^2(Q_Y) \supseteq (\mu^0 +\l_{r_0}^0) \otimes \l_3^1$. The   restriction to $L_X'$ of the second tensor factor contains 
$(\om_{s-2} + 2\om_s + 7\om_{s+1} +2\om_{s+2})^2$.
  It follows that $V^2$ contains a repeated composition factor with $S$-value at least $4+9, 8+9$, or $8+9$ according as $\mu^0 =  0, \l_1^0$, or $\l_{r_0}$ and this is a contradiction. 
 
 It still remains to show that $\la \l, \g_k \ra = 0$.  Suppose otherwise.  Then $V^2(Q_Y)$ contains
$(\mu^0 \otimes \cdots \otimes \mu^{k-2}) \otimes \l_{r_{k-1}}^{k-1} \otimes \l_1^k$.  By \cite[7.1.5]{MF}, the
restriction to $L_X'$ of $\l_{r_{k-1}}^{k-1} \otimes \l_1^k$ contains a repeated composition factor
of $S$-value at least $2b-2 \ge 6$.  As $S(V^1) = S(\mu^0 \otimes \cdots \otimes\mu^{k-2})$ this is a contradiction.
 \hal
 
 \begin{lem}\label{muk=0,odd}  Assume $b \ge 5$ is odd.  Then $\mu^k = 0$ and $\la \l, \g_{k+1} \ra = 0$. 
 
 \end{lem}
 
 \pf  Assume $\mu^k \ne 0$.  
By Lemma \ref{setup}, $Y$ is symplectic if $s$ is even and orthogonal if $s$ is odd. Here $C^k$ is of type $A$ and the embedding of $L_X'$ in $C^k$ is given by $\frac{b-1}{2}\om_s + \frac{b+1}{2}\om_{s+1}$. In particular
the natural module for $Y$ has even dimension in either case.
Also  Lemma 7.3.1 of \cite{MF} and Lemma \ref{stemb} imply that $\mu^i = 0$ for $1 \le i \le k-1$ and $\mu^0$ is $0,\l_1^0$ or $\l_{r_0}^0$.
 
 The main result of \cite{MF} implies that $\mu^k = \l_1^k$ or $\l_{r_k}^k$.
 In the former case $V^2$   contains $(\mu^0\downarrow L_X') \otimes (\frac{b+3}{2}\om_s + \frac{b-3}{2}\om_{s+1}) \otimes \wedge^2(\frac{b-1}{2}\om_s + \frac{b+1}{2}\om_{s+1})$.  The latter tensor factor contains $b\om_s +(b-1)\om_{s+1} +\om_{s+2}$ so $V^2$ has a repeated tensor factor of highest weight  $(\mu^0\downarrow L_X') + (\om_{s-1}+\frac{3b+1}{2}\om_s + \frac{3b-7}{2}\om_{s+1} +2\om_{s+2})^2$.  Therefore $V^2$ has a repeated composition factor of $S$-value at least $S(\mu^0) + 3b-3$ whereas $S(V^1) = S(\mu^0) + b$.  This is a contradiction.
 
 Now assume $\mu^k = \l_{r_k}^k$.  Then $V^2 \supseteq (\mu^0\downarrow L_X') \otimes \wedge^3(\frac{b+1}{2}\om_s + \frac{b-1}{2}\om_{s+1})$ or $(\mu^0\downarrow L_X') \otimes ((\l_{r_k-1}^k +\l_{r_k}^k)\downarrow L_X')$, according to whether $Y$ is $SO$ or $Sp$. Indeed in the former case the irreducible summand is afforded by $\l - \b_{k-1}^k-\b_k^k-\g_{k+1}$ while in the latter case it is afforded by $\l - \b_k^k-\g_{k+1}$. Using the fact that  $\l_{r_k-1}^k +\l_{r_k}^k =  (\l_{r_k-1}^k \otimes \l_{r_k}^k) - \l_{r_k-2}^k  $ we count weights and find that in either case $V^2 \supseteq (\mu^0\downarrow L_X') \otimes  (2\om_{s-1}+\frac{3b-3}{2}\om_s + \frac{3b-3}{2}\om_{s+1}+\om_{s+2})^2$ so that $V^2$ has a repeated composition factor of $S$-value at least $S(\mu^0) + 3b-3$, a contradiction.  Therefore $\mu^k = 0$.  

Finally assume that $\la \l, \g_{k+1} \ra > 0$.   Then $V^2(Q_Y)$ has a composition factor $(\mu^0 \otimes \cdots \otimes \mu^{k-1}) \otimes (\l_{r_k-1}^k)$ or $(\mu^0 \otimes \cdots \otimes \mu^{k-1}) \otimes 2\l_{r_k}^k$ depending on whether $Y = SO$ or $Sp$.  Here we claim that the last tensor factor contains  a repeated composition factor of highest weight $\om_{s-1} + (b-1)\om_s + (b-3)\om_{s+1} +\om_{s+2}$  or $2\om_{s-1} + (b-1)\om_s + (b-3)\om_{s+1} +2\om_{s+2}$, respectively.  For example consider  the wedge square $\l^k_{r_k-1}$.  Working in $A_4$ we use Magma to see that for $b = 5$ and $7$ the wedge square contains a composition factor of multiplicity 2 and highest weight $(0,b+1,b-1,0) - (\a_1 + 3\a_2 + 3\a_3 + \a_4) = (1,b-1,b-3,1).$ 
It then follows that such a composition factor also occurs for $b>7$ by observing
 that weight spaces for higher dominant weights are the same as for $b = 7$ (using \cite{cavallin} and a parabolic argument in the usual way).
Similar considerations apply for the symmetric square $2\l_{r_k}^k$.
As $S(V^1) = S(\mu^0 \otimes \cdots \otimes \mu^{k-1})$ we obtain an $S$-value contradiction. \hal

\begin{lem}\label{doneexceptb=2,3} The conclusion of Theorem $\ref{0...0b0...0}$ holds
if $b \ge 4$.
\end{lem}

\pf Assume $b\ge 4$.  Then the above lemmas show that $\mu^k = 0$, and that $\la \l, \g_k \ra = 0$ if $b$ is even and 
$\la \l, \g_{k+1} \ra = 0$ if $b$ is odd. We shall adopt arguments in Section 9 of \cite{MF}, which  show that  $\mu^i = 0$ and $\la \l, \g_i \ra = 0$ for various values of $i$.  However the local arguments used  to achieve these equalities assumed that
the relevant factors of $L_Y'$ are of type $A$, so we must adjust the considerations to include factors of type $B,C,$ or $D$.

When $b$ is even, Lemma \ref{mu(k-1)=0even} shows that $\mu^{k-1}  = 0$, while if $b$ is odd then $C^k$ has type $A$. It follows that in any case we can use the  arguments in Section 9 of \cite{MF} to obtain $\mu^1 = \cdots = \mu^{k-1} = 0$  and $\la \l, \g_1 \ra  = \cdots =  \la \l, \g_k \ra = 0$ in all cases. 

 It remains to consider $ \mu^0$, as this will determine $\l$.  By
 \cite{MF}, $\mu^0$ or its dual is one of $\l_1^0$, $\l_2^0,$ $\l_3^0$ (for $s=1$ and $b=4,5,6$),  $\l_4^0$ (for $s=1$ and $b=4$), $2\l_1^0$, $3\l_1^0$ ($s = 1$ and  $b= 4, 5$), or $\l_1^0 +\l_{r_0}^0$.
The dual versions, and also $\l_1^0 +\l_{r_0}^0$, are ruled out as in Section 9 of \cite{MF}.  Also $\l$ is as in the conclusion of the theorem  if $\mu^0 = \l_1, \l_2$ or $2\l_1$. 
 Therefore we must rule out the cases $\l = \l_3,  \l_4$ or $3\l_1$ where in each case $X = A_3$.  Using
 Magma for $b = 4,5,6$ we find that $\wedge^3(\d) \supseteq (3\d - (2\a_1+5\a_2+2\a_3))^2 = (1,3b-6,1)^2$, 
$ \wedge^4(\d) \supseteq (4\d - (2\a_1+5\a_2+2\a_3))^2 = (1,4b-6,1)^2$, and $S^3(\d) \supseteq (3\d - (\a_1+4\a_2+\a_3))^2 = (1,3b-6,1)^2$.  It  follows using \cite{cavallin} that the same holds for all $b \ge 4$. Then from Lemma \ref{symwed},  we see that $V\downarrow X$ is not MF.  \hal
 
 \subsection{The case where $b=3$}

 In this subsection we prove Theorem \ref{0...0b0...0} in the case where $b = 3$. Here  $L_Y' = C^0 \times C^1$ and $L_X'$ is embedded in the factors via the highest weights $3\om_{s+1}$ and $\om_s + 2\om_{s+1}$, respectively.
 
Now \cite{MF}  implies that $\mu^1$ or its dual is in $\{0, \l_1^1, \l_2^1, 2\l_1^1\}$ whereas
 $\mu^0$ or its dual is in $\{0, \l_i^0, c\l_1^0 (c\le 4), \l_1^0 +\l_{r_0}^0, \l_1^0+\l_2^0\}$, where in the $\mu^0$ case some of the weights only occur if $s = 1$.
 
 \begin{lem}\label{b=3mu1=0}  $\mu^1 =0$.
  \end{lem}
  
  \pf Suppose $\mu^1 \ne 0$.
  Then by Lemma \ref{stemb} together with \cite[Lemma 7.3.1]{MF}, we have  $\mu^0 = 0$,  $\l_1^0$, or $\l_{r_0}^0$ and $S(\mu^0) = 0,3$ or 3, respectively.
  
 Assume $\mu^1 = \l_1^1$.  Then $S(V^1) = S(\mu^0) +3$ and 
  $V^2 \supseteq ((\mu^0 + \l_{r_0}^0)\downarrow L_X') \otimes \wedge^2(\om_s + 2\om_{s+1})$.
 The first factor has a composition factor $3\om_s$, $3\om_s+3\om_{s+1}$, or $6\om_s$ according to whether $\mu^0 = 0, \l_1^0$ or $ \l_{r_0}^0$.  And the second tensor factor
 has summands $(\om_{s-1} +\om_s + 3\om_{s+1} + \om_{s+2})$ and  $(3\om_s + 2\om_{s+1} + \om_{s+2})$.  It follows that $V^2$ has a repeated composition factor of highest weight
 $(\om_{s-1} +4\om_s + 3\om_{s+1} + \om_{s+2})$, $(2\om_{s-1} +3\om_s + 5\om_{s+1} + 2\om_{s+2})$ or  $(\om_{s-1} +7\om_s + 3\om_{s+1} + \om_{s+2})$ according as
 $\mu^0 = 0,\l_1^0$ or $\l_{r_0}^0$.  A consideration of $S$-values gives a contradiction.
 
 Next suppose that $\mu^1 = 2\l_1^1$.  Here  $S(V^1) = S(\mu^0) +6$ and 
  $V^2 \supseteq ((\mu^0 + \l_{r_0}^0)\downarrow L_X') \otimes ((\l_1^1+\l_2^1)\downarrow L_X')$.
 A Magma check shows that the second factor contains $(\om_{s-1} + 2\om_s +5\om_{s+1} +\om_{s+2})^2$.  Therefore $V^2$ has a repeated composition factor of $S$-value at least $3+7$,
 $6+7$, or $6+7$ according as $\mu^0 = 0$,  $\l_1^0$, or $\l_{r_0}^0$.  This is a contradiction.
 
 If $\mu^1 = \l_2^1$, then $S(V^1) \le S(\mu^0) + 6$.  Now
  $V^2 \supseteq ((\mu^0 + \l_{r_0}^0)\downarrow L_X') \otimes \wedge^3(\om_s + 2\om_{s+1})$.
  The second factor contains $(\om_{s-1} + 3\om_s + 3\om_{s+1} + 2\om_{s+2})^2$ so $V^2$ has a repeated composition factor of $S$ value at least $3+6$, $6+6$, $6+6$, respectively, and we get the usual contradiction.

  If $\mu^1 = \l_{r_1}^1$ then $S(V^1) = S(\mu^0)+3$.  Now $V^2(Q_Y) \supseteq  \mu^0 \otimes \l_{r_1-2}^1$  or $ \mu^0 \otimes (\l_{r_1-1}^1 +\l_{r_1}^1)$ according to whether $Y$ is $SO$ or $Sp$.  Consequently as in the previous two paragraphs, $V^2$ has a repeated composition factor of $S$-value at least $S(\mu^0) + 6$ or $S(\mu^0) + 7$, according to whether  $Y$ is $SO$ or $Sp$. This is a contradiction.
  
  Next suppose $\mu^1 = \l_{r_1-1}^1$.  Then $S(V^1) \le S(\mu^0)+6$.  Here $V^2(Q_Y) \supseteq  \mu^0 \otimes (\l_{r_1-2}^1 + \l_{r_1}^1)$ no matter if  $Y$ has type $Sp$ or $SO$.   The second tensor factor can be expressed as $(\l_{r_1-2}^1 \otimes \l_{r_1}^1) -  \l_{r_1-3}^1$.  Therefore $V^2 \supseteq (\mu^0 \downarrow L_X') \otimes ((\wedge^3(2\om_s + \om_{s+1}) \otimes (2\om_s + \om_{s+1})) - \wedge^4(2\om_s + \om_{s+1}))$ and Magma shows that $V^2 \supseteq (\mu^0 \downarrow L_X')  \otimes ( 8\om_s + \om_{s+1} + \om_{s+2})^2$.  An $S$-value argument gives a contradiction.  
  
  Finally, assume that $\mu^1 = 2\l_{r_1}^1$.   Then $S(V^1) = S(\mu^0) + 6$.  If $Y = SO$ then $V^2(Q_Y) \supseteq \mu^0 \otimes (\l_{r_1-2}^1 + \l_{r_1}^1)$ and we get a contradiction as in the previous paragraph.  Suppose $Y = Sp$.  In this case
 $V^2(Q_Y) \supseteq \mu^0 \otimes (\l_{r_1-1}^1 + 2\l_{r_1}^1)$. The second tensor factor can be expressed as $(\l_{r_1-1}^1 \otimes  2\l_{r_1}^1) - (\l_{r_1-2}^1 + \l_{r_1}^1)$ and Magma shows that $V^2 \supseteq (\mu^0 \downarrow L_X') \otimes (\om_{s-1} + 7\om_s + 3\om_{s+1} + \om_{s+2})^2$ and we get the usual contradiction.  \hal
 
 \begin{lem}\label{mu0=m1=0} We have  $\mu^0 \ne 0$. 
  \end{lem} 
  
 \pf  Suppose $\mu^0 = 0$. By Lemma \ref{b=3mu1=0} we  have $\mu^1 = 0$  so that $V^1 = 0$. We must consider 
$\la \l,\g_1\ra$ and $\la \l,\g_2\ra$.  If $\la \l, \g_1 \ra \ne 0$, then $V^2(Q_Y) \supseteq \l_{r_0}^0 \otimes \l_1^1$ and $V^2 \supseteq 3\om_s \otimes (\om_s + 2\om_{s+1}) \supseteq  (4\om_s + 2\om_{s+1}) \oplus  (\om_{s-1} + 2\om_s + 3\om_{s+1})$.  And if  $\la \l, \g_2 \ra \ne 0$, then $V^2(Q_Y)$ contains $ 2\l_{r_1}^1$ or $\l_{r_1-1}^1$, according to whether $Y$ is $Sp$ or $SO$ and hence $V^2$ also contains $S^2(2\om_s + \om_{s+1})$ or  $\wedge^2(2\om_s + \om_{s+1})$, accordingly. Consequently $V^2$
 contains $(4\om_s + 2\om_{s+1})^2$, respectively $(\om_{s-1} + 2\om_s + 3\om_{s+1})^2$, which gives a contradiction if  $\la \l, \g_1 \ra \ne 0 \ne  \la \l, \g_2 \ra $.
 
 Suppose  $\la \l, \g_1 \ra  > 0$ and  $\la \l, \g_2 \ra = 0$.  Then $V^2= 3\om_s \otimes (\om_s + 2\om_{s+1})$.  Now $V^3(Q_Y) $ contains a summand afforded by $\l- \b_{r_0}^0 -2\g_1 - \b_1^1$  and so $V^3$ contains $\wedge^2(3\om_s) \otimes \wedge^2( \om_s + 2\om_{s+1})$ and a Magma check shows that this contains $(2\om_{s-1} + 6\om_s +2\om_{s+1}+2\om_{s+2})^2$.
 As $S(V^2) = 6$ an $S$-value argument gives a contradiction.  Therefore $\la \l, \g_1 \ra  = 0$. 
 
 Finally assume $\langle \lambda,\gamma_2 \rangle >0$.  First suppose $Y=SO$. If $\langle\lambda,\gamma_2\rangle = 1$,  then $V^2 = \wedge^2(2\omega_s+\omega_{s+1})$ and $V^3\supseteq \wedge^4(2\omega_s+\omega_{s+1})$, afforded by $\lambda-\beta_{r_1-2}-2\beta_{r_1-1}-\beta_{r_1}-2\gamma_2$.
Using Magma for the cases $s=1$ and $s=2$, we find that $V^2$ has $S$-value 5 if $s=1$ and $S$-value at most $6$ if $s\geq 2$. On the other hand, $V^3\supseteq (43)^2$, when $s=1$ and $(1171)^2$ for $s=2$. The latter weight is afforded by $8\omega_s+4\omega_{s+1}-2\alpha_{s-1}-5\alpha_s-\alpha_{s+1}$. The usual parabolic argument and a comparison of $S$-values provides the required contradiction.

Now if $\langle\lambda,\gamma_2\rangle>1$, we have $V^3(Q_Y)\supseteq 2 \lambda_{r_1-1}^1$ and so $V^2\supseteq
S^2(\wedge^2(2\omega_s+\omega_{s+1})) - \wedge^4(2\omega_s+\omega_{s+1})$. This contains $(43)^6$ if $s=1$ and $(1171)^4$ if $s=2$. The above argument now gives the result.

  Now suppose $Y = Sp$, which forces $s$ even and so $s\ge 2$.  Here we check in the usual way using   Magma together with \cite{cavallin} that $V^2 (Q_Y) = 2\l_{r_1}^1$ so that $V^2 = S^2(2\om_s+\om_{s+1})$ of $S$-value 6.
 
 First assume that $\la \l, \g_2 \ra  \ge 2$.  Then $V^3(Q_Y) \supseteq 4\l_{r_1}^1$ so that $V^3 \supseteq 
  S^4(2\om_s+\om_{s+1})$.  Using Magma for $A_4$ and the usual parabolic argument, we see that $V^3$ has a  repeated composition factor of highest weight $4(2\om_s+\om_{s+1}) - (2\a_s+2\a_{s+1}) = 2\om_{s-1}+6\om_s+2\om_{s+1}+2\om_{s+2}$. This gives a contradiction from an $S$-value comparison.
 
Finally assume $\la \l, \g_2 \ra  =1$.  Then $V^3(Q_Y)$ has a composition factor $0\cdots 020$ afforded by 
$\l - 2\b_{r_1}-2\g_2$ and therefore 
 $V^3 \supseteq S^2(\wedge^2(2\om_s+\om_{s+1})) - \wedge^4(2\om_s+\om_{s+1})$. A Magma check shows that  $V^3 \supseteq (2\om_{s-1} + 6\om_s + 2\om_{s+1} + 2\om_{s+2})^2$.  Once again an $S$-value comparison
gives a contradiction. This completes the proof of the lemma. \hal  
  
 \begin{lem}\label{l1l2l3} If $\mu^0 = \l_1^0,\l_2^0,$ or $\l_3^0$, then $\l=\l_2$. 
 \end{lem}
 
 \pf Write $\mu^0 = \l_i^0$ with $i = 1,2,3$ so that $S(V^1) = 3$, $6$ ($5$ if $s = 1$), or $9$
  ($6$ if $s = 1$), respectively. 
 We first claim that $\la \l, \g_1 \ra = 0 = \la \l, \g_2 \ra$.  
 
 Suppose
$\la \l, \g_1 \ra \ne 0$. Then $V^2(Q_Y)  \supseteq ((\l_i^0 + \l_{r_0}^0) \oplus \l_{i-1}^0) \otimes \l_1^1 = \l_i^0 \otimes  \l_{r_0}^0 \otimes \l_1^1$, with the understanding
that  $ \l_{i-1}^0 = 0$ if $i = 1$.  Therefore $V^2 \supseteq \wedge^i(3\om_{s+1}) \otimes 3\om_s \otimes (\om_s + 2\om_{s+1})$.  It follows that $V^2$ has a repeated composition factor of highest weight $(\om_{s-1} + 3\om_s + 4\om_{s+1} + \om_{s+2})$, $(\om_{s-1} + 4\om_s + 5\om_{s+1} + 2\om_{s+2})$, $(\om_{s-1} + 6\om_s + 4\om_{s+1} + 4\om_{s+2})$, 
according as $i = 1, 2, 3$.  These have $S$-values at least $7, 9, 10$ (15 in the last case when $i \ge 2$), respectively.   This is a contradiction and so $\la \l, \g_1 \ra = 0$.

Next suppose that $\la \l, \g_2 \ra \ne 0$.  Then $V^2(Q_Y)$ has a composition factor afforded by $\l - \g_2$ and
hence $V^2$  contains $V^1 \otimes \wedge^2(2\om_s + \om_{s+1})$ or  $V^1 \otimes S^2(2\om_s + \om_{s+1})$ according to whether $Y = SO$ or $Y = Sp$.  The second tensor factor
contains summands  $\om_{s-1} + 3\om_s +\om_{s+1} +\om_{s+2}$ and $ \om_{s-1} + 2\om_s +3\om_{s+1}$  if $Y = SO$ and summands  $\om_{s-1} + 3\om_s +\om_{s+1} +\om_{s+2}$ and
$ 4\om_s +2\om_{s+1}$ if $Y = Sp$. 

If $\mu^0 = \l_1^0$ then $V^1 = 3\om_{s+1}$ and  from the above paragraph $V^2$ contains
$(\om_{s-1} + 3\om_s +4\om_{s+1} +\om_{s+2})^2$  for both possibilities of $Y$.   If $i = 2$ or $3$, 
then $V^1$ contains $\om_s +4\om_{s+1} +\om_{s+2}$ or $3\om_s +3\om_{s+1} +3\om_{s+2}$,
respectively, and these have maximum $S$-value in $V^1$.  Tensoring with the composition factor
$\om_{s-1} + 3\om_s +\om_{s+1} +\om_{s+2}$ above, and using \cite[7.1.7]{MF}, we see that $V^2$ has a repeated composition factor of $S$-value at least $S(V^1) + 2$, a contradiction.  This establishes the claim.

We now have $\l = \l_1, \l_2,$ or $\l_3$.  By assumption $\l\ne \l_1$, and $\l_2$ is in the conclusion, so
suppose $\l = \l_3$.  Magma shows that for $s = 1$, $\wedge^3(030) \supseteq 131^2$ and $131 =090 - (2\a_1 +5\a_2 +2\a_3)$.  It follows that for $s > 1$, $\wedge^3(\l_1)$ contains
$(2\om_{s-1}+\om_s +3\om_{s+1} + \om_{s+2}+ 2\om_{s+3})^2$ and hence (using also Lemma \ref{symwed}) $V_Y(\l) \downarrow X$ is not MF.  \hal

\begin{lem}\label{al_r0}  $\mu^0 \ne a\l_{r_0}^0$ with $1\le a \le 4$.
\end{lem}

\pf  Suppose $\mu^0= a\l_{r_0}^0$ with $1\le a \le 4$.  Then  $S(V^1) = 3a$ and $V^2(Q_Y) \supseteq (\l_{r_0 -1}^0 +(a-1)\l_{r_0}^0) \otimes \l_1^1$. One then checks that $V^2 \supseteq (\om_{s-1} + (3a+1)\om_s + \om_{s+1}) \otimes (\om_s + 2\om_{s+1})$
which by \cite[7.1.5]{MF} contains  $(2\om_{s-1} + (3a+1)\om_s + 2\om_{s+1}+ \om_{s+2})^2$.  As this has $S$-value at least $3a+3$  this is a contradiction.  \hal

\begin{lem}\label{l(r_0-1)+l_r0} $\mu^0 \ne \l_{r_0-1}^0 +\l_{r_0}^0$. 
\end{lem}

\pf  Suppose $\mu^0 = \l_{r_0-1}^0 +\l_{r_0}^0$, so that $s=1$ (by \cite{MF}). Then $V^1(Q_Y) = (\l_{r_0}^0 \otimes \l_{r_0 -1}^0) - \l_{r_0 -2}^0$.  Therefore $V^1 \subseteq  \wedge^2(3\om_1) \otimes 3\om_1$ and $S(V^1) \le 8$.  
Now $V^2(Q_Y) \supseteq (2\l_{r_0 -1}^0 \oplus (\l_{r_0 -2}^0 +\l_{r_0}^0)) \otimes \l_1^1$  and the
first tensor factor can be expressed as $(\l_{r_0 -1}^0  \otimes \l_{r_0 -1}^0 ) - \l_{r_0 -3}^0 $.
Therefore $V^2 \supseteq ((\wedge^2(30) \otimes \wedge^2(30) )- \wedge^4(30)) \otimes (12)$.  Now $\wedge^2(30) \otimes \wedge^2(30)$  contains a composition factor of highest weight $82$ which does not appear in $\wedge^4(30)$.  It follows that $V^2 \supseteq (83)^2$ and an $S$-value consideration gives a contradiction. \hal

\begin{lem}\label{l1+l_r0} $\mu^0 \ne \l_1^0 +\l_{r_0}^0$.
\end{lem}

\pf  Suppose $\mu^0 = \l_1^0 +\l_{r_0}^0$ so that $S(V^1) = 6$.  Then $V^2(Q_Y) \supseteq ((\l_1^0 + \l_{r_0 -1}^0) \oplus \l_{r_0}^0)) \otimes \l_1^1$  which is $(\l_1^0 \otimes \ \l_{r_0 -1}^0) \otimes \l_1^1$.  Therefore $V^2$ contains $3\om_{s+1} \otimes \wedge^2(3\om_s) \otimes (\om_s + 2\om_{s+1})$ which contains $(2\om_{s-1} + 4\om_s +5\om_{s+1}+\om_{s+2})^2$ and an $S$-value comparison gives a contradiction.  \hal

\begin{lem}\label{li} $\mu^0 \ne \l_i^0 $ with $4\le i < r_0$. 
\end{lem}

\pf Suppose false.  Then by \cite{MF}, either $\mu^0 = \l_{r_0-1}^0$ or $s = 1$ and
$4 \le i \le 7$.  First suppose $s \ge 2$ so that  $\mu^0 = \l_{r_0-1}^0$.  Then $S(V^1) \le 6$ and
$V^2 \supseteq \wedge^3(3\om_s) \otimes (\om_s + 2\om_s) \supseteq  (3\om_{s-1} + 3\om_s + 3\om_{s+1}) \otimes (\om_s + 2\om_{s+1})$. By \cite[7.1.5]{MF} this contains $(4\om_{s-1} + 3\om_s + 4\om_{s+1} + \om_{s+2})^2$ of $S$-value $12$ and we have a contradiction.

At this point we can assume $s = 1$,  $L_X' = A_2$,  $V^1 = \wedge^i(03)$, and $V^2 \supseteq \wedge^{i-1}(03) \otimes (12)$. Magma yields the following
\[
\begin{array}{ll}
i = 4: & V^1 = 25+03+22+ 33+60; V^2 \supseteq (23)^5 \\
i = 5: &  V^1 = 25+03+52+14+41+30; V^2 \supseteq (34)^5 \\
i = 6: &  V^1 = 22+30+06+33+52; V^2 \supseteq (34)^5 \\
i = 7: &  V^1 = 60+33+22+00; V^2 \supseteq (23)^7 \\
i = 8: &  V^1 = 41+03; V^2 \supseteq (23)^4. 
\end{array}
\]
In each case we have a contradiction using \cite[5.1.5]{MF}.  \hal

\begin{lem}\label{al_1} If $\mu^0= a\l_1^0$ with $a = 2,3,$ or $4$, then $\l=2\l_1$.  
\end{lem}

\pf  Suppose $\mu^0= a\l_1^0$ with $a = 2,3,$ or $4$ so that $S(V^1) = 3a$.  Note that if $s \ge 2$, then only the case $a =2$ is possible.  We claim that $\la \l, \g_1 \ra = 0 = \la \l, \g_2 \ra$.
 
By way of contradiction suppose $\la \l, \g_1 \ra > 0$.   Then $V^2(Q_Y) \supseteq ((a\l_1^0 + \l_{r_0}^0) \oplus  (a-1)\l_1^0) \otimes \l_1^1$ and the first tensor factor can be expressed
as  $a\l_1^0 \otimes \l_{r_0}^0$.  Therefore $V^2 \supseteq 3a\om_{s+1} \otimes 3\om_s \otimes
(\om_s + 2\om_{s+1}) \supseteq (3\om_s + 3a\om_{s+1}) \otimes (\om_s + 2\om_{s+1}) $ and hence $V^2$ contains  $(\om_{s-1} + 3\om_s + (3a+1)\om_{s+1}+ \om_{s+2})^2$ (using \cite[7.1.5]{MF} as usual).  Therefore $V^2$ has a repeated composition factor of $S$-value at least $3a+4$ and this is a contradiction as
$S(V^1) = 3a$.  Therefore $\la \l, \g_1 \ra = 0$. 

Now suppose $\la \l, \g_2 \ra > 0$.   Then $V^2(Q_Y)$ contains $a\l_1^0 \otimes \l_{r_1-1}^1$ or  $ a\l_1^0 \otimes 2\l_{r_1}^1$ according as  $Y = SO$ or $Y = Sp$. Therefore
$V^2 $ contains $S^a(3\om_{s+1}) \otimes \wedge^2(2\om_s + \om_{s+1})$ or
 $S^a(3\om_{s+1}) \otimes S^2(2\om_s + \om_{s+1})$, respectively.
Magma shows that  $\wedge^2(2\om_s + \om_{s+1}) \supset(5\om_s + \om_{s+2})\oplus ( \om_{s-1} + 3\om_s + \om_{s+1} + \om_{s+2})$, so tensoring with $3a\om_{s+1}$ produces
$(\om_{s-1} + 4\om_s + (3a-1)\om_{s+1} + 2\om_{s+2})^2$ which yields a contradiction by comparing $S$-values.
Similarly, Magma shows that $S^2(2\om_s + \om_{s+1}) \supseteq ( 4\om_s + 2\om_{s+1})\oplus  (\om_{s-1} + 3\om_s + \om_{s+1} + \om_{s+2})$ and here  tensoring with $3a\om_{s+1}$ produces
$(\om_{s-1} + 3\om_s + (3a+1)\om_{s+1} + \om_{s+2})^2$, again a contradiction.
This establishes the claim.

Therefore $\l = a\l_1$. This is excluded by assumption if $a=1$, and is in the conclusion if $a = 2$.  On the other hand for $a = 3$ and $a =4$ we have $X = A_3$ and $Y$ is orthogonal.  Therefore $3\l_1 = S^3(\l_1) - \l_1$  while $4\l_1 = S^4(\l_1) - S^2(\l_1)$.  In either case Magma shows that $V_Y(\l)\downarrow X$ fails to be MF.  \hal
 
 \begin{lem}\label{l1+l2}   $\mu^0 \ne \l_1^0 + \l_2^0$.
\end{lem}

\pf  Suppose $\mu^0 = \l_1^0 + \l_2^0$, so that $s=1$ (by \cite{MF}).  Then $V^1(Q_Y)  = (\l_1^0 \otimes \l_2^0)-\l_3^0$ and Magma shows that $V^1 = (03) \otimes \wedge^2(03) - \wedge^3(03) = 25+14+33+11+22+41+17$.  On the other hand $V^2(Q_Y) \supseteq (2\l_1^0 \oplus \l_2^0) \otimes \l_1^1 = (\l_1^0 \otimes \l_1^0) \otimes \l_1^1$.  Restricting to $L_X'$ this is $03 \otimes 03 \otimes 12$ and this contains $(23)^6$.   At most 3 such composition factors can arise from
$V^1$ so this establishes the result.  \hal

We have now established Theorem \ref{0...0b0...0} when $b=3$.

 \subsection{The case where $b=2$}
 
  In this subsection we complete the proof of  Theorem \ref{0...0b0...0} 
by handling  the remaining case  where $b = 2$. Here   $L_Y' = C^0 \times C^1$ and $L_X'$ is embedded in the factors via the highest weights $2\om_{s+1}$ and $\om_s + \om_{s+1}$, respectively.  In view of the Lemma \ref{setup}, we have $Y = SO$.  Note that when $s = 1$, $C^1 = D_4$ and the image of $L_X'$ in $C^1$ is invariant under triality.
 
 We begin by listing the possibilities for $\mu^0$ and $\mu^1$.  These follow from \cite{MF} in the first case and by the Inductive Hypothesis in the second case.  There are more possibilities if $s = 1$.  
 
\begin{lem}\label{poss12}
 One of the following holds:
\begin{itemize}
\item[{\rm (i)}] $s = 1$,  $C^0 = A_5$, $C^1=D_4$ and
\begin{itemize}
\item[{\rm (a)}]   $\mu^0$ or its dual is one of
\[
 c\l_1^0,\, \l_i^0, \, 2\l_2^0,\, 3\l_2^0,\,\l_i^0+\l_j^0,\, b\l_1^0 + \l_5^0\, (b \le 3),\, b\l_1^0 + \l_2^0\, (b \le 3).
\]
\item[{\rm (b)}]  up to a graph automorphism of $C^1=D_4$, $\mu^1$ is one of
\[
 c\l_1^1\,(c\le 5),\,\l_2^1,\, \l_3^1+\l_4^1. 
\]
\end{itemize}
\item[{\rm (ii)}] $s \ge 2$,  $\mu^0$ or its dual is one of 
\[
0,\, \l_1^0,\, 2\l_1^0,\, \l_2^0,\, \l_3^0\, (s=2),\, \l_1^0 + \l_{r_0}^0 
\]
and  $\mu^1$  is one of
\[
0, \,\l_1^1, \,2\l_1^1,\, \l_2^1.
\]
\end{itemize}
\end{lem}

  \begin{lem}\label{mu0ormu1not0} Either $\mu^0 = 0$ or   $\mu^1 = 0$. 
 \end{lem}
 
 \pf Assume false. Lemma 7.3.1 of \cite{MF} implies that there is a composition factor of $\mu^1 \downarrow L_X'$ with highest weight having at least two nonzero labels.  Then Lemma \ref{stemb} and another application of \cite[Lemma 7.3.1]{MF} shows that either $s \ge 2$ and $\mu^0 \in  \{ \l_1^0, \l_{r_0}^0 \}$ or
 $s = 1$ and  $\mu^0 \in \{ \l_1^0, \l_5^0, \l_3^0, 2\l_1^0, 2\l_5^0 \}$.  We note that in the  $s = 1$ cases  $\mu^0 \downarrow L_X' = 02, 20, 03+30, 04+20, 40+02$, respectively.
 
We first show that $\mu^1 = \l_1^1$ or an image under a graph automorphism if  $C^1 = D_4$.   Suppose $\mu^1 = \l_2^1$, so that $\mu^1 \downarrow  L_X' = \wedge^2(\om_s+\om_{s+1})$.
 If $s = 1$ then $V^1 = (\mu^0\downarrow L_X') \otimes \wedge^2(11)$ and this contains $21^3$, $12^3$, $22^6$, $12^5$, $32^4$ for  the respective possibilities so $V^1$ is  not MF in these cases. On the other hand if $s\ge 2$,  then $V^1 \supseteq (2\om_{s+1}) \otimes ((\om_{s-1} +3\om_{s+1}) \oplus (3\om_s + \om_{s+2}))$ or $V^1 \supseteq 2\om_s \otimes ((\om_{s-1} +3\om_{s+1}) \oplus (3\om_s + \om_{s+2}))$ according as $\mu^0 = \l_1^0$ or $\l_{r_0}^0$.  Hence
  $V^1$ contains $(\om_{s-1} +2\om_s + \om_{s+1} + 2\om_{s+2})^2$ or $(2\om_{s-1} +\om_s + 2  \om_{s+1} + \om_{s+2})^2$, respectively.  Therefore $\mu^1 \ne \l_2^1$.

If $\mu^1 = \l_3^1+\l_4^1$ , then $s=1$ and  $V^1 = (\mu^0\downarrow L_X') \otimes \wedge^3(11)$ and this contains $21^4$, $12^4$, $11^6$, $34^2$, $43^2$, respectively, and again $V^1$ is  not MF.  The same holds for an image of    
$ \l_3^1+\l_4^1$ under a graph automorphism of $C^1 = D_4$, since the image of $L_X'$ is invariant under graph automorphisms.  

The remaining cases are where $s\ge 2$ with $\mu^1 =2\l_1^1$, and
  $s = 1$ with  $\mu^1 = c\l_1^1$ for $2 \le c \le 5 $ (or the image under a graph automorphism).  In the $s = 1$ cases  $(\mu^1 \downarrow L_X') = S^c(11) - S^{c-2}(11)$ and in all cases Magma shows that $V^1$ is not MF. Now suppose $s \ge 2$ with $\mu^1 =2\l_1^1$.  If $\mu^0 = \l_1^0$, then  $V^1 \supseteq 2\om_{s+1} \otimes ((2\om_s + 2\om_{s+1}) \oplus (\om_{s-1} +\om_s + \om_{s+1} + \om_{s+2}))$ and this contains 
  $(\om_{s-1} +\om_s + 3\om_{s+1} + \om_{s+2})^2$.  Similar considerations apply if $\mu^0 = \l_{r_0}^0$.

Therefore $\mu^1 = \l_1^1$ or an image under a graph automorphism if $C^1 = D_4$.  Using \cite[7.3.1]{MF} together with Lemma \ref{stemb} and Magma (for the case where $s=1$),  we see that the only possibilities where $V^1$ is MF are  $\mu^0 = \l_1^0$ or $\l_{r_0}^0$. 

Suppose $\mu^0 = \l_{r_0}^0$. If $s=1$, then $V^1 = 20 \otimes 11 = 31+12+20+01$ whereas $V^2 \supseteq \wedge^2(20) \otimes (S^2(11) - 0) $  if $\mu^1 = \l_1^1$ and $V^2 \supseteq \wedge^2(20) \otimes \wedge^3(11)$  if  $\mu^1 = \l_3^1$ or $\l_4^1$. Magma shows that this contains $(32)^3$ in either case.  Only one such composition factor can arise from $V^1$ so this is a contradiction.

If instead $s\ge 2$, then $V^2 \supseteq 4\om_s \otimes ((3\om_s + \om_{s+2}) \oplus (\om_{s-1} +\om_s + \om_{s+1} + \om_{s+2}))$ and this contains $(\om_{s-1} +5\om_s + \om_{s+1} + \om_{s+2})^2$ which has $S$-value $8$, whereas $S(V^1) =4$, again a contradiction.

Now suppose that $\mu^0 = \l_1^0$.  
Then either $V^2(Q_Y) \supseteq (\l_1^0 + \l_{r_0}^0) \otimes \l_2^1$ or $s = 1$ and
$V^2(Q_Y) \supseteq (\l_1^0 + \l_{r_0}^0) \otimes \l_a^1$ for $a = 3$ or 4.    In the former case, with $s \ge 2$,
$V^2 \supseteq (2\om_s + 2\om_{s+1}) \otimes \wedge^2(\om_s + \om_{s+1})$ and
this contains $(2\om_s + 2\om_{s+1}) \otimes (\om_{s-1} +\om_s + \om_{s+1} +\om_{s+2})$;
this contains $ (2\om_{s-1} +2\om_s + 2\om_{s+1} +2\om_{s+2})^2$ which has $S$-value 8, which is a contradiction as $S(V^1) = 4$.  Finally, if $s = 1$ then Magma shows that
$V^2 \supseteq (22)^4$ whereas $V^1 = 02 \otimes 11 = 13+21+02+10$,  again a contradiction.
\hal
 
 \begin{lem}\label{mu0=0mu1not0}  Suppose $\mu^0 = 0$ and $\mu^1 \ne 0$.  Then $s=1$, $X=A_3$, $Y=D_{10}$ and $\l=\l_9$ or $\l_{10}$.
 \end{lem}
 
  \pf  Suppose $\mu^0 = 0$ and $\mu^1 \ne 0$. 
 
   First assume $\mu^1 = \l_1^1$.  
 Then $V^1 = \om_s + \om_{s+1}$ and $V^2 \supseteq 2\om_s \otimes \wedge^2(\om_s + \om_{s+1})$.  If $s = 1$, then Magma shows that $V^2 \supseteq (12)^3$ a contradiction.
 If $s \ge 2$ then $\wedge^2(\om_s + \om_{s+1}) \supseteq (\om_{s-1} +\om_s + \om_{s+1} +\om_{s+2}) \oplus  (3\om_s +\om_{s+2}) $ and $V^2 \supseteq (\om_{s-1} +3\om_s + \om_{s+1} +\om_{s+2})^2$ and an $S$-value argument gives a contradiction. 
 
 Next suppose $\mu^1 = \l_2^1$.  Then $V^1 = \wedge^2(\om_s + \om_{s+1})$ and $V^2 \supseteq 2\om_s \otimes \wedge^3(\om_s + \om_{s+1})$. If $s = 1$, then Magma shows that $V^2 \supseteq (31)^3$, whereas $V^1 = 11+30+03$.  This is a contradiction.  If $s \ge 2$, then $\wedge^3(\om_s + \om_{s+1}) \supseteq (\om_{s-1} +2\om_s + 2\om_{s+1} +\om_{s+2}) + (2\om_{s-1}  + 3\om_{s+1} +\om_{s+2})$ and hence $V^2 \supseteq (2\om_{s-1} +2\om_s + 3\om_{s+1} +\om_{s+2})^2$.  An $S$-value argument gives a contradiction.
 
 Assume $\mu^1 = 2\l_1^1$ and $s \ge 2$.  Here $S(V^1) = 4$ and $V^2(Q_Y) \supseteq \l_{r_0}^0 \otimes (\l_1^1 +\l_2^1)$.  Now  Lemma \ref{tensorl1}(ii)  implies that $\l_1^1 +\l_2^1 = (\l_1^1 \otimes \l_2^1) - \l_3^1 -\l_1^1$.  Now
  $\wedge^2(\om_s +\om_{s+1})\supseteq (3\om_s+\om_{s+2}) \oplus (\om_{s-1}+3\om_s)$ 
and  the summands arise from the weights $2(\om_s + \om_{s+1}) - \a_{s+1}$ and $2(\om_s + \om_{s+1}) - \a_s$, respectively.  It follows that the restriction of  $\l_1^1 +\l_2^1$ to $L_X'$ contains composition factors of highest weights  $3(\om_s +\om_{s+1}) -\a_{s+1}$  and $3(\om_s +\om_{s+1}) -\a_s$,
as these cannot occur in the restriction of $\l_1^1$ or $\l_3^1$.  Tensoring with $2\om_s$ 
we conclude that $V^2 \supseteq (2\om_{s-1} +2\om_s + 3\om_{s+1} +\om_{s+2})^2$.  An $S$-value argument gives a contradiction.  

 For the rest of this proof we can assume $s = 1$. First assume $\mu^1 = c\l_i^1$ with $2 \le c \le 5$ and $i = 1,3,4$.  Then $V^1 = S^c(11)-S^{c-2}(11)$ and has $S$-value $2c$. Suppose $i = 1$. Then
 $V^2(Q_Y) \supseteq \l_5^0 \otimes ((c-1)\l_1^1 + \l_2^1)$.  Magma calculations show that the restriction of the second tensor factor contains $(cc)^2$ and hence $V^2 \supseteq ((c+2)c)^2$ which gives an $S$-value contradiction.  
 
 Since $L_X'$ is invariant under graph automorphisms
 we can now assume $i = 3$ so that $V^2(Q_Y) \supseteq \l_5^0 \otimes ((c-1)\l_3^1 + \l_4^1)$.
 Magma shows that $ (c-1)\l_3^1 + \l_4^1 = ((c-1)\l_3^1 \otimes  \l_4^1) - (\l_1^1 +(c-2)\l_3^1)$. Using this and $S$-values we use Magma to see that $V^2$ contains $(31)^3$, $(15)^2$, $(26)^2,$ or  $(37)^2$, according as $c = 2,3,4$ or $5$, respectively. The composition factors of
 $V^1$ are given below:
 
 $c = 2: 22+11$.
 
 $c = 3: 33 + 22+ 03 + 30 + 00$.
 
 $c = 4: 44 + 33 + 41+14 +22+11$.
 
 $c = 5: 55+44+33+22+11+41+14+52+25$.
 
 \noindent In each case we see that there are not enough repeated composition factors of $V^2$ arising from $V^1$, 
 so this is a contradiction.
 
 Assume $\mu^1 = \l_3^1 + \l_4^1$.  Then $V^1 = \wedge^3(11) = 22 + 11 + 30 + 03 + 00$
 and $V^2(Q_Y) \supseteq (\l_5^0 \otimes 2\l_3^1) \oplus (\l_5^0 \otimes2\l_4^1)$.  Therefore $V^2 \supseteq 20 \otimes ((22)^2 + (11)^2) \supseteq (42)^2$.  As $42$ cannot arise from $V^1$ this is a contradiction.
 Now assume that $\mu^1 = \l_1^1 + \l_3^1$ or $\l_1^1 + \l_4^1$.  It suffices to settle the former.
 Then $V^1 = \wedge^3(11)$ as above and $V^2(Q_Y) \supseteq \l_5^0 \otimes (\l_2^1 + \l_3^1)$.
 Magma shows that $\l_2^1 + \l_3^1 = (\l_2^1 \otimes \l_3^1) - (\l_1^1 + \l_4^1) - \l_3^1$ and that $V^2$ contains $(42)^3$ again a contradiction.
 
 The remaining cases are where $\mu^1 = \l_3^1$ or $\l_4^1$ and it will suffice to settle the first
 of these.  We claim that $\la \l, \g_1 \ra = 0$.  Otherwise $V^2 \supseteq 20 \otimes \wedge^3(11) \supseteq (31)^2$.  However, $V^1 = 11$, so this is a contradiction.  This gives the claim.
 At this point $V_Y(\l)$ is a spin module, as in the conclusion. \hal

By Lemmas \ref{mu0ormu1not0} and  \ref {mu0=0mu1not0}, we may assume that $\mu^1 = 0$. We now consider the various possibilities for $\mu^0$ given by Lemma \ref{poss12}.
 
 \begin{lem}\label{l123} Assume $\mu^0 = \l_1^0$, $\l_2^0$,  $\l_3^0$ or $2\l_1^0$.  Then $\l = \l_2$, 
$2\l_1$ or $\l_3$ ($s = 1$).
\end{lem}

\pf  Assume the hypotheses.  In the case $\mu^0 = \l_3^0$ we have $s \le 2$ by Lemma \ref{poss12}.
We first claim that $\la \l, \g_1\ra = 0$.  Suppose otherwise.  Then $S(V^1) \le 2, 4, 6,$ or $4$, respectively
($2,3,3$ or 4 if $s=1$).  On the other hand  by  Lemma 5.4.1 of \cite{MF}, $V^2(Q_Y) \supseteq (\mu^0 \otimes \l_{r_0}^0) \otimes \l_1^1$.
As in other places we will carry out the computations assuming $s>1$ but we can simply drop terms to get the $s=1$ result.

If $\mu^0 = \l_1^0$  then $V^2 \supseteq 2\om_{s+1} \otimes  2\om_s \otimes (\om_s +\om_{s+1})$ and this contains $(\om_{s-1} + 2\om_s + 2\om_{s+1}+\om_{s+2})^2$ which has
$S$-value at least 4 and this is a contradiction.  If   $\mu^0 =2\l_1^0$, then 
$V^2 \supseteq 4\om_{s+1} \otimes  2\om_s \otimes (\om_s +\om_{s+1})$ and this contains $(\om_{s-1} + 2\om_s + 4\om_{s+1}+\om_{s+2})^2$, again a contradiction.

The remaining cases are $\mu^0 = \l_2^0$ and $\l_3^0$  and $V^2 \supseteq (\mu^0 \downarrow L_X') \otimes (3\om_s + \om_{s+1})$.  Now $\mu^0 \downarrow L_X'$ contains $\om_s +2\om_{s+1}+\om_{s+2}$ or  $2\om_s+3\om_{s+1} + \om_{s+3}$, respectively.  Therefore,
$V^2 \supseteq  (\om_{s-1} + 3\om_s +2\om_{s+1} +2\om_{s+2})^2$ or  $(\om_{s-1} + 4\om_s +3\om_{s+1} +\om_{s+2}+ \om_{s+3})^2$, respectively.  Now we get a contradiction in each case
using $S$-values.  

This establishes the claim that $\la \l, \g_1\ra = 0$. Hence $\l = \l_1, \l_2, \l_3$ ($s = 1,2$) or $2\l_1$. 
Note finally that $\l_1$ is excluded by assumption, and $V_Y(\l_3)\downarrow X$ is not MF when $s=2$ by \cite{MF}. This completes the proof. \hal
 
  \begin{lem}\label{l4}  $\mu^0 \ne \l_4^0$. 
\end{lem}

\pf  Assume $\mu^0 = \l_4^0$.  Then Lemma \ref{poss12} implies that $X = A_3$ and so $V^1 = \wedge^2(20) = 21$.  If $\la \l, \g_1 \ra > 0$, then by  Lemma 5.4.1 of \cite{MF}, $V^2(Q_Y) \supseteq 
 \l_4^0 \otimes \l_5^0  \otimes \l_1^1$.  But then $V^2 = \wedge^2(20) \otimes 20 \otimes 11 \supseteq (41)^2$. Comparing  $S$-values gives a contradiction.  Therefore  $\la \l, \g_1 \ra = 0$ and $\l = \l_4$.  Now the main result of \cite{MF}  gives
 the assertion. \hal

\begin{lem} $\mu^0 \ne \l_{r_0-1}^0$ or $ \l_{r_0-2}^0$.
\end{lem}

\pf Assume false. If $s = 1$ these cases are $\mu^0 = \l_4^0$ and $\mu^0 =\l_3^0$. These are covered by 
Lemmas \ref{l4} and \ref{l123}.

Now assume $s\ge 2$. The case $\l_{r_0-2}$ only occurs for $s=2$, where we obtain a contradiction using Magma. Now consider 
$\l_{r_0-1}$. Here $V^2 \supseteq \wedge^3(2\om_s) \otimes (\om_s + \om_{s+1})$.
This contains $(3\om_s + 2\om_{s+1}) \otimes (\om_s +\om_{s+1})$, which contains $(\om_{s-1}+3\om_s + 2\om_{s+1} +\om_{s+2})^2$, of $S$-value 7. On the other hand $V^1 = \wedge^2(2\om_s)$ has
$S$-value 4 and we again obtain a contradiction. \hal

   \begin{lem}\label{cl1}  If $\mu^0 = c\l_1^0$ for $c \ge 3$, then $\l=3\l_1$ and $X=A_3$.  
 \end{lem}
 
 \pf  Suppose $\mu^0 = c\l_1^0$ with $c\ge 3$ so that $V^1 = S^c(2\om_{s+1})$. Then Lemma \ref{poss12} implies
 that $X = A_3$.  We  claim that $\la \l, \g_1 \ra = 0$.  
 Otherwise Lemma 5.4.1 of \cite{MF} shows that $V^2(Q_Y) \supseteq c\l_1^0 \otimes \l_{r_0}^0  \otimes \l_1^1$.  Therefore $V^2 \supseteq (\om_{s-1} + 2\om_s +2c\om_{s+1} + \om_{s+2})^2$
 which has $S$-value at least $2+2c$.  This is a contradiction proving the claim.  
 
It follows that $\l = c\l_1$.  
 Suppose that $c \ge 4$.  Then $V^2 \supseteq S^{c-1}(02) \otimes 11 \supseteq (2(2c-6)) \otimes 11 \supseteq (2(2c-6))^2$
 and $(2(2c-6))$ has $S$-value $2c-4$.  But using \cite[\S 3.1]{Howe} (or \cite[\S 6.5]{MF}), we see that the only composition factors of $S$-value at least $2c-5$  in $V^1$ are $(0(2c)), (2(2c-4))$ and $(4(2c-8))$ and none of these contributes a term $(2(2c-6))$ to $V^2$.  Hence $c=3$, completing the proof of the lemma. \hal

\begin{lem}\label{cl5}  $\mu^0 \ne c\l_{r_0}$ with $c>0$.  
\end{lem}

\pf Assume $\mu^0 = c\l_{r_0}$ with $c>0$. Here $V^1 = S^c(2\om_s)$ and $V^2(Q_Y)\supseteq (\l_{r_0-1}^0 +(c-1)\l_{r_0}^0) \otimes \l_1^1$.

First assume that $c = 1$.  We  claim that $\la \l, \g_1 \ra = 0$.  If not then Lemma 5.4.1 of \cite{MF} shows that $V^2 \supseteq 2\om_s \otimes 2\om_s \otimes (\om_s +\om_{s+1})$.  If $s = 1$ Magma shows that $V^2$ contains $(21)^4$ and only one such factor can arise from $V^1$.  If $s \ge 2$  then $V^2$ contains $(2\om_{s-1} + 2\om_s+\om_{s+1}+ \om_{s-1})^2$  which has $S$-value 6 contradicting the
fact that $V^1$ has $S$ value 2.  This establishes the claim and hence $\l = \l_{r_0}$.  Then Lemma \ref{symwed} together with \cite{MF} shows that $V_Y(\l)\downarrow X$ is not MF.  

From now on assume $c \ge 2$.
Corollary 4.1.3 of \cite{MF} shows that the first tensor factor of $V^2(Q_Y)$ equals 
$(\l_{r_0-1}^0 \otimes (c-1)\l_{r_0}^0)-(\l_{r_0-2}^0 +(c-2)\l_{r_0}^0)$.  
The first tensor product restricts to $L_X'$ as $S^{c-1}(2\om_s) \otimes \wedge^2(2\om_s)$ and the restriction of the subtracted term
is contained in $S^{c-2}(2\om_s) \otimes \wedge^3(2\om_s)$.

First assume that $s = 1$.  Then $S^{c-1}(2\om_s) \otimes \wedge^2(2\om_s)$  contains $((2c-2)0) \otimes (21) \supseteq ((2c)1) + ((2c-2)2)$.
These summands have $S$-value at least $2c$.  On the other hand the $S$-value of $S^{c-2}(2\om_s) \otimes \wedge^3(2\om_s)$ is  $2(c-2) + 3 = 2c-1$.  It follows that $V^2 \supseteq ((2c)1) + (2c-2)2) \otimes 11 \supseteq ((2c)1)^3$.
Only one term $(2c)1$ can arise from $V^1 = S^c(20)$ as all summands other than $(2c)0$ have $S$-value strictly less than $2c$.  It follows that  $V_Y(\l)\downarrow X$ is not MF.

Finally we must consider the case where $c \ge 2$ and $s \ge 2$.  The only possibility is where $c = 2$.  Therefore
$V^1 = S^2(2\om_s)$ where all irreducible summands have $S$-value at most 4.  On the other hand
$V^2 \supseteq (2\om_s \otimes \wedge^2(2\om_s))\otimes (\om_s + \om_{s+1}) -(\wedge^3(2\om_s) \otimes (\om_s + \om_{s+1}))$.  Now $\wedge^2(2\om_s)$
contains $(\om_{s-1} + 2\om_s + \om_{s+1})$, and from this, using the parabolic argument together with Magma,  we see that the first summand contains $(\om_{s-1}+5\om_s+2\om_{s+1}-\a_s-\a_{s+1})^4$ which is equal to $(2\om_{s-1} + 4\om_s + \om_{s+1} + \om_{s+2})^4$, while only one such composition factor can appear in the subtracted term.   Therefore  $V^2$ contains a repeated composition factor of  $S$-value $8$, giving the final contradiction.  \hal

\begin{lem}\label{bl1+l5}  $\mu^0 \ne b\l_1^0 + \l_{r_0}^0$ with $1 \le b \le 3$.
\end{lem}

\pf Suppose $\mu^0 = b\l_1^0 + \l_{r_0}^0$ with $1 \le b \le 3$. We note that $s = 1$ if $b > 1$.
Then $V^2(Q_Y) \supseteq ((b\l_1^0 + \l_{r_0 -1}^0) \oplus  ((b-1)\l_1^0 + \l_{r_0}^0)) \otimes \l_1^1$ and this
equals $b\l_1^0 \otimes \l_{r_0 -1}^0 \otimes \l_1^1$. 

Restricting to $L_X'$ we see that $V^2 \supseteq S^b(2\om_{s+1}) \otimes \wedge^2(2\om_s) \otimes (\om_s + \om_{s+1})$.  First suppose $b = 1$.  Then $V^1 = (2\om_{s+1} \otimes 2\om_s)-0$ which has $S$-value 4, and
$V^2 \supseteq (2\om_{s+1} \otimes \wedge^2(2\om_s)) \otimes (\om_s + \om_{s+1}) \supseteq
(2\om_{s-1} + 2\om_s + 3\om_{s+1} + \om_{s+2})^2$.  If $s\ge 2$, this  has $S$-value 8, which gives a contradiction.  On the other hand if $s=1$, then $V^2 \supseteq (31)^5$ while $V^1 = (22) + (11)$ and this is again a contradiction.

What remains are the cases where $s = 1$ and $b = 2$ or $3$.   Then $V^2 \supset
 S^2(02) \otimes (21) \otimes (11)$ or $V^2 \supseteq S^3(02) \otimes (21) \otimes 11$, respectively.  A Magma computation  shows that $V^2 \supseteq (22)^{12}$ or $(13)^{16}$, respectively.  On the other hand $V^1 \subseteq  S^b(02) \otimes 20$ and another Magma argument shows that in either case at most 4 of these composition factors can arise from $V^1$.  This is a contradiction. \hal

\begin{lem}\label{l5+bl1} $\mu^0 \ne \l_1^0 + b\l_5^0$ with $b = 2$ or $ 3$.  
\end{lem}

\pf  Assume $\mu^0 = \l_1^0 + b\l_5^0$ with $b = 2$ or $b = 3$, so that $s = 1$.
Then $V^2(Q_Y) \supseteq ((\l_1^0 + \l_4^0 +(b-1)\l_5^0) \oplus b\l_5^0) \otimes \l_1^1$.  
Now  $\l_1^0 \otimes (\l_4^0 +(b-1)\l_5^0) = (\l_1^0 + \l_4^0 +(b-1)\l_5^0) \oplus  b\l_5^0 \oplus (\l_4^0 +(b-2)\l_5^0)$.  It follows from Table 7.2 of \cite{MF} and an $S$-value argument that $(\l_1^0 + \l_4^0 +(b-1)\l_5^0) \downarrow L_X' \supseteq ((2b)3)+((2b+1)1)$.  Now tensoring with
$\l_1^1 \downarrow L_X' = (11)$ we see that $V^2 \supseteq ((2b)3)^3$.  At most one such
summand can arise from $V^1$, so this establishes the result.  \hal

\begin{lem}\label{2l4,3l4}  $\mu^0 \ne 2\l_4^0$,  $3\l_4^0$, $2\l_2^0$ or $3\l_2^0$.  
\end{lem}

\pf Assume false. Then $s = 1$ by Lemma \ref{poss12}.  First assume  $\mu^0 = 2\l_4^0$ or $3\l_4^0$.    It follows from  Table 7.2 of \cite{MF} that $V^1 = 42+31+04+20$ or $63+25+52+60+33+41+22+30+03$, respectively.

Now $V^2(Q_Y) \supseteq (\l_3^0 + \l_4^0) \otimes \l_1^1$ or $(\l_3^0 + 2\l_4^0) \otimes \l_1^1$,
respectively.  Restricting to $L_X'$ and using  Magma we find that $V^2 \supseteq (51)^3$ or $(64)^3$, respectively. (In the second case
we argue that $\l_3^0 + 2\l_4^0 = (S^2(\l_4^0) \otimes \l_3^0)  - (\l_1+\l_4) - (\l_2+\l_4+\l_5)$ and
upon restriction to $L_X'$ and tensoring with $(11)$ $S$-value considerations show that  the term (64) does not arise from the first  subtracted term and  can occur at most once in the second subtracted term.) But $V^1$ can contribute at most one such composition factor,
so this completes the argument for these cases.  

Now assume  $\mu^0 = 2\l_2^0$ or $3\l_2^0$, so that $V^1$ is the dual of what appears
in the first paragraph.  Magma computation shows that $V^2 \supseteq (14+22+11) \otimes (11) \supseteq (22)^3$ or $(42+26+34+23+20+12+31+15+04^2) \otimes (11) \supseteq (15)^6$, respectively.  We obtain a contradiction in each case.  \hal

\begin{lem}\label{bl1+l2} $\mu^0 \ne b\l_1^0+\l_2^0$ or $\l_4^0+b\l_5^0$ with $1 \le b \le 3$.  
\end{lem}

\pf Assume false. Then $s = 1$.  We begin with the case $\mu^0 = b\l_1^0+\l_2^0$.  Table 7.2 
of \cite{MF} shows that $V^1  = 14+22+11,$ $02+21+24+32+13+16$, or $34+23+12+31+20+15+26+04+42+18$, according as $b = 1$, $2$, or $3$.

On the other hand $V^2(Q_Y) \supseteq ((b+1)\l_1^0 \oplus ((b-1)\l_1^0 + \l_2^0)) \otimes \l_1^1$
and this equals $b\l_1^0 \otimes \l_1^0 \otimes \l_1^1$.  Restricting to $L_X'$
this becomes $S^b(02)\otimes (02) \otimes (11)$ and Magma shows that this contains
$(12)^4$, $(14)^5$, or $(13)^8$ according as $b = 1$, $2$, or $3$.  However,  we see from the above that in each case $V^1$ can contribute at most $2, 2,$ or $3$ such composition factors according as $b = 1,2,$ or $3$. Therefore we have a contradiction.  

Now assume that $\mu^0 = \l_4^0+b\l_5^0$.  Here $V^1$ is the dual of that given in the first paragraph and  $V^2(Q_Y) \supseteq (2\l_4^0 + (b-1)\l_5^0) \otimes \l_1^1$. If $b \ge 2$, then starting from $2\l_4^0 \otimes (b-1)\l_5^0$ and using  Magma together with $S$-value arguments we find that $V^2 \supseteq ((2b+2)2 + (2b)3) \otimes 11) \supseteq ((2b+2)2)^3$, whereas at most one such composition factor can arise from $V^1$.  Suppose $b = 1$.  Then $V^2(Q_Y) \supseteq (2\l_4^0 \oplus (\l_3^0 +\l_5^0)) \otimes \l_1^1$.
Restricting to $L_X'$  we see that $V^2 \supseteq (42 + 31) \otimes 11$ and this contains
$(42)^3$ and we have the same contradiction.  \hal

\begin{lem}\label{l1+l3,l3+l5} $\mu^0 \ne \l_i^0+\l_j^0$ for $(i,j) \in \{(1,3),\,(1,4),\,
(2,3),\,(2,4),\,(2,5),\,(3,4),\,(3,5)\}$.
\end{lem}

\pf  Assume false, so $s = 1$. For each possibility $\mu^0 = \l_i^0+\l_j^0$ with $i\le 2$, we can read off the composition factors in $V^1$ from \cite[Table 7.2]{MF}, and those for $i=3$ are duals of these. The results are listed in Table \ref{mu0lists}. Also in this table are a summand of $V^2(Q_Y)$, and a repeated composition factor of $V^2$, computed using Magma. In all cases, this gives a contradiction using Proposition \ref{induct}. \hal

\begin{table}[h!] 
\caption{}\label{mu0lists}
\[
\begin{array}{|cccc|}
\hline
\mu^0 & V^2(Q_Y) \supseteq & V^1 & V^2 \supseteq \\
\hline
\l_1^0+\l_3^0 & \l_1^0\otimes \l_2^0\otimes \l_1^1 & 05 + 21 + 32+13+10 & 22^6 \\
\l_1^0+\l_4^0 & \l_1^0 \otimes \l_3^0 \otimes \l_1^1 & 23+31+12+01 & 21^7 \\
\l_2^0+\l_3^0 & ((\l_2^0 \otimes \l_2^0)-\l_4^0) \otimes \l_1^1 & 15+42+23+04+31+12+20 & 32^7 \\
\l_2^0+\l_4^0 & ((\l_2^0 \otimes \l_3^0)-\l_5^0) \otimes \l_1^1 & 33+41+14+22+30+03+11 & 12^{11} \\
\l_2^0+\l_5^0 &  ((\l_2^0 \otimes \l_4^0)-0) \otimes \l_1^1 & 32+13+21+10 & 22^{11} \\
\l_3^0+\l_4^0 & ((\l_3^0 \otimes \l_3^0)-(\l_1^0 \otimes \l_5^0)) \otimes \l_1^1 & 51+24+32+40+13+21+02 & 22^{13} \\
\l_3^0+\l_5^0 & ((\l_3^0 \otimes \l_4^0) -\l_1^0) \otimes \l_1^1 & 50 + 12 + 23+31+01 & 21^{11} \\
\hline
\end{array}
\]
\end{table}

 
 \begin{lem}\label{mu0=0}  $\mu^0 \ne 0.$
 \end{lem}
 
 \pf  Suppose $\mu^0 = 0$ so that $V^1 = 0$.  Here we allow all values of $s$.  Let $c = \la \l, \g_1 \ra > 0$.  If $c = 1$, then
 $\l = \l_{r_0+1}$, $V_Y(\l) \downarrow X = \wedge^{r_0+1}(\d)$ and \cite{MF} shows that
 this is not MF.  So now assume $c > 1$. Then $V^2 = 2\om_s  \otimes (\om_s + \om_{s+1})$ and all composition factors have $S$-value at most 4. On the other hand $V^3(Q_Y)$ has a composition factor afforded by $\l - 2\g_1$ which
 restricts to $S^2(2\om_s) \otimes (S^2(\om_s + \om_{s+1})) - 0)$.  If $s > 1$ this contains $(4\om_s \oplus (2\om_{s-1} + 2\om_{s+1})) \otimes (2\om_s + 2\om_{s+1})$ and hence $V^3 \supseteq (2\om_{s-1} + 2\om_s + 4\om_{s+1})^2$ of $S$-value $8$.  And if $s = 1$, $V^3$ contains $S^2(20) \otimes (S^2(11) - 0)$ and this contains $(24)^2$ of $S$-value 6.  In either case we have a contradiction.  \hal
 
 This completes the analysis of the $b=2$ case and hence completes the proof of Theorem \ref{0...0b0...0}.

\section{Proof of Theorem \ref{mainthm}, part I: even rank} \label{secevenrk}

In this section we complete the proof of our main Theorem \ref{mainthm} in the case where $X=A_{2s}$ has even rank.
Assume that $X = A_{2s}$ and let $W = V_X(\d)$ be a self-dual module for $X$, so that 
\begin{equation}\label{evenwt}
\d = (a_1,\cdots ,a_s,a_s,\cdots ,a_1).
\end{equation}
In view of Theorems \ref{b0...0b} and \ref{omi+istar}, we may assume that
\begin{equation}\label{noteq}
\d \ne b\o_1+b\o_{2s} \hbox{ or }\o_i+\o_{2s+2-i}.
\end{equation}
In particular, $s\ge 2$ and $X$ has rank at least 4. 
Note that $Y$ is an orthogonal group, by Lemma \ref{setup}.

We shall prove

\begin{thm}\label{evenrk} Let $X=A_{2s}$, let $W=V_X(\d)$ with $\d$ as in $(\ref{evenwt})$ and $(\ref{noteq})$, 
and let $X<Y=SO(W)$. If $\l \ne 0,\l_1$ is a dominant weight for $Y$, then $V_Y(\l)\downarrow X$ is not MF.
 \end{thm}

We now embark on the proof. 
Assume the hypotheses of the theorem, so that $X=A_{2s}$ and $W = V_X(\d)$ with $\d$ as in $(\ref{evenwt})$ and $(\ref{noteq})$.  Suppose $\l \ne 0,\l_1$ is a dominant weight for $Y$ such that $V_Y(\l) \downarrow X$ is MF. We aim for a contradiction.

By Lemma \ref{setup}(i), we have $L_Y'  = C^0\times \cdots \times C^k$ with $k=\sum_1^sa_i\ge 2$. Define $j$ to be maximal such that $a_j\ne 0$.

\begin{lem}\label{wk1cf} $C^k$ is an orthogonal group and $V_{C^k}(\l_1^k) \downarrow L_X' \supseteq \nu \oplus \nu'$, where 
\[
\begin{array}{l}
\nu = (a_1, \cdots ,a_{s-1},(2a_s),a_{s-1}, \cdots ,a_1), \\
\nu' = \left\{ \begin{array}{l} (a_1,\cdots,a_{j-1}+1,a_j-1,0,\cdots ,0,a_j-1,a_{j-1}+1,\cdots,a_1), \hbox{ if }j<s \\
(a_1,\cdots,a_{s-1}+1,2(a_s-1),a_{s-1}+1,\cdots,a_1), \hbox{ if }j=s 
\end{array}
\right.
\end{array}
\]
\end{lem}

\pf As in the proof of Lemma \ref{setup}, $\nu$ arises from application of $\bar f_j^{a_j}\cdots \bar f_1^{a_1}$ to a maximal vector, and similarly $\nu'$ arises from applying $f_j\bar f_j^{a_j-1}\cdots \bar f_1^{a_1}$. \hal

\begin{lem}\label{muk01s} We have $\mu^i=0$ for $1\le i\le k-2$ and also $\mu^k=0$.
\end{lem}

\pf Assume $\mu^k \ne 0$.  The summand of highest weight $\nu$ occurs with multiplicity 1 and hence the image of $L_X'$ in $C^k$ is contained in a product of orthogonal groups $R \times S$, where the restriction of the natural module for $R$ to  $L_X'$ affords $\nu$.  Then applying the argument of Lemma \ref{possmu0mu1} and the Inductive Hypothesis we see that $\mu^k = \l_1^k, \l_2^k$, or $2\l_1^k$. 
If $\mu^k = \l_2^k$ or $2\l_1^k$, the argument of \cite[Lemma 14.2.1]{MF} shows that $V^1$ is not MF, a contradiction. 

Hence $\mu^k = \l_1^k$.   We have $S(V^1) = S(\mu^0 \otimes \cdots \otimes \mu^{k-1} \downarrow L_X' ) + S(\nu)$. Denote the first summand by $z$.  Now $V^2(Q_Y)$ has a summand $\mu^0 \otimes \cdots \otimes \mu^{k-2} \otimes (\mu^{k-1} + \l_{r_{k-1}}^{k-1}) \otimes \l_2^k$.  Using Lemma 3.9  of \cite{MF} we see that the $S$-value of the restriction to $L_X'$ is at least $z + S(\l_2^k \downarrow L_X')$.  Then \cite[14.2.1]{MF} shows that $S(\l_2^k \downarrow L_X')$ contains a repeated summand of $S$-value $2S(\nu)$. Hence Lemma \ref{prop38} gives $z + 2S(\nu) \le z + S(\nu)+1$, a contradiction.

Hence $\mu^k=0$. Now the proof of \cite[Theorem 14.2]{MF} shows that also $\mu^i=0$ for $1\le i\le k-2$.  \hal

\begin{lem} \label{mk10} We have $\mu^{k-1}=0$.
\end{lem}

\pf Assume $\mu^{k-1}\ne 0$. Then \cite[Lemma 14.1.5]{MF} implies that either $\mu^{k-1}$ is the natural or dual natural module for $C^{k-1}$, or $s=2$ and  level $k-1$ restricts to $L_X'$ as $020\oplus 101$. However, when $s=2$ we have $\d = a_1a_2a_2a_1$ with $a_2 \ge 1$, and level $k-1$ contains a summand $(a_1,2a_2-1,a_1+1)$ (from application of $f_3^{a_2-1}f_4^{a_1}$). Hence $\mu^{k-1}$ is the natural module or its dual.

The arguments in the proof of \cite[Theorem 14.2]{MF} rely only on the local analysis of summands of $V^2$ arising from the weight $\l - \g_{k-1} - \b_1^{k-1}$. This  gives a contradiction in the case where 
$\mu^{k-1}$ is the natural module $\l_1^{k-1}$. So assume now that $\mu^{k-1} = \l_{r_{k-1}}^{k-1}$. Writing $W^{i+1}$ for the restriction of the $i^{th}$ level to $L_X'$ (i.e. for $V_{C^i}(\l_1^i)\downarrow L_X'$), we see that 
\[
V^2 \supseteq (\mu^0\downarrow L_X') \otimes \wedge^2(W^k)^* \otimes W^{k+1}.
\]
If $\nu_1$ denotes a summand of $W^k$ of highest $S$-value, then \cite[Lemma 14.2.1]{MF} shows that 
$\wedge^2(W^k)^*$ has a repeated summand of $S$-value at least $2S(\nu_1)-1$. Hence $V^2$ has a repeated summand of $S$-value at least $S(\mu^0\downarrow L_X')+2S(\nu_1)+1$, so Lemma \ref{prop38} implies that $2S(\nu_1)+1 \le S(\nu_1)+1$, a contradiction. \hal

\begin{lem} \label{gams} We have $\la \l,\g_y\ra = 0$ for all $y$.
\end{lem}

\pf  By Lemmas \ref{muk01s} and \ref{mk10}, we have $V^1 = \mu^0\downarrow L_X'$. 
Section 16 of \cite{MF} shows that $\la \l,\g_y\ra = 0$ for $y\le k-1$. Now suppose $\la \l,\g_k\ra \ne 0$. 

Assume first that $a_s \ne 0$. Then there is an $L_X'$-summand  at level $k-1$ of weight 
$\nu_2 = (a_1,\cdots,a_{s-1},2a_s-1,a_{s-1}+1,a_{s-2}, \cdots a_1)$.  This is afforded by $\bar f_s^{a_s -1}\bar f_{s-1}^{a_{s -1}}\cdots \bar f_1^{a_1}$.  Note that $\nu_2$ has maximal $S$-value at level $k-1$ since we are not applying a power of $f_1$ (see \cite[Cor. 5.1.4]{MF}). 
Hence, considering the summand of $V^2(Q_Y)$ afforded by $\l-\g_k$, we see that 
$V^2 \supseteq (\mu^0\downarrow L_X') \otimes \nu_2^*\otimes \nu'$ (where $\nu'$ is the summand of level $k$ given in Lemma \ref{wk1cf}). By \cite[Lemma 7.1.7]{MF}, this has a repeated summand of $S$-value at least 
$S(\mu^0\downarrow L_X')+S(\nu_2)+S(\nu')-2$, and this is greater than $S(V^1)+1$, contradicting Lemma \ref{prop38}.

Now assume that $a_s=0$, and recall that $j$ is maximal such that $a_j\ne 0$. There is an $L_X'$-summand in level $k$ of highest weight $\nu = (a_1,\cdots ,a_j,0,\cdots ,0,a_j,\cdots ,a_1)$ (see Lemma \ref{wk1cf}), and one in level $k-1$ of highest weight $\nu_3=(a_1,\cdots ,a_j,0,\cdots ,0,a_j-1,a_{j-1}+1,\cdots ,a_1)$ afforded by $\bar f_j^{a_j -1}\bar f_{j-1}^{a_{j -1}}\cdots \bar f_1^{a_1}$.  As $j > 1$ this has maximal $S$-value at level $k-1$  (see \cite[Cor. 5.1.4]{MF} again).  
Hence $V^2$ contains a summand $(\mu^0\downarrow L_X') \otimes \nu_3^* \otimes \nu$. By \cite[Lemma 7.1.7]{MF}, 
$\nu_3^* \otimes \nu$ has a repeated summand of $S$-value at least $S(\nu_3)+S(\nu)-2$, and this leads to the usual contradiction to Lemma \ref{prop38}. \hal 

\vspace{2mm}
At this point we can complete the proof of Theorem \ref{evenrk}. By the previous lemmas we have $\l = \mu^0$. Recall from Section 3 that $\d' = (a_1,\cdots ,a_s,a_s,\cdots ,a_2)$ and $s\ge 2$. By hypothesis, either $L(\d')\ge 3$, or $L(\d')=2$ and $a_j\ge 2$. Hence the main result of \cite{MF} implies that $\mu^0 = \l_1^0$ or $\l_{r_0}^0$. The latter case is ruled out by Section 15 of \cite{MF}, and the former case gives $\l = \l_1$, which is excluded in the hypothesis of the theorem. This completes the proof.

\section{Proof of Theorem \ref{mainthm}, part II: odd rank}

In this section we complete the proof of our main Theorem \ref{mainthm}, by handling the case where $X=A_{2s+1}$ has odd rank.
Assume that $X = A_{2s+1}$ and let $W = V_X(\d)$ be a self-dual module for $X$, so that 
\begin{equation}\label{oddwt}
\d = (a_1,\cdots ,a_s,\,c,\,a_s,\cdots ,a_1).
\end{equation}
In view of the results in the previous sections, we may assume that
\begin{equation}\label{noteq1}
\d \ne b\o_{s+1},\; b\o_1+b\o_{2s+1} \hbox{ or }\o_i+\o_{2s+2-i}.
\end{equation}
By Lemma \ref{setup}, $Y$ is a symplectic group if $c$ is odd and $s$ is even; otherwise $Y$ is orthogonal.

We shall prove

\begin{thm}\label{oddrk} Let $X=A_{2s+1}$, let $W=V_X(\d)$ with $\d$ as in $(\ref{oddwt})$ and $(\ref{noteq1})$, 
and let $X<Y=Sp(W)$ or $SO(W)$. If $\l \ne 0,\l_1$ is a dominant weight for $Y$, then $V_Y(\l)\downarrow X$ is not MF.
 \end{thm}

We now embark on the proof. 
Assume the hypotheses of the theorem, so that $X=A_{2s+1}$ and $W = V_X(\d)$ with $\d$ as in $(\ref{oddwt})$ and $(\ref{noteq1})$. Note that this implies $s\ge 1$. Suppose $\l \ne 0,\l_1$ is a dominant weight for $Y$ such that $V_Y(\l)\downarrow X$ is MF. We aim for a contradiction.

By Lemma \ref{setup}(ii), we have $L_Y'  = C^0\times \cdots \times C^k$ with 
$k=\sum_1^sa_i +\lfloor \frac{c}{2} \rfloor$. Also $C^k$ is an orthogonal group if $c$ is even, and is of type $A$ if $c$ is odd.

\begin{lem}\label{muk0lem} Assume the following conditions hold:
\begin{itemize}
\item[{\rm (i)}] $\mu^k = 0$, and 
\item[{\rm (ii)}] if $C^k$ has type $A$ (i.e. if $c$ is odd), then $\la \l,\g_{k+1}\ra =0$.
\end{itemize}
Then $\mu^i=0$ for $1\le i\le k$, and also $\la \l,\g_{y}\ra =0$ for all $y$.
\end{lem}

\pf If $c$ is odd, then all the $C^i$ are of type $A$, and the conclusion follows from Theorems 14.2 and 16.1 of \cite{MF}.

Now suppose that $c$ is even and $c>0$, so $C^k$ is orthogonal. Then by the hypothesis (\ref{noteq1}), we have $L(\d)\ge 3$. By \cite[Thm. 14.2]{MF} we have $\mu^i=0$ for $0<i<k-1$. We indicate below summands of highest weights $\e$ and $\nu$ at level $k$ and $k-1$ respectively.  The first of these is given by Lemma \ref{setup} and the second is afforded by $f_{s+1}^{\frac{c-2}{2}}\bar f_s^{a_s} \cdots \bar f_1^{a_1}$:
\[
\begin{array}{l}
\e = (a_1,\cdots ,a_{s-1},a_s+\frac{c}{2},a_s+\frac{c}{2},a_{s-1},\cdots,a_1), \\
\nu = (a_1,\cdots ,a_{s-1},a_s+\frac{c-2}{2},a_s+\frac{c+2}{2},a_{s-1},\cdots,a_1).
\end{array}
\]
Since $L(\nu)\ge 4$, the first paragraph of the proof of \cite[Lemma 14.1.5]{MF} applies to show that $\mu^{k-1}$ is one of $0,\,\l_1^{k-1}$ or $\l_{r_{k-1}}^{k-1}$. Also $S(V^1) = S(\mu^0\downarrow L_X')+S(\nu) = S(\mu^0\downarrow L_X')  +2a_1+\cdots +2a_s+c$.    The argument of \cite[Thm. 14.2]{MF} rules out the case where $\mu^{k-1}=\l_1^{k-1}$. Now suppose 
$\mu^{k-1} = \l_{r_{k-1}}^{k-1}$. Then $V^2 \supseteq (\mu^0\downarrow L_X') \otimes \wedge^2\nu^* \otimes \e$. Now
\[
\wedge^2\nu^* \supseteq (2a_1,\cdots ,2a_{s-2}, 2a_{s-1}+1,2a_s+c, 2a_s+c-1,2a_{s-1},\cdots ,2a_1), 
\]
so using \cite[7.1.7]{MF}, we have 
\[
\wedge^2\nu^* \otimes \e \supseteq (3a_1,\cdots ,3a_{s-2}, 3a_{s-1}+2,3a_s+
\frac{3c-2}{2},3a_s+\frac{3c-4}{2},3a_{s-1}+1,3a_{s-2},\cdots ,3a_1)^2.
\]
It follows that $V^2$ has  a repeated summand of $S$-value at least $S(\mu^0\downarrow L_X') + 6\sum_1^sa_i+3c$, and this is greater than  $S(V^1)+1$, a contradiction. Therefore $\mu^{k-1}=0$. 

Now consider $\la \l,\g_y\ra$. Section 16 of \cite{MF} analyzes the terms $\g_y$ between adjacent Levi factors of type $A$. Since $C^i$ is of type $A$ for $0\le i\le k-1$, the arguments for \cite[Thm. 16.1]{MF} give $\la \l,\g_y\ra = 0$ for $y\le k-1$. Now suppose $\la \l,\g_k\ra \ne 0$. Then (again using \cite[7.1.7]{MF}, this time applied to $\nu^*\otimes \epsilon$), 
we see that $\l-\g_k$ affords a summand of $V^2(Q_Y)$ whose restriction $V^2$ to $L_X'$ contains a repeated summand of highest weight $(\mu^0\downarrow L_X') + \chi$, where 
\[
\chi = \left\{\begin{array}{l} (2a_1,\cdots ,2a_{s-2},2a_{s-1}+1,2a_s+c,2a_s+c-2,2a_{s-1}+1,2a_{s-2},\cdots ,2a_1), \hbox{ if }a_s>0 \\
(2a_1,\cdots ,2a_{j-2},2a_{j-1}+1,2a_j-1,0,\cdots ,0,c+1,c-1,0,\cdots ,0,2a_j-1,2a_{j-1}+1,\\
         2a_{j-2},\cdots ,2a_1), \hbox{ if }a_s=0
\end{array}
\right.
\]
(where, as in Section \ref{secevenrk}, $j$ is maximal such that $a_j\ne 0$). Since $S(V^1) = S(\mu^0\downarrow L_X')$, it follows that  $V^2$ has a repeated summand of $S$-value greater than $S(V^1)+1$, a contradiction. Hence $\la \l,\g_k\ra = 0$. 

It remains to handle the case where $c=0$. Here the argument is just the same as above, replacing $\e,\nu$ and $\chi$ by
\[
\begin{array}{l}
\e = (a_1,\cdots ,a_s,a_s,\cdots,a_1), \\
\nu = (a_1,\cdots ,a_j,0,\cdots,0,a_j-1, a_{j-1}+1,a_{j-2},\cdots,a_1), \\
\chi = (2a_1,\cdots ,2a_{j-2}+1,2a_{j-1}+1,2a_j-1,0,\cdots ,0,2a_j-1,2a_{j-1}+1,\cdots ,2a_1).
\end{array}
\]
This completes the proof. \hal

\begin{lem}\label{ceven1} If $c$ is even, then $\mu^k=0$.
\end{lem}

\pf Assume that $c$ is even and $\mu^k\ne 0$. By Lemma \ref{setup},  level $k$ restricted to $L_X'$ contains a summand of highest weight $\e = (a_1,\cdots ,a_{s-1},a_s+\frac{c}{2},a_s+\frac{c}{2},a_{s-1},\cdots,a_1)$. The image of $L_X'$ in the orthogonal group $C^k$ embeds in a product $R\times S$ of orthogonal subgroups such that the natural module for $S$ restricts to $L_X'$ as $\e$. 
Hence the Inductive Hypothesis, together with the argument of Lemma \ref{possmu0mu1} shows that either 
$\mu^k = \l_1^k$, or $s=1$, $\mu^k=\lambda_2^k$ and $(a_s,c)=(1,2)$. In the latter case, $k = 2$ and
$L_X'$ acts on the natural module for $C^2$ as $22+30+03+11$.  But  
$\wedge^2(22+30+03+11) \supseteq (33)^2$, contradicting the fact that $V^1$ is MF. Hence $\mu^k = \l_1^k$.

Now combining Lemmas 14.2.1 and 3.9 of \cite{MF}, we see that $V^2$ has a repeated summand of $S$-value at least $2S(\e)-1+\sum_0^{k-1}S(\mu^i\downarrow L_X')$. By Lemma \ref{prop38} this must be less than $S(V^1)+2$, and hence 
$2S(\e)-1 < S(\e)+2$, forcing $S(\e)\le 2$, a contradiction. \hal

\begin{lem}\label{thceven} Theorem $\ref{oddrk}$ holds if $c$ is even.
\end{lem}

\pf Assume $c$ is even. By the previous lemmas we have $\mu^i=0$ for $1\le i\le k$, and $\la \l,\g_j \ra = 0$ for all $j$. Hence $\l = \mu^0$, and so $\mu^0 \ne \l_1^0$. Now $V_{C^0}(\l_1^0)\downarrow L_X' \supseteq (a_1,\cdots,a_s,c,a_s,\cdots,a_2)$.  Also $C^0$ is of type $A$ and $V^1$ is MF, so \cite{MF} and (\ref{noteq1}) imply that either $\mu^0$ or its dual is one of $\l_2^0$, $2\l_1^0$, and $\d = (1,0,\cdots,0,c,0,\cdots,0,1)$, or $\mu^0$ is the dual of $\l_1^0$. Section 15 of \cite{MF} rules out the duals. Hence $\l = \l_2$ or $2\l_1$. But then $V_Y(\l) \downarrow X$ is not MF by \cite[Lemma 7.1.8]{MF}, a contradiction. \hal

At this point we may assume that $c$ is odd.

\begin{lem}\label{codd3} Suppose $c$ is odd, and either $c\ge 3$, or $a_s+\frac{c-1}{2} \ge 2$. Then $\mu^k=0$ and $\la \l,\g_{k+1}\ra = 0$.
\end{lem}

\pf Assume $\mu^k\ne 0$. Recall that $C^k$ is of type $A$, and also that level $k$ restricted to $L_X'$ (that is, $W^{k+1} = V_{C^k}(\l_1^k)\downarrow L_X'$) has a summand of highest weight 
\[
\nu = (a_1, \cdots a_{s-1},(a_s +\frac{c-1}{2}),(a_s +\frac{c+1}{2}),a_{s-1}, \cdots ,a_1)
\]
and this has maximal $S$-value (see Lemma \ref{setup}(ii)). Since by hypothesis, either $c\ge 3$, or $a_s+\frac{c-1}{2} \ge 2$, it follows from \cite[14.1.1]{MF} that $\mu^k = \l_1^k$ or $\l_{r_k}^k$.

If $\mu^k = \l_1^k$, then \cite[Lemmas 3.9,  14.2.1]{MF} implies that $V^2$ has a repeated summand of $S$-value at least $\sum_0^{k-1}S(\mu^i\downarrow L_X')+2S(\nu)-1$. Hence Lemma \ref{prop38} gives $2S(\nu)-1 \le S(\nu)+1$, a contradiction.

Now suppose $\mu^k = \l_{r_k}^k$. Here we have $V^2 \supseteq (\mu^0\otimes \cdots \otimes \mu^{k-1})\downarrow L_X' \otimes Z$, where $Z = \wedge^3\nu^*$ if $Y$ is orthogonal, and $Z = (\wedge^2\nu^* \otimes \nu^*)-\wedge^3\nu^*$ if $Y$ is symplectic. It now follows from Lemma \ref{extab} that $V^2$ has a repeated summand of $S$-value greater than $S(V^1)+1$, contradicting Lemma \ref{prop38}.

Hence $\mu^k=0$. Now suppose $\la \l,\g_{k+1}\ra >0$. Then considering the weight $\l-\g_{k+1}$, we see that $V^2$ 
has a summand $(\mu^0\otimes \cdots \otimes \mu^{k-1})\downarrow L_X' \otimes Z'$, where $Z' = \wedge^2\nu^*$ if $Y$ is orthogonal, and $Z' = S^2\nu^*$ if $Y$ is symplectic. This gives a contradiction as above, using 
\cite[Lemma 14.2.1]{MF}. Hence $\la \l,\g_{k+1}\ra =0$. \hal 

\begin{lem}\label{thcodd3} Theorem $\ref{oddrk}$ holds if $c$ is odd and either $c\ge 3$, or $a_s+\frac{c-1}{2} \ge 2$.
\end{lem}

\pf Suppose $c$ is odd and either $c\ge 3$, or $a_s+\frac{c-1}{2} \ge 2$. Then by Lemmas \ref{muk0lem} and \ref{codd3}, we have $\l = \mu^0$. Now we complete the proof exactly as in Lemma \ref{thceven}. \hal

In view of Lemma \ref{thcodd3}, we may now assume that $c=1$ and also that $a_s\le 1$.

\begin{lem}\label{d010} Suppose that $c=1$ and $a_s=0$. Then  $\mu^k=0$ and $\la \l,\g_{k+1}\ra = 0$.
\end{lem}

\pf Suppose $\mu^k\ne 0$. By Lemma \ref{setup}(ii), level $k$ restricted to $L_X'$ has a summand of highest weight $\nu = (a_1,\cdots ,a_{s-1},0,1,a_{s-1},\cdots ,a_1)$ of maximal $S$-value. As in the proof of Lemma \ref{codd3}, this forces 
$\mu^k = \l_1^k$ or $\l_{r_k}^k$, and the first possibility is excluded using \cite[Lemma 14.2.1]{MF}. The second possibility is also excluded as in the proof of Lemma \ref{codd3}, with the adjustment that instead of Lemma \ref{extab} we use Lemma \ref{extbd}, applied to $(a_j,0,\cdots ,0,1,a_{s-1})^*$ and then extended to higher rank, where $j$ is maximal such that $a_j>0$.

Hence $\mu^k=0$. We argue that $\la \l,\g_{k+1}\ra = 0$ exactly as in the proof of Lemma \ref{codd3}. \hal

\begin{lem}\label{d010} Suppose that $c=1$ and $a_s=1$. Then  $\mu^k=0$ and $\la \l,\g_{k+1}\ra = 0$.
\end{lem}

\pf Assume $\mu^k \ne 0$. Here level $k$ has a summand $\nu = (a_1,\cdots,a_{s-1},1,2,a_{s-1},\cdots,a_1)$. If $a_i\ne 0$ for some $i\le s-1$, then \cite{MF} implies that $\mu^k = \l_1^k$ or $\l_{r_k}^k$. Otherwise, $\mu^k$ could also be $\l_2^k$, $2\l_1^k$ or the dual of one of these; however in these cases \cite[Lemma 14.2.1]{MF} shows that $V^1$ is not MF, a contradiction. Hence $\mu^k = \l_1^k$ or $\l_{r_k}^k$.

If $\mu^k = \l_1^k$, we can use \cite[14.2.1]{MF} in the usual way to get a repeated summand of $V^2$ of $S$-value greater than $S(V^1)+1$. And if $\mu^k = \l_{r_k}^k$, then use of Lemma \ref{extab} gives the same conclusion.

Thus $\mu^k=0$. Finally, the same argument as in previous proofs gives $\la \l,\g_{k+1}\ra = 0$. \hal

\begin{lem}\label{thcodd} Theorem $\ref{oddrk}$ holds if $c=1$.
\end{lem}

\pf The previous two lemmas and Lemma \ref{muk0lem} show that $\l = \mu^0$. As $V_{C^0}(\l_1^0)\downarrow L_X'$ has a summand $\d' = (a_1,\cdots ,a_s,1,a_s,\cdots ,a_2)$, application of \cite{MF} implies that $\mu^0$ or its dual is one of 
$\l_1^0$, $\l_2^0$, $2\l_1^0$, $\l_3^0$ and $3\l_1^0$ (the last two only for $a_s = s=1$). The duals are ruled out by Section 15 of \cite{MF}. Hence $\l =  \l_2$, $2\l_1$, $\l_3$ or $3\l_1$ (the last two with $X=A_3$ and $\d = 111$). Now \cite{MF} (and Magma in the last case) shows that $V_Y(\l)\downarrow X$ is not MF, a contradiction. \hal

This completes the proof of Theorem \ref{oddrk}, and hence also of Theorem \ref{mainthm}.


  \end{document}